\documentclass{amsart}

\usepackage{style}

\title{Geometric Construction of Quiver Tensor Products}
\author{Daigo Ito and John S.\ Nolan}
\date{\today}

\begin{document}

\begin{abstract}
    By a classic theorem of Beilinson, the perfect derived category $\operatorname{Perf}(\mathbb{P}^n)$ of projective space is equivalent to the category of derived representations of a certain quiver with relations.
    The vertex-wise tensor product of quiver representations corresponds to a symmetric monoidal structure $\otimes_{\mathsf{Q}}$ on $\operatorname{Perf}(\mathbb{P}^n)$.
    We prove that, for a certain choice of equivalence, the symmetric monoidal structure $\otimes_{\mathsf{Q}}$ may be described geometrically as an \emph{extended convolution product} in the sense that the Fourier--Mukai kernel is given by the closure of the torus multiplication map in $(\mathbb{P}^n)^3$.
    We also set up a general framework for such problems, allowing us to generalize the extended convolution description of quiver tensor products to the case where $\mathbb{P}^n$ is replaced by any smooth complete toric variety of Bondal--Ruan type.
    Under toric mirror symmetry, this extended convolution product corresponds to the tensor product of constructible sheaves on a real torus.
    As another generalization of our results for $\mathbb{P}^n$, we show that any finite-dimensional algebra $A$ gives rise to a monoidal structure $\star_A'$ on $\operatorname{Perf}(\mathbb{P}(A))$, providing insights into the moduli of monoidal structures on $\operatorname{Perf}(\mathbb{P}^n)$.
\end{abstract}

\maketitle
\tableofcontents

\section{Introduction}
By \cite{beilinson1978coherent}, the (dg-enriched) derived category $\Perf(\PP^n)$ of perfect complexes on projective space is equivalent to the dg-category of $k$-linear functors $\Qsf_n\op \to \Perf(k)$, where $\Qsf_n$ is the category depicted in \cref{fig:beilinson_quiver} and $k$ is our ground field (assumed algebraically closed of characteristic zero).
We think of $\Qsf_n$ as a quiver with relations, so $k$-linear functors $\Qsf_n\op \to \Perf(k)$ are nothing but (derived) representations of this quiver.

\begin{figure}[h!]
    \centering
    \begin{tikzcd}
        q_0 \ar[r, shift left=2, "x_0"] \ar[r, phantom, "\cdots" description] \ar[r, shift right=2, "x_n", swap] & q_1 \ar[r, shift left=2, "x_0"] \ar[r, phantom, "\cdots" description] \ar[r, shift right=2, "x_n", swap] & \dots \ar[r, shift left=2, "x_0"] \ar[r, phantom, "\cdots" description] \ar[r, shift right=2, "x_n", swap] & q_n
    \end{tikzcd}
    \caption{The Beilinson quiver for $\PP^n$ (with relations $x_i x_j = x_j x_i$ for all $i, j$).}
    \label{fig:beilinson_quiver}
\end{figure}

The usual (vertexwise) tensor product of quiver representations defines a symmetric monoidal structure $\otimes_\Qsf$ on the category $\Perf \PP^n$.\footnote{We will write $\otimes_\Qsf$ both for the quiver tensor product on $\Fun(\Qsf_n\op, \Perf(k))$ and for the corresponding symmetric monoidal structure on $\Perf \PP^n$ (as well as for the variants that arise later), trusting in context to make the notation clear.}
More generally, for any smooth complete toric variety $X$ with a \emph{full strong exceptional collection of line bundles} in $\Perf(X)$ (cf. \cite{king_conj} for definitions and some key examples), there is an analogous quiver description of $\Perf(X)$ and an analogous quiver tensor product $\otimes_\Qsf$ on $\Perf(X)$.
Our goal in this project is to obtain a geometric understanding of the quiver tensor products $\otimes_\Qsf$ and some of their many generalizations.
Along the way we will make contact with several interesting subjects and phenomena, including:
\begin{itemize}
    \item Windows in geometric invariant theory,
    \item Homological mirror symmetry for toric varieties,
    \item Tensor triangular geometry, 
    \item ``Categorical compactifications'' of linear algebraic groups, and
    \item ``Moduli theory'' of monoidal structures on a fixed stable $\infty$-category / pretriangulated dg-category.
\end{itemize}

We claim that the quiver tensor product $\otimes_\Qsf$ on $\Perf \PP^n$ can be described as an ``extended convolution'' (or ``EC'') product of sheaves.
More precisely, let
\[
    T = \bset{[x_0 : \dots : x_n]}{x_i \neq 0 \textrm{ for all $i$}} \cong \GG_m^n
\]
be the standard torus in $\PP^n$.
For a morphism $f: T^a \to T^b$, let $Z_f$ be the closure in $(\PP^n)^{a+b}$ of the graph of $f$.
In particular, if $\mu: T \times T \to T$ is the standard coordinatewise multiplication map, there is a correspondence
\[
\begin{tikzcd}
    & Z_\mu \ar[dl, "r_1 \times r_2", swap] \ar[dr, "r_3"] & \\
    \PP^n \times \PP^n & & \PP^n
\end{tikzcd}
\]
Using Beilinson's resolution of the diagonal, it is not too hard to compute that
\[
    \Fsc \otimes_\Qsf \Gsc = r_{3*}(r_1^* \Fsc \otimes r_2^* \Gsc)
\]
for all $\Fsc, \Gsc \in \Perf(\PP^n)$.\footnote{
We use ``implicitly derived'' notation for derived functors and for derived categories of sheaves, e.g.\ $\QC$ is used for (dg-enriched) derived categories of quasicoherent sheaves.
We will often write $=$ for a natural / preferred equivalence.}
In particular, if $\Fsc$ and $\Gsc$ are skyscraper sheaves supported on $T$, then $\Fsc \otimes_\Qsf \Gsc = \mu_*(\Fsc \boxtimes \Gsc)$, the convolution product of $\Fsc$ and $\Gsc$, justifying the name ``extended convolution.''

The above description of EC products has one troublesome feature: \emph{it is not clear \emph{a priori} that the push-pull operation $(\Fsc, \Gsc) \mapsto r_{3*}(r_1^* \Fsc \otimes r_2^* \Gsc)$ is associative!}
This traces back to the non-functoriality of graph closures: for $f: T^a \to T^b$ and $g: T^b \to T^c$, we typically have $Z_{g \circ f} \neq Z_g \times_{\PP^b} Z_f$, even when we try to fix issues of non-flatness by taking derived fiber products.
In particular, as $Z_{\mu \times \id}  \times_{\PP^2} Z_\mu \neq Z_{\id \times \mu}  \times_{\PP^2} Z_\mu$, it is not clear that the two-fold multiplication map given by pulling back and pushing forward along the topmost path of the diagram  
\[
\begin{tikzcd}
    & & Z_{\mu \times \id}  \times_{(\PP^n)^2} Z_\mu \ar[dl, swap] \ar[dr] & & \\
    & Z_{\mu \times \id} \ar[dl, swap] \ar[dr] & & Z_\mu \ar[dl, swap] \ar[dr] & \\
    (\PP^n)^3 & & (\PP^n)^2 & & \PP^n \\
\end{tikzcd}
\]
agrees with the push-pull operation for the analogous diagram with $\mu \times \id$ replaced by $\id \times \mu$.
The operations here do in fact agree, essentially because the pushforwards of the corresponding structure sheaves to $(\PP^n)^3 \times \PP^n$ are naturally isomorphic.

It is not immediately obvious how to generalize the construction of EC products to analogous examples or how to check that the EC product upgrades to a symmetric monoidal structure at the $\infty$-categorical level.
Thus we would like a general construction of EC products for which the associativity and higher coherence data is ``obvious.''

One of our main results is that such a construction is possible for smooth complete toric varieties of \emph{Bondal-Ruan type} (see \cref{dfn:bondal_thomsen}), which include $\PP^n$ and many toric Fano varieties.
For such a variety $X$, the \emph{Bondal-Thomsen collection} $\Theta$ gives a full strong exceptional collection of line bundles on $X$ and thus an equivalence $\Perf(X) \simeq \Fun\big(\Qsf_{\Theta}, \Perf(k)\big)$ for some quiver with relations $\Qsf_\Theta$.
By \cite{bondal2006derived}, this equivalence may also be understood through homological mirror symmetry as follows.

Let $M$ be the cocharacter lattice of the dense torus in $X$.
Then $M$ is a free abelian group of rank $\dim X$, $M_\RR = M \otimes_\ZZ \RR$ is a real vector space of dimension $\dim X$, and $M_\RR / M$ is a real torus of dimension $\dim X$.
Then there is a natural equivalence
\begin{equation} \label{eq:bondal_ruan_hms}
    \Perf(X) \simeq \Fun\big(\Qsf_{\Theta}, \Perf(k)\big) \simeq \Sh_Z^\perf(M_\RR / M)
\end{equation}
where $\Sh_Z^\perf(M_\RR / M)$ denotes the derived category of constructible sheaves of $k$-vector spaces on the real torus $M_\RR / M$ which:
\begin{itemize}
    \item are constructible with respect to a certain real stratification $Z$ determined by the fan of $X$, and
    \item have perfect stalks.
\end{itemize}
In particular, because we are working with sheaves which are constructible with respect to a fixed stratification, $\Sh_Z^\perf(M_\RR / M)$ is closed under stalk-wise tensor product $\otimes_k$ of constructible sheaves.
We may upgrade \cref{eq:bondal_ruan_hms} to a symmetric monoidal equivalence:

\begin{thm}[{\cref{mainthm_toric}, \cref{prop:hms_tensor}}]
    Let $X$ be a smooth complete toric variety of Bondal-Ruan type.
    Write $Z_\mu$ for the closure in $X^3$ of the graph of the binary multiplication on the dense torus in $X$.
    Then push-pull along the correspondence
    \[
        \begin{tikzcd}
              & Z_\mu \ar[dl] \ar[dr] & \\
            X \times X & & X
        \end{tikzcd}
    \]
    extends to a symmetric monoidal structure $\star'_X$ on $\Perf(X)$.
    There are symmetric monoidal equivalences
    \[
        \big(\Perf(X), \star'_X\big) \simeq \Big(\Fun\big(\Qsf_{\Theta}, \Perf(k)\big), \otimes_\Qsf\Big) \simeq \big(\Sh_{\Lambda}^\perf(M_\RR / M), \otimes_k\big).
    \]
\end{thm}

\begin{rmk}
    We caution readers that our conventions for toric mirror symmetry here are dual to those commonly used in the literature.
    See \cref{rmk:dual_conventions} and \cref{sec:hms} for a more careful discussion.
\end{rmk}

For smooth complete toric varieties $X$ which admit full strong exceptional collections of line bundles but which are not of Bondal-Ruan type, we are still able to construct EC products as Fourier-Mukai transforms defined in terms of resolutions of the diagonal on $X$ (\cref{prop:ec_fm_computation}).
These are still equivalent to the quiver tensor products $\otimes_\Qsf$.
However, in this generality, we are unable to give a simple ``push-pull'' description of the EC products -- the Fourier-Mukai kernel on $X^3$ does not have an algebra structure \emph{a priori}, so we cannot expect it to be a pushforward of the structure sheaf of some scheme (or stack).
We suspect that such an algebra structure does not exist without additional hypotheses on the full strong exceptional collection.

We may use similar methods to produce monoidal structures on the category $\Perf \PP(A)$ whenever $A$ is a finite-dimensional $k$-algebra.

\begin{thm}[\cref{mainthm_algebra}, \cref{prop:algebra_reconstruction}] \label{mainthm_algebra_introduction}
    Let $A$ be a nonzero finite-dimensional $k$-algebra.
    Let $\AA(A)^\times$ be the subgroup of units in $\AA(A)$.
    Write $Z_A$ for the closure in $\PP(A)^3$ of the graph of the binary multiplication on the group scheme $[\AA(A)^\times / \GG_m]$.
    Then:
    \begin{enumerate}
        \item Push-pull along the correspondence
        \[
            \begin{tikzcd}
                & Z_A \ar[dl] \ar[dr] & \\
                \PP(A) \times \PP(A) & & \PP(A)
            \end{tikzcd}
        \]
        defines a monoidal structure $\star'_A$ on $\Perf (\PP(A))$.
        \item If $j': [\AA(A)^\times / \GG_m] \hookrightarrow \PP(A)$ is the inclusion, then the pushforward functor 
        \[
            j'_*: \big(\QC([\AA(A)^\times / \GG_m]), \star_{[\AA(A)^\times / \GG_m]}\big) \to \big(\QC(\PP(A)), \star'_A\big)
        \]
        is monoidal.
        \item The construction of $\star'_A$ is functorial in surjections of finite-dimensional $k$-algebras.
        \item The construction of $\star'_A$ is essentially injective in the sense that if $(\Perf (\PP(A')), \star_{A'}')$ is monoidally equivalent to $(\Perf (\PP(A)), \star_{A}')$ for a finite-dimensional $k$-algebra $A$, then $A'$ is isomorphic to $A$ as $k$-algebras.
    \end{enumerate}
    When $A$ is commutative, we may replace ``monoidal'' with ``symmetric monoidal'' throughout.
\end{thm}

\cref{mainthm_algebra_introduction} may be understood as giving a map from the moduli space of $(n+1)$-dimensional (commutative) $k$-algebras (cf.\ \cite{poonen2008moduli}) to the ``moduli stack of (symmetric) monoidal structures on $\Perf(\PP^n)$.''
This map is injective on geometric points, and we expect that it parametrizes a component of the latter moduli stack. Moreover, this component corresponds to monoidal structures with zero-dimensional Balmer spectrum (\cref{prop:balmer_spectrum}), whereas the sheaf tensor product sits in a component with $n$-dimensional Balmer spectrum.
By a version of the Bondal--Orlov reconstruction theorem (\cite[Corollary 1.4]{toledo2024tensor}), we expect that the component of the moduli stack of symmetric monoidal structures on $\Perf(\PP^n)$ containing the usual tensor product contains a single geometric point, though we are not aware of a complete proof of this claim.

One may also use \cref{mainthm_algebra_introduction} to obtain ``categorical compactifications'' of many familiar groups and algebras -- see \cref{sec:Pn} for details.
It would be interesting to understand which groups admit categorical compactifications in general.

\subsection{Approach}

Let us explain our method of constructing the EC product on $\Perf \PP^n$ in a way which makes the associativity ``obvious.'' 
The strategy is similar for other examples, and in the body of the paper we introduce the notion of \emph{(geometric) EC setup} (\cref{dfn:ec}) to handle all of our examples at once.

Note that coordinatewise multiplication does not define a commutative monoid structure on $\PP^n$, so we cannot view our EC product as a genuine convolution product.
To fix this, we first consider the quotient stack $[\AA^{n+1} / \GG_m]$, which contains $\PP^n$ as an open substack and admits a well-defined commutative monoid structure.

\begin{rmk}
    Recall that the derived category of quasicoherent sheaves $\QC([\AA^{n+1} / \GG_m])$ is equivalent to the derived category of $\GG_m$-equivariant quasicoherent sheaves on $\AA^{n+1}$.
    Although we use the language of stacks for convenience, readers more comfortable with equivariant geometry may use the latter language without losing much.
\end{rmk}

The (suitably $\GG_m$-equivariant) coordinatewise product on $\AA^n$ descends to give a commutative monoid structure on the quotient stack $[\AA^{n+1} / \GG_m]$ and thus a symmetric monoidal convolution product $\star_{[\AA^{n+1}/\GG_m]}$ on $\QC([\AA^{n+1} / \GG_m])$.
The category $\QC([\AA^{n+1} / \GG_m])$ is equivalent to the category of functors from $\Qsf_{n,\infty}\op$ to $\Dsf(k)$, where $\Qsf_{n,\infty}$ is the ``infinite Beilinson quiver'' depicted in \cref{fig:infinite_beilinson_quiver}.
We show (as a special case of \cref{thm:basic_comparison}) that this equivalence upgrades to a symmetric monoidal equivalence
\[
    \big(\QC([\AA^{n+1} / \GG_m]), \star_{[\AA^{n+1}/\GG_m]}\big) \simeq \big(\Fun(\Qsf_{n,\infty}\op, \Dsf(k)), \otimes_\Qsf\big).
\]
That is, we may identify the convolution product on $[\AA^{n+1} / \GG_m]$ with the quiver tensor product on $\Fun(\Qsf_{n,\infty}\op, \Dsf(k))$.
This extends the main theorem of \cite{moulinos_geometry}.

\begin{figure}[h!]
    \centering
    \begin{tikzcd}
        \dots \ar[r, shift left=2, "x_0"] \ar[r, phantom, "\cdots" description] \ar[r, shift right=2, "x_n", swap] & q_{-1} \ar[r, shift left=2, "x_0"] \ar[r, phantom, "\cdots" description] \ar[r, shift right=2, "x_n", swap] & q_0 \ar[r, shift left=2, "x_0"] \ar[r, phantom, "\cdots" description] \ar[r, shift right=2, "x_n", swap] & q_1 \ar[r, shift left=2, "x_0"] \ar[r, phantom, "\cdots" description] \ar[r, shift right=2, "x_n", swap] & \dots
    \end{tikzcd}
    \caption{The infinite Beilinson quiver for $[\AA^{n+1} / \GG_m]$ (with relations $x_i x_j = x_j x_i$ for all $i, j$).
    Vertices are indexed by $\weight(\GG_m) = \ZZ$.}
    \label{fig:infinite_beilinson_quiver}
\end{figure}

Now let us return to the case of $\PP^n$.
There is an embedding of the usual Beilinson quiver $\Qsf_n$ into the infinite Beilinson quiver $\Qsf_{n,\infty}$ (given by $q_i \mapsto q_i$).
Left Kan extension along this embedding, i.e.\ universally filling in the diagram
\[
\begin{tikzcd}
    \Qsf_n \rar \dar[hook] & \Dsf(k), \\
    \Qsf_{n,\infty} \ar[ur, dashed]
\end{tikzcd}
\]
gives a functor $\Fun(\Qsf_n\op, \Dsf(k)) \hookrightarrow \Fun(\Qsf_{n,\infty}\op, \Dsf(k))$.
Using our aforementioned equivalences, this functor may be rewritten as $W: \QC(\PP^n) \hookrightarrow \QC([\AA^{n+1} / \GG_m])$.
Here $W$ is an example of a \emph{window} in geometric invariant theory (see \cite{hl_derived}).

The right adjoint $H: \QC([\AA^{n+1} / \GG_m]) \to \QC(\PP^n)$ to $W$, which we call the \emph{Hitchcock functor},\footnote{Our name comes from the classic Alfred Hitchcock film \emph{Rear Window}: the Hitchcock functor $H$ ``sees $\QC([\AA^{n+1} / \GG_m])$ through the window $W$.''}
corresponds to pullback of functors / restriction of quiver representations along the embedding $\Qsf_n \hookrightarrow \Qsf_{n,\infty}$.
It follows that the quiver tensor product $\otimes_\Qsf$ on $\QC(\PP^n)$ is the unique symmetric monoidal structure  such that the composite
\[
    \big(\Fun(\Qsf_{n,\infty}\op, \Dsf(k)), \otimes_\Qsf\big)  \simeq \big(\QC([\AA^{n+1} / \GG_m]), \star_{[\AA^{n+1}/\GG_m]}\big) \xrightarrow{H} \big(\QC(\PP^n), \otimes_\Qsf\big)
\]
is symmetric monoidal (stated more generally as \cref{prop:ec_exists}).

In the case of $\PP^n$ (and other smooth complete toric varieties of Bondal-Ruan type), it turns out that the Hitchcock functor $H$ can be described entirely in terms of geometry.
More precisely, $H: \QC([\AA^{n+1} / \GG_m]) \to \QC(\PP^n)$ is given by push-pull along the correspondence
\[
    \begin{tikzcd}
        & \ol{\Delta}_{\PP^n} \ar[dr] \ar[dl] & \\
        {[\AA^{n+1} / \GG_m]} & & \PP^n.
    \end{tikzcd}
\]
where $\ol{\Delta}_{\PP^n} \subset [\AA^{n+1} / \GG_m] \times \PP^n$ is the closure of the diagonal of $\PP^n$.
Thus we may \emph{define} the EC product $\star'_{\PP^n}$ as the unique symmetric monoidal structure on $\QC(\PP^n)$ such that the Hitchcock functor
\[
H: \big(\QC([\AA^{n+1} / \GG_m]), \star_{[\AA^{n+1}/\GG_m]}\big) \to \big(\QC(\PP^n), \star'_{\PP^n}\big)
\]
is symmetric monoidal, so associativity of $\star'_{\PP^n}$ is immediate.
(Showing that $\star'_{\PP^n}$ is well-defined does take some effort, but we are able to reduce the verification to \HA{Proposition 2.2.1.9}.)

It is clear that:
\begin{itemize}
    \item There is a symmetric monoidal equivalence $\big(\QC(\PP^n), \star'_{\PP^n}\big) \simeq \big(\Fun(\Qsf_n\op, \Dsf(k)), \otimes_\Qsf\big)$.
    \item This restricts to a symmetric monoidal equivalence $\big(\Perf(\PP^n), \star'_{\PP^n}\big) \simeq \big(\Fun(\Qsf_n\op, \Perf(k)), \otimes_\Qsf\big)$.
    \item The geometry of $\star'_{\PP^n}$ may be understood entirely in terms of the geometry of $\PP^n$ and $[\AA^{n+1} / \GG_m]$.
\end{itemize}
We are also able to recover the aforementioned push-pull description of extended convolution as a consequence of this definition of $\star'_{\PP^n}$.

\begin{rmk}
    In the body of the paper we make use of $\infty$-categorical techniques and techniques from derived algebraic geometry due to the added flexibility these methods provide.
    By remembering extra data about the homotopical structure of morphisms, we are able to make cleaner statements and more conceptual arguments throughout.
\end{rmk}

\subsection{Outline}

In \cref{sec:qc_dag} we provide a brief review of concepts and techniques from derived algebraic geometry and $\infty$-category theory that we will use throughout the paper.
The material in this section is all known to the experts.
We include it to assist readers not well-versed in higher algebra and to include some useful statements we are not able to locate in the literature.

In \cref{sec:quiver} we begin our discussion of quiver tensor products.
We introduce \emph{transparent collections of weights}, a variation on the notion of ``full strong exceptional collection of line bundles'' that allows us to construct windows of a quiver-theoretic nature.
We also discuss the structure of these windows and their associated Hitchcock functors.

\cref{sec:basic} contains a proof of the equivalence of convolution products on certain $\EE_n$-monoid derived stacks and the corresponding ``quiver tensor products.''
This generalizes the above discussion of convolution on $[\AA^{n+1} / \GG_m]$.

In \cref{sec:ec}, we combine the results of \cref{sec:quiver} and \cref{sec:basic} to construct EC products and prove their equivalence with quiver tensor products.
We use toric varieties of Bondal-Ruan type and projectivizations of algebras $\PP(A)$ as running examples.

\cref{sec:Pn} contains some first consequences of the existence of EC products.
We introduce ``categorical compactifications'' of groups and algebras and give some examples, though the question of their existence in general remains open.
We also compute some classical invariants of EC products.

Finally, in \cref{sec:hms}, we discuss Bondal-Ruan mirror symmetry.
We prove that EC products are mirror to the tensor product of constructible sheaves for toric varieties of Bondal-Ruan type.

In \cref{app:monoidal_adjunctions}, we include some useful technical results on the relationship between monoidal structures and adjunctions in $\infty$-category theory.

\subsection{Acknowledgments}

We are deeply grateful to our PhD advisors, David Nadler (advising DI) and Constantin Teleman (advising JSN), for years of useful conversations and patient feedback.
We also benefited greatly from conversations with 
Colleen Delaney,
David Favero,
Martin Gallauer,
Swapnil Garg,
Daniel Halpern-Leistner,
Kimoi Kemboi,
Yuji Okitani, and
Ed Segal.
JSN was supported in part by the Simons Collaboration on Global Categorical Symmetries.

\subsection{Notation}
We follow the following convetions/notations throughout the paper. 
\begin{itemize}
    \item Our base field $k$ is assumed to be algebraically closed of characteristic $0$.
    (In some cases, this hypothesis can be dropped; we will indicate when it is useful or interesting to do so.)
    \item Algebras are assumed to be unital and associative but not necessarily assumed to be commutative unless otherwise specified. 
    \item We use ``implicitly derived'' notation for derived functors and for derived categories of sheaves, e.g.\ $\QC$ is used for (dg-enriched) derived categories of quasicoherent sheaves.
    Derived categories will always be assumed to be enriched (i.e.\ dg-categories or stable $\infty$-categories).
    \item The symbol $=$ will often be used to denote a natural / preferred choice of equivalence.
    \item If $\Qsf$ is a small $k$-linear category, we write $\Dsf(\Qsf\op) := \Fun(\Qsf, \Dsf(k))$ for the (enriched) derived category of right $\Qsf$-modules.
    \item Write $(-)^\simeq$ for the maximal subgroupoid functor $\Cat_\infty \to \Sc$, where $\Cat_\infty$ is the (large) $\infty$-category of small $\infty$-categories and $\Sc$ is the $\infty$-category of spaces (cf. \cref{sec:basics_infty_cats}).
    In particular, we view every $\Cat_\infty$-enriched category as an $\infty$-category by applying $(-)^\simeq$ on all mapping $\infty$-categories.
    \item Full subcategories are always assumed to be strictly full (i.e.\ closed under equivalence).
\end{itemize}

\section{Quasicoherent sheaves in derived algebraic geometry} \label{sec:qc_dag}

To construct extended convolution products at an $\infty$-categorical level, we will need to use some methods from derived (and spectral) algebraic geometry.
In the interest of keeping this paper comprehensible to readers from a more 1-categorical background, and in the interest of writing down things we have not seen explicitly stated in the literature, we provide an overview of some definitions, results, and methods we will use.
We recommend \cite{bzfn_integral} and \cite{scholze_6ff} for a more extended (but still elementary and conceptual) discussion of many of the ideas mentioned here.
Everything in this section is well-known to the experts.

\subsection{Presentability and compact generation}\label{sec:basics_infty_cats}

To make our claims about homotopy coherence rigorous, we use the language of $\infty$-categories as developed in \cite{htt} and \cite{ha} (among other references).
This subsection is devoted to a rapid review of the theory of \emph{presentable} $\infty$-categories discussed in \cite{htt}.

\begin{rmk}
    There is a (classical) analogous theory of presentable 1-categories, but it is less central to applications.
    Many of the consequences of the classical theory (e.g.\ existence of adjoint functors) are easy to check ``by hand'' as needed.
    By contrast, proving results ``by hand'' in the $\infty$-categorical context is often much more difficult, and we are left with no recourse but to use categorical methods.
\end{rmk}

The essential idea of $\infty$-category theory is to replace sets in ordinary category theory by \emph{spaces}, a.k.a.\ (weak) \emph{homotopy types}.
The category of spaces is denoted $\Sc$.
There are many equivalent point-set models one can use to understand and construct $\Sc$.
For example, one may obtain $\Sc$ from the category of CW-complexes and continuous maps by ``inverting homotopy equivalences'' in a suitable sense.

Write $\Cat_\infty$ for the (large) $\infty$-category of small $\infty$-categories and $\widehat{\Cat}_\infty$ for the (very large) $\infty$-category of large $\infty$-categories.\footnote{Standard techniques from the theory of \emph{Grothendieck universes} allow us to deal with most set-theoretic ``size issues'' when they arise.}
For any $\Cc_0 \in \Cat_\infty$, we define the $\infty$-category of \emph{presheaves on $\Cc_0$} as $\PSh(\Cc_0) = \Fun(\Cc\op, \Sc)$.
As usual, we have a Yoneda embedding $\Cc_0 \hookrightarrow \PSh(\Cc_0)$.

Presheaf categories are particular cases of presentable $\infty$-categories.
An $\infty$-category $\Cc$ is \emph{presentable} if there exists an $\infty$-category $\Cc_0$ and an functor $L: \PSh(\Cc_0) \to \Cc$ such that $L$ has a fully faithful right adjoint which preserves $\kappa$-filtered colimits for some regular cardinal $\kappa$.
Such a category $\Cc$ possesses many pleasing potential features of large $\infty$-categories -- in particular, $\Cc$ is complete and cocomplete.
However, the behavior of $\Cc$ is still controlled by that of the small $\infty$-category $\Cc_0$.

The adjoint functor theorem \HTT{Corollary 5.5.2.9} states that a functor $F: \Cc \to \Cc'$ between presentable $\infty$-categories has a right adjoint if and only if $F$ preserves colimits.
Motivated by this definition, we let $\Pr^L$ denote the $\infty$-category of presentable $\infty$-categories and colimit-preserving functors.
We write $\Fun^L(\Cc, \Cc')$ for the $\infty$-categories of functors in $\Pr^L$ from $\Cc$ to $\Cc'$.
Letting $\Pr^R$ be the $\infty$-category of presentable $\infty$-categories and functors which preserve limits and $\kappa$-filtered colimits (for some regular cardinal $\kappa$), the adjoint functor theorem upgrades to an equivalence $(\Pr^L)\op \xrightarrow{\sim} \Pr^R$ which is the identity on objects.

The category $\Pr^L$ has a natural symmetric monoidal structure, the \emph{Lurie tensor product} $\otimes$, defined so that colimit-preserving functors $\Cc_1 \otimes \Cc_2 \to \Cc_3$ are functors $\Cc_1 \times \Cc_2 \to \Cc_3$ which preserve colimits in each variable separately.
Our later variants of $\Pr^L$ will all inherit a corresponding Lurie tensor product.

If $\Cc \in \Pr^L$, we say that:
\begin{itemize}
    \item An object $c_0 \in \Cc_0$ is \emph{compact} if the functor $\Hom_\Cc(c_0, -)$ preserves filtered colimits.
    \item A small full subcategory $\Cc_0 \subset \Cc$ \emph{generates} $\Cc$ if a morphism $f: c \to c'$ is an isomorphism if and only if the induced map $\Hom_\Cc(c_0, c) \to \Hom_\Cc(c_0, c')$ is an isomorphism for all $c_0 \in \Cc$.
    Equivalently, every object of $\Cc$ may be expressed as a colimit of objects of $\Cc_0$.
\end{itemize}
If the objects of $\Cc_0$ are compact and $\Cc_0$ generates $\Cc$, we say that $\Cc$ is \emph{compactly generated} by $\Cc_0$.
In this case, there is a natural equivalence $\Cc = \Indsf(\Cc_0)$, where the functor $\Ind: \Cat_\infty \to \Pr^L$ ``freely adjoins filtered colimits'' to its input (see \HTT{\S 5.3.5} for a precise definition).
In particular, every object of $\Cc$ can be obtained (canonically) as a colimit of objects of $\Cc_0$.
We let $\Pr^L_{\omega}$ denote the $\infty$-category of compactly generated $\infty$-categories and functors which preserve colimits and compact objects.

Let $(\Vc, \otimes)$ be a \emph{presentably symmetric monoidal $\infty$-category}, i.e.\ an object of $\CAlg(\Pr^L, \otimes)$.
There exists a rich theory of ``$\Vc$-enriched $\infty$-categories'' -- see e.g.\ \cite[Appendix A]{mgs_universal} for a highly readable account of the subject.
The above results extend (with some mild modifications) to the $\Vc$-enriched setting.
In fact, by \cite[Theorem A.3.8]{mgs_universal}, a $\Vc$-enriched presentable $\infty$-category is the same as a module over $(\Vc, \otimes)$ in $\Pr^L$.
In particular, we note that $\Vc$ is automatically enriched over itself, and that any small $\Vc$-enriched category $\Cc$ admits a Yoneda embedding $\Cc \hookrightarrow \Fun_\Vc(\Cc\op, \Vc)$.
We may use this result to reduce the study of presentable $\infty$-categories enriched in $(\Vc, \otimes)$ to the study of unenriched presentable $\infty$-categories.

Two key examples of the above are as follows:
\begin{itemize}
    \item When $(\Vc, \otimes) = (\Dsf(k), \otimes_k)$, where $\Dsf(k)$-enriched $\infty$-categories are the same as \emph{$k$-linear $\infty$-categories}.
    We write $\Cat_k$ for the large category of small $k$-linear $\infty$-categories.
    \item We may also take $(\Vc, \otimes) = (\Sp, \otimes)$, the $\infty$-category of \emph{spectra}.
    Here $\otimes$ is the \emph{smash product}.
\end{itemize}

\subsection{$\EE_n$-algebras}

The $\infty$-categorical theory of (commutative) algebras in a (symmetric) monoidal category is subsumed by the theory of \emph{$\EE_n$-algebras} (for $1 \leq n \leq \infty$), or more generally \emph{algebras over an $\infty$-operad}, as developed in \cite{ha}.
We'll focus on the $\EE_n$-case here for simplicity.
In fact, in this work, we are only truly interested in the case $n = 1$ or $n = \infty$.
We use the terminology of $\EE_n$-algebras primarily as an efficient method of covering both commutative and noncommutative cases with the same statement.

An $\EE_n$-algebra in a symmetric monoidal $\infty$-category can be thought of as an algebra object with $n$ compatible associative multiplications.
The existence of these multiplications enforces a sort of commutativity on the operations involved.
In particular:
\begin{itemize}
    \item $\EE_1$-algebras are the same as associative algebras.
    \item $\EE_\infty$-algebras are the same as commutative algebras.
\end{itemize}
Note a key difference with the classical case: commutativity is no longer a \emph{property} of the multiplication but an extra \emph{structure} witnessed by the infinitely many compatible multiplications.

\begin{rmk}\label{rmk:en-algebra-in-1-cat}
    When working in a $1$-category, $\EE_n$-algebras for $n \geq 2$ are always commutative (i.e.\ the same as $\EE_\infty$-algebras).
    The differences between $\EE_2$-algebras and $\EE_\infty$-algebras appear only when we allow nontrivial $2$-morphisms.
\end{rmk}

\begin{rmk}\label{rmk:en-algebra-in-2-cat}
    When working in $\Cat_\infty$ (or similar $\infty$-categories such as $\Cat$, the 2-category of discrete categories):
    \begin{itemize}
        \item $\EE_1$-algebras are the same as monoidal categories.
        \item $\EE_2$-algebras are the same as braided monoidal categories.
        \item $\EE_\infty$-algebras are the same as symmetric monoidal categories.
    \end{itemize}
    In $\Cat_\infty$, we may take this as a definition of ``monoidal / braided monoidal category'' (though we need to bootstrap the definitions so that we may consider $\Cat_\infty$ as a symmetric monoidal category).
\end{rmk}

We write $\Alg_{\EE_n}(\Cc) = \Alg_{\EE_n}(\Cc, \otimes)$ for the $\infty$-category of $\EE_n$-algebras in a symmetric monoidal $\infty$-category $(\Cc, \otimes)$.\footnote{This definition actually makes sense when $(\Cc, \otimes)$ is only assumed to be $\EE_n$-monoidal itself.}
Lax $\EE_n$-monoidal functors induce functors between the corresponding $\infty$-categories of algebras (this is automatic from the definition, \HA{2.1.3.1}).
There is also an analogous statement for $\EE_n$-monoidal adjunctions (\HA{7.3.2.13}, restated as \cref{prop:monoidal_adjunction_algebra} here and strengthened in \cref{prop:partial_monoidal_adjunction_algebra}).

\subsection{Stable $\infty$-categories}

Recall that an $\infty$-category $\Cc$ is \emph{stable} if (\HA{Proposition 1.1.3.4}):
\begin{itemize}
    \item $\Cc$ admits finite limits and finite colimits, and
    \item A commutative square
    \[
        \begin{tikzcd}
            A \rar \dar & B \dar \\
            C \rar & D
        \end{tikzcd}
    \]
    is a pullback square if and only if it is a pushout square.
\end{itemize}
These conditions imply that $\Cc$ has a zero object $0$.
If the diagram
\[
    \begin{tikzcd}
            A \rar["f"] \dar & B \dar["g"] \\
            0 \rar & C
    \end{tikzcd}
\]
is a pullback square, we say that $A$ is the \emph{fiber} of $g$ (denoted $\fib(g)$) and that $C$ is the \emph{cofiber} of $f$ (denoted $\cofib(g)$).
There is a natural \emph{shift autoequivalence} $[1]: \Cc \to \Cc$ defined by $A[1] = \cofib(A \to 0)$.\footnote{When working in the $\infty$-categorical context, such colimits are typically nontrivial.}
Every stable $\infty$-category is canonically enriched over $(\Sp, \otimes)$.

The $\infty$-category of pretriangulated dg-categories over $k$ is equivalent to the $\infty$-category of stable $\infty$-categories enriched over $\Dsf(k)$.
Via this equivalence, cones in a pretriangulated dg-category correspond to cofibers in the corresponding stable $\infty$-category.
In other words, the \emph{homological algebra} of cones, shifts, etc.\ in dg-categories translates into the \emph{homotopical algebra} of (co)limits in stable $\infty$-categories.
The stable $\infty$-categories of primary interest to us are enriched over $\Dsf(k)$, so readers will not lose much by thinking of these stable $\infty$-categories as pretriangulated dg-categories.

\begin{rmk}
    The \emph{homotopy category} of a stable $\infty$-category / pretriangulated dg-category is a triangulated category.
    We work in the context of stable $\infty$-categories to avoid various difficulties in the theory of triangulated categories, e.g.\ non-existence of tensor products of categories, non-functoriality of cones, poor behavior of Fourier-Mukai transforms, \textrm{etc.}
\end{rmk}

Write $\Pr^L_\st$ for the $\infty$-category of stable, presentable $\infty$-categories and colimit-preserving functors.
Let $\Pr^L_k$ be the category of $\Dsf(k)$-modules in $\Pr^L_\st$.
To avoid repeating ourselves too much, we will state the following results for $\Pr^L_k$ (though the analogues for $\Pr^L_\st$ also hold).

The theory of compact generation admits some simplifications in the stable setting.
For $\Cc \in \Pr^L_k$, we have the following (straightforward) results:
\begin{itemize}
    \item An object $c \in \Cc$ is compact if $\Hom(\Cc, -)$ commutes with infinite direct sums (and hence with all colimits).
    \item A full subcategory $\Cc_0 \subset \Cc$ generates $\Cc$ if, whenever $c \in \Cc$ satisfies $\Hom(c_0, c) = 0$ for all $c_0 \in \Cc_0$, we must have $c = 0$.
\end{itemize}
If $\Cc$ is compactly generated by $\Cc_0$, there is a natural equivalence $\Cc = \Indsf(\Cc_0)$, where we use a $k$-linear version of $\Indsf$.
In fact, we have $\Indsf(\Cc_0) = \Fun_k(\Cc_0, \Dsf(k))$ in the stable $k$-linear setting.

Let $\Pr^L_{k,\omega}$ be the $\infty$-category of stable compactly generated $k$-linear $\infty$-categories with morphisms given by $k$-linear functors that preserve colimits and compact objects.
Letting $\Pr^R_{k,\omega}$ be the category with the same objects but with morphisms given by $k$-linear functors that preserve both limits and colimits, the adjoint functor theorem gives an equivalence $(\Pr^L_{k,\omega})\op = \Pr^R_{k,\omega}$.

The study of $\Pr^L_{k,\omega}$ may be reduced to the study of certain small $k$-linear $\infty$-categories as follows.
Following \cite{ag_brauer}, we say that a $k$-linear $\infty$-category $\Cc_0$ is \emph{perfect} if $\Cc_0$ is stable, and every idempotent $p: c_0 \to c_0$ in $\Cc_0$ induces a direct sum decomposition $c_0 = c_0' \oplus c_0''$ with $c_0', c_0'' \in \Cc_0$.
Let $\Cat_k^\perf$ denote the category of perfect $k$-linear $\infty$-categories and exact (i.e.\ finite (co)limit-preserving) $k$-linear functors.
For $\Cc \in \Pr^L_{k,\omega}$, let $\Cc^\omega$ be the full subcategory consisting of all compact objects in $\Cc$.
Then $(-)^\omega$ defines an equivalence $\Pr^L_{k,\omega} \xrightarrow{\sim} \Cat_k^\perf$, with inverse given by $\Ind: \Cat_k^\perf \to \Pr^L_{k,\omega}$.

\begin{rmk}
    Using \cite[Theorem 1.10]{bgt_universal}, we may obtain the category $\Cat_k^\perf$ from $\Cat_k$ by ``localizing along the Morita equivalences.''
    That is, the functor 
    \begin{align*}
        \Cat_k &\to \Cat_k^\perf \\
        \Cc_0 &\mapsto \Indsf(\Cc_0)^\omega
    \end{align*}
    is the universal functor out of $\Cat_k$ inverting Morita equivalences of $k$-linear $\infty$-categories (including e.g.\ Morita equivalences of $k$-algebras).
    Its right adjoint is the inclusion $\Cat_k^\perf \hookrightarrow \Cat_k$, which is necessarily fully faithful.
\end{rmk}

\subsection{Derived algebraic geometry and perfect stacks}

Derived algebraic geometry is the study of (geometric) stacks on a suitable site of \emph{derived affine schemes}.
As in classical algebraic geometry, the $\infty$-category $\dAff_k$ of derived affine schemes over $k$ is (defined to be) the opposite of the $\infty$-category of \emph{derived commutative $k$-algebras} $\dCAlg_k$.
Write $\Spec: \dCAlg_k\op \xrightarrow{\sim} \dAff_k$ for the equivalence.

What one means by ``derived commutative rings'' varies depending on one's goals, but a few standard definitions include:
\begin{itemize}
    \item The $\infty$-category of connective commutative dg-algebras.
    \item The $\infty$-category of simplicial commutative rings.
    \item The $\infty$-category of connective $\EE_\infty$-ring spectra.
\end{itemize}
When working over a field $k$ of characteristic zero (as we do generally), all of these approaches are equivalent, so we will not concern ourselves much with the differences between these.
Readers will not lose anything of much significance if they restrict to the case of commutative dg-algebras.

We shall equip $\dAff_k$ with the \'etale topology (see e.g.\ \cite[\S 1.2.2]{gr}).
A \emph{derived stack} over $k$ is a presheaf $\Xfr \in \PSh(\dAff_k)$ which satisfies descent for the \'etale topology.
Write $\dStk_k$ for the $\infty$-category of derived stacks.

\begin{rmk}
    As usual for the theory of stacks, there are many other topologies one could consider (e.g.\ smooth, fppf, fpqc, \dots).
    However, most stacks that appear in practice (including those we will deal with) satisfy descent for all of the standard choices of topology.
\end{rmk}

We may define the (derived) $\infty$-category of \emph{quasicoherent sheaves} on a derived stack $\Xfr$ as\footnote{This limit is \emph{a priori} large but can be reduced to a small limit by writing $\Xfr$ as a small colimit of derived affine schemes as in \cite[\S 3.1]{bzfn_integral}.}
\[
    \QC(\Xfr) = \lim_{\Spec R \to \Xfr} \Dsf(R).
\]
Within this $\infty$-category there is a full subcategory of \emph{perfect complexes}
\[
    \Perf(\Xfr) = \lim_{\Spec R \to \Xfr} \Perf(R),
\]
where $\Perf(R) = \Dsf(R)^\omega$. 
Both $\QC(\Xfr)$ and $\Perf(\Xfr)$ are $k$-linear stable $\infty$-categories.
The $\infty$-category $\QC(\Xfr)$ is presentable, and $\Perf(\Xfr)$ is perfect.
From the definition as limits, one sees that both $\QC(\Xfr)$ and $\Perf(\Xfr)$ admit natural tensor products $\otimes_{\Osc_\Xfr}$, and a morphism $f: \Xfr \to \Yfr$ induces a symmetric monoidal functor $f^*: \big(\QC(\Yfr), \otimes_{\Osc_\Yfr}\big) \to \big(\QC(\Xfr), \otimes_{\Osc_\Xfr}\big)$.
For notational simplicity, we shall often write $\otimes_\Osc$ for $\otimes_{\Osc_\Xfr}$ when the stack $\Xfr$ is clear from context.

Furthermore (cf.\ \cite[\S 3.1.5]{gr}), the usual $t$-structures on $\Dsf(R)$ (where $\Dsf^{\leq 0}(R)$ consists of connective modules) induce $t$-structures on the categories $\QC(\Xfr)$ and $\Perf(\Xfr)$.
The functors $\otimes_{\Osc}$ and $f^*$ are right $t$-exact for the usual $t$-structure on $\QC(-)$.
Note that the hearts of these $t$-structures may behave poorly: if $R$ is a derived ring, there is no ``abelian category of $R$-modules'' equivalent to $\Dsf(R)^\heartsuit$.
We shall not make heavy use of these $t$-structures, and we mention them only so that we may appeal to \cite[Theorem 1.3]{bhl_tannaka} later on.

In general, perfect complexes may not be compact in $\QC(\Xfr)$, and they may fail to generate $\QC(\Xfr)$.
Following \cite{bzfn_integral}, we say a derived stack $\Xfr$ is \emph{perfect} if $\Xfr$ has affine diagonal and $\Perf(\Xfr)$ compactly generates $\QCoh(\Xfr)$.
Write $\dStk_k^\perf$ for the full subcategory of $\dStk_k$ consisting of perfect stacks.

Most ``small'' derived stacks in characteristic zero are perfect -- see \cite[\S 3.3]{bzfn_integral} for sufficient criteria.
In particular, when $\charop k = 0$, the classifying stack $BG$ of any linear algebraic group $G$ is perfect.
Furthermore, if $\Yfr$ is perfect and $f: \Xfr \to \Yfr$ is affine, then $\Xfr$ is perfect.

Morphisms of perfect stacks enjoy many of the standard sheaf-theoretic identities of algebraic geometry.
If $f: \Xfr \to \Yfr$ is a morphism of perfect stacks, then the pullback functor $f^*: \QC(\Yfr) \to \QC(\Xfr)$ admits a colimit-preserving right adjoint $f_*: \QC(\Xfr) \to \QC(\Yfr)$.
These satisfy base change and the projection formula by \cite[Proposition 3.10]{bzfn_integral}.
There is also an equivalence $\QC(X \times Y) = \QC(\Xfr) \otimes_k \QC(\Yfr)$ by \cite[Theorem 4.7]{bzfn_integral}.
In \cref{sub:3ff_qc} we will construct a ``three-functor formalism'' on $\dStk_k$, ensuring that all of the higher homotopy coherence relations for these functors behave ``as expected.''

Perfect stacks also support a good theory of Fourier-Mukai transforms by \cite[Theorem 1.2]{bzfn_integral}.
More precisely, if $\Xfr$ and $\Yfr$ are perfect stacks, then there is a pair of mutually inverse equivalences
\[
    \Ksc_{-} : \Fun^L_k\big(\QC(\Xfr), \QC(\Yfr)\big) \xleftrightarrow{\sim} \QC(\Xfr \times \Yfr) : \Phi_{-},
\]
where $\Phi_\Fsc = \pi_{2*}(\pi_1^*(-) \otimes \Fsc)$ is the \emph{Fourier-Mukai transform} associated with $\Fsc \in \QC(\Xfr \times \Yfr)$, and $\Ksc_F$ is the \emph{Fourier-Mukai kernel} associated with $F \in \Fun^L_k\big(\QC(\Xfr), \QC(\Yfr)\big)$.
The analogous claim for perfect complexes holds when $\Xfr$ and $\Yfr$ are smooth and proper over $k$.
Note the distinction between the above result and the triangulated theory -- here every reasonable functor is automatically uniquely / functorially a Fourier-Mukai transform!

\subsection{Resolutions of the diagonal} \label{sub:resolution}

Suppose $\Xfr$ is a perfect stack.
The relationship between generators of $\Xfr$ and resolutions of the diagonal of $\Xfr$ is well-known and classical (going back to \cite{beilinson1978coherent}).
We shall review this relationship in modern language for future reference.

Let $\Delta_\Xfr: \Xfr \to \Xfr \times \Xfr$ be the diagonal morphism of $\Xfr$.
The identity functor $\id_\Xfr: \QC(\Xfr) \to \QC(\Xfr)$ may be understood as the Fourier-Mukai transform with kernel $\Delta_{\Xfr*} \Osc_\Xfr \in \QC(\Xfr \times \Xfr)$.

\begin{lem} \label{lem:identity_resolution_diagonal}
    Let $\Xfr$ be a perfect stack, and suppose that $\Delta_{\Xfr*} \Osc_\Xfr = \colim_{i \in I} \Asc_i \boxtimes \Bsc_i$ for some families $\{ \Asc_i \}_{i \in I}, \{ \Bsc_i \}_{i \in I} \subset \QC(\Xfr)$.
    Then, for any $\Fsc \in \QC(\Xfr)$, we have
    \[
        \Fsc = \colim_{i \in I} \Gamma(\Xfr, \Fsc \otimes \Asc_i) \otimes_k \Bsc_i.
    \]
\end{lem}

\begin{proof}
    For $i = 1, 2$, let $\pi_i: \Xfr \times \Xfr \to \Xfr$ be projection onto the $i$th factor.
    Note that $\pi_{2*}$ preserves colimits because $\pi_2$ is a morphism of perfect stacks (\cite[Proposition 3.10]{bzfn_integral}).
    Thus we may compute
    \begin{align*}
    \Fsc &= \pi_{2*}(\pi_1^* \Fsc \otimes \Delta_{\Xfr*} \Osc_\Xfr) \\
    &= \pi_{2*}\bigg(\pi_1^* \Fsc \otimes_\Osc \colim_{i \in I}(\Asc_i \boxtimes \Bsc_i)\bigg) \\
    &= \colim_{i \in I} \pi_{2*}\big((\Fsc \otimes_\Osc \Asc_i) \boxtimes \Bsc_i\big) \textrm{ because all functors involved commute with colimits} \\
    &= \colim_{i \in I} \pi_{2*}\big(\pi_1^*(\Fsc \otimes_\Osc \Asc_i) \otimes_\Osc \pi_2^* \Bsc_i\big) \\
    &= \colim_{i \in I} \pi_{2*}\big(\pi_1^*(\Fsc \otimes_\Osc \Asc_i)\big) \otimes_k \Bsc_i \textrm{ by the projection formula} \\
    &= \colim_{i \in I} \Gamma(\Xfr, \Fsc \otimes_\Osc \Asc_i) \otimes_k \Bsc_i \textrm{ by base change.} \qedhere
\end{align*}
\end{proof}

In the situation of \cref{lem:identity_resolution_diagonal}, we see that the family $\{ \Bsc_i \}_{i \in I}$ generates $\QC(\Xfr)$.
Conversely, if $\{ \Bsc_i \}_{i \in I}$ generates $\QC(\Xfr)$, then $\{ \Bsc_i \boxtimes \Bsc_j \}_{i, j \in I}$ generates $\QC(\Xfr \times \Xfr) = \QC(\Xfr) \boxtimes \QC(\Xfr)$, so we may write $\Delta_{\Xfr*} \Osc_\Xfr$ as a colimit over these sheaves.

We may use this perspective to rewrite any colimit-preserving functor $F: \QC(\Xfr) \to \QC(\Yfr)$ (where $\Xfr$ and $\Yfr$ are perfect stacks) in terms of a given resolution of the diagonal of $\Xfr$.
More precisely:

\begin{prop} \label{prop:fm_resolution_diagonal}
    Let $\Xfr$ and $\Yfr$ be perfect stacks, and write $\Delta_{\Xfr*} \Osc_\Xfr = \colim_{i \in I} \Asc_i \boxtimes \Bsc_i$.
    If $F: \QC(\Xfr) \to \QC(\Yfr)$ is a colimit-preserving functor, then $\Ksc_F = \colim_{i \in I} \Asc_i \boxtimes F(\Bsc_i)$. 
\end{prop}

\begin{proof}
    This is a direct computation using \cref{lem:identity_resolution_diagonal}:
    \begin{align*}
        F(\Fsc) &= F\bigg(\colim_{i \in I} \Gamma(\Yfr, \Fsc \otimes \Asc_i) \otimes_k \Bsc_i\bigg) \\
        &= \colim_{i \in I} \Gamma(\Yfr, \Fsc \otimes \Asc_i) \otimes_k F(\Bsc_i) \\
        &= \colim_{i \in I} \pi_{2*}\big((\Fsc \otimes \Asc_i) \boxtimes F(\Bsc_i)\big) \\
        &= \pi_{2*}\bigg(\pi_1^* \Fsc \otimes \colim_{i \in I} \big(\Asc_i \boxtimes F(\Bsc_i)\big)\bigg). \qedhere
    \end{align*}
\end{proof}

\subsection{Three-functor formalisms in general} \label{sub:3ff}

To rigorously construct monoidal convolution products on categories of sheaves, it is useful to have access to a ``three-functor formalism'' for said categories.
One can define this loosely as follows.
Fix some $\infty$-category of ``geometric spaces,'' e.g.\ $\dAff_k$ or $\dStk_k$.
A three-functor formalism should assign to each object $X$ a category $\Dsf(X)$ of ``sheaves on $X$'' together with:
\begin{itemize}
    \item For every $X$, a symmetric monoidal tensor product $\otimes$ on $\Dsf(X)$,
    \item For every $f: X \to Y$, a $*$-pullback functor $f^*: \Dsf(Y) \to \Dsf(X)$, and
    \item For ``nice'' morphisms $f: X \to Y$, a $!$-pushforward functor $f_!: \Dsf(X) \to \Dsf(Y)$,
\end{itemize}
satisfying base change and the projection formula.

Three-functor formalisms, as well as the stronger notion of six-functor formalisms, are made mathematically precise (in the language of $(\infty, 1)$-categories) in \cite[Appendix A.5]{mann_6ff}, building on previous works \cite{lz_enhanced, gr}.
We recall (a slightly simplified version of) Mann's definition.

Let $(\Cc, E)$ be a pair where $\Cc$ is an $\infty$-category with finite limits and $E$ is a collection of morphisms in $\Cc$ which contains all isomorphisms and is stable under homotopy, composition, and pullback.
One can define an $\infty$-category $\Corr(\Cc, E)$ of \emph{correspondences in $\Cc$ with right leg in $E$}.
The objects of $\Corr(\Cc, E)$ are the objects of $\Cc$.
Morphisms $X_1 \to X_2$ in $\Corr(\Cc, E)$ are given by correspondences (also called ``spans'') $X_1 \xleftarrow{f} Y \xrightarrow{g} X_2$ where $f$ is a morphism in $\Cc$ and $g$ is a morphism in $E$.\footnote{In particular, the morphism spaces in $\Corr(\Cc, E)$ are typically no longer small.
Any set-theoretic difficulties this presents can be circumvented by standard universe-based arguments, and we will ignore such difficulties for simplicity of exposition.}
Composition of correspondences is given by fiber products\footnote{Because fiber products are only unique up to natural equivalence, composition is only well-defined up to coherent homotopy.
In particular, the rigorous definition of $\Corr(\Cc, E)$ is somewhat technical (as usual for the subject).}
\[
\begin{tikzcd}
     & & Y_1 \times_{X_2} Y_2 \ar[dl] \ar[dr] & & \\
     & Y_1 \ar[dl] \ar[dr] & & Y_2 \ar[dl] \ar[dr] & \\
    X_1 & & X_2 & & X_3.
\end{tikzcd}
\]
The category $\Corr(\Cc, E)$ has a symmetric monoidal structure $\times$ defined using the Cartesian monoidal structure on $\Cc$.

\begin{dfn}[{\cite[Definition A.5.6]{mann_6ff}}] \label{dfn:6ff}
    A \emph{three-functor formalism} on $(\Cc, E)$ is a lax symmetric monoidal functor $\Sh: \big(\Corr(\Cc, E), \times\big) \to \big(\widehat{\Cat}_\infty, \times\big)$, where $\widehat{\Cat}_\infty$ is the $\infty$-category of large $\infty$-categories (in some universe).
    This induces the following ``three functors:''
    \begin{itemize}
        \item The inclusion map $\Cc\op \hookrightarrow \Corr(\Cc, E)$, given on morphisms by $(X \xrightarrow{f} Y) \mapsto (Y \xleftarrow{f} X \xrightarrow{\id_X} X)$, is symmetric monoidal when viewed as a functor $(\Cc\op, \times) \to (\Corr(\Cc, E), \times)$.
        Every object of $\Cc\op$ is uniquely a commutative algebra for $\times$ (the coproduct in $\Cc\op$) by \HA{2.4.3.10}.
        Thus every category $\Sh(X)$ for $X \in \Cc$ has a natural symmetric monoidal structure $\otimes_{\Sh(X)}$.
        \item Similarly, for a morphism $f: X \to Y$ in $\Cc$, we obtain a symmetric monoidal functor $f^*: (\Sh(Y), \otimes_{\Sh(Y)}) \to (\Sh(X), \otimes_{\Sh(X)})$.
        \item Let $\Cc_E$ be the subcategory of $\Cc$ containing all objects of $\Cc$ but only containing the morphisms in $E$.
        Consider the inclusion map $\Cc_E \hookrightarrow \Corr(\Cc, E)$ given on morphisms by $(X \xrightarrow{f} Y) \mapsto (X \xleftarrow{\id_X} X \xrightarrow{f} Y)$.
        For a morphism $f: X \to Y$ in $E$, we obtain a functor $f_!: \Sh(X) \to \Sh(Y)$.
    \end{itemize}
    If all functors $\otimes$, $f^*$, and $f_!$ admit right adjoints, then $\Sh$ is called a \emph{six-functor formalism}.
\end{dfn}

\begin{rmk}
    We note that, at the level of objects, the lax symmetric monoidal structure on $\Sh$ corresponds to an ``external tensor product'' $\boxtimes: \Sh(X) \times \Sh(Y) \to \Sh(X \times Y)$.
    The fact that $\Sh$ is a lax symmetric monoidal \emph{functor} implies that $\boxtimes$ satisfies the ``expected'' compatibilities with the pullback / pushforward functors, e.g.\ $\Fsc \otimes \Gsc = \Delta_X^* (\Fsc \boxtimes \Gsc)$ for $X \in \Cc$, $\Fsc, \Gsc \in \Sh(X)$, and $\Delta_X: X \to X \times X$ the diagonal.
    (We will explicitly mention these compatibilities when needed.)
\end{rmk}

Six-functor formalisms behave best in the presentable case:

\begin{dfn}
    A \emph{presentable six-functor formalism} is a six-functor formalism $\Sh$ such that every $\infty$-category $\Sh(X)$ (for $X \in \Cc$) is presentable.
\end{dfn}

In other words, a presentable six-functor formalism is a lax symmetric monoidal functor $\Sh: \big(\Corr(\Cc, E), \times\big) \to \big(\Pr^L, \times\big)$.
As noted in \cite[proof of Lemma A.5.11]{mann_6ff}, any presentable six-functor formalism defines a lax symmetric monoidal functor\footnote{Here $\otimes$ denotes the Lurie tensor product on $\Pr^L$.}
\[
    \Sh: \big(\Corr(\Cc, E), \times\big) \to \big(\Pr^L, \otimes\big)
\]
by postcomposition with the lax symmetric monoidal identity functor $\big(\Pr^L, \times\big) \to \big(\Pr^L, \otimes\big)$.

In fact, there is even more structure on the categories $\Sh(X)$ when $\Sh$ is presentable.
If $\Sh$ is a presentable six-functor formalism, then all categories $\Sh(X)$ and all functors $\otimes$, $f^*$, and $f_!$ are naturally enriched and tensored over $\Sh(\pt)$, where $\pt$ is the terminal object of $\Cc$ (see \cite[Theorem A.3.8]{mgs_universal}).
This allows us to view $\Sh$ as a lax symmetric monoidal functor
\begin{equation} \label{eq:presentable_6ff}
    \Sh: \big(\Corr(\Cc, E), \times\big) \to \big(\Pr^L_{\Sh(\pt)}, \otimes^L_{\Sh(\pt)}\big),
\end{equation}
where $\Pr^L_{\Sh(\pt)}$ is the symmetric monoidal $\infty$-category of presentably $\Sh(\pt)$-enriched $\infty$ and $\Sh(\pt)$-enriched left adjoint functors.
We will often abuse notation and refer to the lax symmetric monoidal functor of \eqref{eq:presentable_6ff} as a ``presentable six-functor formalism.''

\subsection{Three-functor formalisms and convolution} \label{sub:3ff_convolution}

In this subsection we shall discuss the \emph{convolution products} on categories of sheaves obtained from a three-functor formalism.
We begin with a cautionary remark.

\begin{rmk}
    Let $(\Cc, E)$ be as in \cref{sub:3ff}.
    The category $\Cc_E$ may not have products (and even if it does, these products need not agree with the products in $\Cc$).
    Indeed, for $X_1, X_2 \in \Cc$, there is no reason that the projection morphisms $\pi_i: X_1 \times X_2 \to X_i$ must be in $E$.
    Thus it does not make sense \emph{a priori} to claim that the functor $(\Cc_E, \times) \to (\Corr(\Cc, E), \times)$ is symmetric monoidal.

\end{rmk}

Under the given hypotheses on $E$, we expect that the Cartesian product of $\Cc$ defines a non-Cartesian symmetric monoidal structure (also denoted $\times$) on $\Cc_E$.
Furthermore, the inclusion functor $(\Cc_E, \times) \to (\Corr(\Cc, E), \times)$ should be symmetric monoidal.
A proof of these claims under an additional assumption (which does not hold in our case of interest) is given in \cite[Proposition 2.3.7]{heyer20246ff}.
However, proving these claims in our case of interest appears to involve more tinkering with the ``internal machinery'' of $\infty$-categories of correspondences than we care to do in this paper.

Thus, \emph{in this section, we will assume for simplicity that $E$ consists of all morphisms in $\Cc$}.
(In practice, we may often restrict to this case by replacing $(\Cc, E)$ by $(\Cc', \Mor \Cc')$ where $\Cc'$ is a full subcategory of ``nice objects'' of $\Cc$.)

\begin{notn}
    When $E$ is the class of all morphisms in $\Cc$, we will write $\Corr(\Cc) = \Corr(\Cc, E)$ and call $\Corr(\Cc)$ the \emph{category of correspondences in $\Cc$}.
\end{notn}

The inclusion $(\Cc, \times) \to (\Corr(\Cc), \times)$ is symmetric monoidal.
Indeed, because $\Corr(\Cc)$ is equivalent to its opposite category, it suffices to show that $(\Cc\op, \times) \to (\Corr(\Cc), \times)$ is symmetric monoidal.
But this was noted above.

Let $\Sh: \Corr(\Cc) \to \widehat{\Cat}_\infty$ be a three-functor formalism.
Then the composite functor
\[
    \Sh_!: (\Corr(\Cc), \times) \to \big(\widehat{\Cat}_\infty, \times\big)
\]
is lax symmetric monoidal.
In particular, for an $\EE_n$-algebra object $M$ of $\Cc$, applying $\Sh_!$ to $M$ yields an $\EE_n$-monoidal \emph{convolution product} $\star_M$ on $\Sh(M)$.
As a binary operation, $\star_M$ is given by $\Fsc \star_M \Gsc = \mu_!(\Fsc \boxtimes \Gsc)$ for $\Fsc, \Gsc \in \Sh(M)$, where $\mu: M \times M \to M$ is the binary multiplication.
The construction of $\star_M$ is functorial: if $f: M \to N$ is a homomorphism of $\EE_n$-monoids in $\Cc_E$, then $f_!: (\Sh(M), \star_M) \to (\Sh(N), \star_N)$ is $\EE_n$-monoidal.

Under some additional hypotheses, the convolution product combined with the diagonal map $\Delta_M: M \to M \times M$ equips $\Sh(M)$ with the structure of an $(\EE_n, \EE_\infty)$-bialgebra.
Before proving this, we recall the $\infty$-categorical definition of bialgebras:

\begin{dfn} \label{dfn:bialgebra}
    Let $(\Dc, \otimes_\Dc)$ be a symmetric monoidal $\infty$-category.
    For any $n, m \in \NN \cup \{ \infty\}$, the category of \emph{$(\EE_n, \EE_m)$-bialgebras in $(\Dc, \otimes_\Dc)$} is 
    \[
        \BiAlg_{\EE_n,\EE_m}(\Dc, \otimes_\Dc) := \Alg_{\EE_m}\big((\Alg_{\EE_n}(\Dc, \otimes_\Dc))\op, \otimes_\Dc\big)\op.
    \]
    We will frequently leave the tensor product $\otimes_\Dc$ implicit, writing $\BiAlg_{\EE_n,\EE_m}(\Dc) := \BiAlg_{\EE_n,\EE_m}(\Dc, \otimes_\Dc)$.
\end{dfn}

Informally, an $(\EE_n, \EE_m)$-bialgebra in a symmetric monoidal $\infty$-category $\Dc$ is an object $d$ of $\Dc$ together with:
\begin{itemize}
    \item an $\EE_n$-algebra structure on $d$,
    \item an $\EE_m$-coalgebra structure on $d$, and
    \item extra data controlling the (homotopy coherent) compatibility of the two structures.
\end{itemize}
Because (strong) symmetric monoidal functors preserve $\EE_n$-algebras and $\EE_m$-coalgebras, a symmetric monoidal functor $\Dc_1 \to \Dc_2$ induces a functor $\BiAlg_{\EE_n,\EE_m}(\Dc_1) \to \BiAlg_{\EE_n,\EE_m}(\Dc_2)$.

\begin{prop} \label{prop:6ff_bialgebra}
    Let $\Sh: (\Corr(\Cc), \times) \to (\Pr^L, \otimes)$ be a presentable six-functor formalism.
    Assume that the induced functor $\Sh: (\Cc, \times) \to (\Pr^L_{\Sh(\pt)}, \otimes_{\Sh(\pt)})$ is (strong) symmetric monoidal.
    If $M$ is an $\EE_n$-monoid in $\Cc$, then:
    \begin{enumerate}
        \item $*$-pullback along the structure maps and the diagonal of $M$ make $\Sh(M)$ into an $(\EE_\infty, \EE_n)$-bialgebra in $\Pr^L_{\Sh(\pt)}$.
        \item $!$-pushforward along the structure maps and the diagonal of $M$ make $\Sh(M)$ into an $(\EE_n, \EE_\infty)$-bialgebra in $\Pr^L_{\Sh(\pt)}$.
    \end{enumerate}
    Furthermore, these constructions are functorial in $M$.
\end{prop}

\begin{proof}
    (1) The category $\Alg_{\EE_n}(\Cc)$ is Cartesian monoidal, so by \HA{2.4.3.10} there is an equivalence $\Alg_{\EE_n}(\Cc) \xrightarrow{\sim} \BiAlg_{\EE_n,\EE_\infty}(\Cc)$.
    Thus we may view $M$ as an $(\EE_n, \EE_\infty)$-bialgebra in $\Cc$.
    Reversing the direction of arrows makes $M$ into an $(\EE_\infty, \EE_n)$-bialgebra in $\Cc\op$.
    The composite functor
    \[
        (\Cc\op, \times) \to (\Corr(\Cc), \times) \to (\Pr^L_{\Sh(\pt)}, \otimes_{\Sh(\pt)})
    \]
    is symmetric monoidal, hence must send $(\EE_\infty, \EE_n)$-bialgebras to $(\EE_\infty, \EE_n)$-bialgebras.

    (2) The same argument works here.
\end{proof}

\begin{rmk}
    Part (1) of \cref{prop:6ff_bialgebra} does not require that $E$ is the class of all morphisms in $\Cc$, though part (2) does.
\end{rmk}

We state \cref{prop:6ff_bialgebra} in terms of the induced functor $(\Cc, \times) \to (\Pr^L_{\Sh(\pt)}, \otimes_{\Sh(\pt)})$ rather than $(\Cc, \times) \to (\Pr^L_{\Sh(\pt)}, \otimes_{\Sh(\pt)})$ because the former is rarely symmetric monoidal while the latter often is so.
This will be the case for the six-functor formalism $\QC: (\Corr(\dStk_k^\perf), \times) \to (\Pr^L_k, \otimes_k)$ discussed below (see \cref{prop:3ff_qc}(3)).

\subsection{A three-functor formalism for quasicoherent sheaves} \label{sub:3ff_qc}

For our applications, we would like to use a three-functor formalism for quasicoherent sheaves on stacks in which the ``$!$-pushforward functors'' are the right adjoints to the usual $*$-pullback functors (in reasonable cases, these are the usual $*$-pushforward functors).
The analogous three-functor formalism for qcqs schemes is discussed in \cite[\S 8.3]{scholze_6ff}.
We will need to extend the formalism to general derived stacks -- this can be done using the results of \cite{bzfn_integral}.

We begin by clarifying the notion of ``nice'' morphisms in this context.

\begin{dfn}
    A morphism of stacks $f: \Xfr \to \Yfr$ is \emph{very perfect} if $f$ has affine diagonal and, for every morphism of stacks $g: \Zfr \to \Yfr$ with $\Zfr$ perfect, the fiber product $\Xfr \times_\Yfr \Zfr$ is perfect.
    Let $\VP$ be the class of very perfect morphisms in $\dStk_k$.
\end{dfn}

\begin{lem} \label{lem:very_perfect}
    The class $\VP$ has the following properties:
    \begin{enumerate}
        \item $\VP$ is stable under composition.\footnote{It is not clear \emph{a priori} whether this holds for the \emph{perfect morphisms} of \cite[Definition 3.2]{bzfn_integral}.}
        \item $\VP$ is stable under base change.
        \item $\VP$ contains all affine morphisms.
        \item $\VP$ satisfies a cancellation law: if $f: \Xfr \to \Yfr$ and $g: \Yfr \to \Zfr$ are such that $g \in \VP$ and $g \circ f \in \VP$, then $f \in \VP$.
        \item Any morphism between perfect stacks is in $\VP$.
    \end{enumerate}   
\end{lem}

\begin{proof}
    This is standard diagram chasing and references to \cite{bzfn_integral}.
    We first note that the class of morphisms with affine diagonal is stable under composition and base change, e.g.\ by (the argument of) \cite[Theorem 11.1.2]{vakil}.
    
    (1). Let $f: \Xfr \to \Yfr$ and $g: \Yfr \to \Zfr$ be very perfect.
    To see that $g \circ f$ is very perfect, let $\Wc \to \Zfr$ be a morphism with $\Wc$ perfect.
    Then $\Xfr \times_\Zfr \Wc = (\Yfr \times_\Zfr \Wc) \times_\Yfr \Xfr$.
    Because $g$ is very perfect, the stack $\Yfr \times_\Zfr \Wc$ is perfect.
    Then, because $f$ is very perfect, the stack $\Xfr \times_\Zfr \Wc$ must also be perfect.
    It follows that $g \circ f$ is very perfect.

    (2). Let $\Xfr \to \Yfr$ be very perfect and let $\Zfr \to \Yfr$ be any morphism.
    If $W \to \Zfr$ is a morphism with $W$ perfect, then $(\Xfr \times_\Yfr \Zfr) \times_\Zfr W = \Xfr \times_\Yfr W$ is perfect.
    Thus the base change $\Xfr \times_\Yfr \Zfr \to \Zfr$ is very perfect.

    (3). Any affine morphism has affine diagonal.
    Furthermore, if $\Xfr \to \Yfr$ is affine and $\Zfr \to \Yfr$ is any morphism with $\Zfr$ perfect, then $\Xfr \times_\Yfr \Zfr \to \Zfr$ is affine.
    Then $\Xfr \times_\Yfr \Zfr$ is perfect by \cite[Proposition 3.21]{bzfn_integral}, so $\Xfr \to \Yfr$ must be very perfect.

    (4). This is a consequence of (the argument of) \cite[Theorem 11.1.1]{vakil}, where we use (3) above to see that the diagonal of a very perfect morphism is very perfect.

    (5). By (4) it suffices to show that, if $\Xfr$ is a perfect stack, then the natural map $\pi: \Xfr \to \Spec k$ is very perfect.
    The map $\pi$ has affine diagonal because the same is true for $\Xfr$.
    Furthermore, if $\Yfr$ is a perfect stack, then $\Xfr \times \Yfr$ is perfect by \cite[Proposition 3.24]{bzfn_integral}, giving the claim.
\end{proof}

By parts (1) and (2) of \cref{lem:very_perfect}, the pair $(\dStk_k, \VP)$ satisfies the conditions discussed in \cref{sub:3ff}, so the category $\Corr(\dStk_k, \VP)$ is well-defined.
We may now define our three-functor formalisms of interest:

\begin{prop} \label{prop:3ff_qc}
    \begin{enumerate}
        \item There is a three-functor formalism
    \[
        \QC: (\Corr(\dStk_k, \VP), \times) \to (\widehat{\Cat}_\infty, \times)
    \]
    sending a correspondence $\Xfr \xleftarrow{f} \Zfr \xrightarrow{g} \Yfr$ with $g$ very perfect to
    \[
        \QC(\Xfr) \xrightarrow{f^*} \QC(\Zfr) \xrightarrow{g_*} \QC(\Yfr).
    \]
        \item The same formula defines a presentable six-functor formalism
        \[
            \QC: (\Corr(\dStk_k^\perf), \times) \to (\Pr^L, \otimes)
        \]
        \item The induced functor
        \[
            \QC: (\Corr(\dStk_k^\perf), \times) \to (\Pr^L_k, \otimes_k)
        \]
        is (strong) symmetric monoidal.
    \end{enumerate}
\end{prop}

\begin{proof}
    (1). We have already constructed $\QC$ as a functor $\dStk_k\op \to \CAlg(\widehat{\Cat}_\infty)$ (sending $f: \Xfr \to \Yfr$ to $f^*: \QC(\Yfr) \to \QC(\Xfr)$.
    By \cite[Proposition A.5.10]{mann_6ff},\footnote{More precisely: \cite[Proposition A.5.10]{mann_6ff} applies to a ``suitable decomposition'' of $\VP$ into classes of morphisms $I$ and $P$ satisfying certain conditions.
    Here we simplify things by taking $P = \VP$ and taking $I$ to be the class of isomorphisms in $\dStk_k$.} to extend $\QC$ to a three-functor formalism, it suffices to show that:
    \begin{itemize}
        \item $\VP$ satisfies the cancellation law of \cref{lem:very_perfect} (4).
        \item For $f: \Xfr \to \Yfr$ in $\VP$ and any $g: \Zfr \to \Yfr$, let $f': \Xfr \times_\Yfr \Zfr \to \Zfr$ and $g': \Xfr \times_\Yfr \Zfr \to \Xfr$ be the induced maps.
        Then the base change formula $g^* f_* = (f')_* (g')^*$ holds.
        \item The projection formula holds for $f \in \VP$.
    \end{itemize}
    The former condition holds by \cref{lem:very_perfect} (4), while the latter two conditions hold by \cite[Proposition 3.10]{bzfn_integral}.

    (2). Note that the categories $\QC(\Xfr)$ are automatically presentable (even if $\Xfr$ is not perfect).
    Thus it suffices to show that, if $f: \Xfr \to \Yfr$ is a morphism of perfect stacks, then the functors $\otimes_{\Osc_\Xfr}$, $f^*$, and $f_*$ all admit right adjoints.
    For $\otimes_{\Osc_\Xfr}$ and $f^*$ this is standard (and does not use perfectness of the stacks involved).
    For $f_*$, we know by \cite[Proposition 3.10]{bzfn_integral} that $f_*$ preserves all colimits, so the existence of a right adjoint to $f_*$ follows from the adjoint functor theorem \HTT{5.5.2.9}.

    (3). This follows from \cite[Theorem 4.7]{bzfn_integral} and the above discussion.
\end{proof}

\begin{cor} \label{cor:QC_convolution_bialgebra}
    If $\Mfr$ is an $\EE_n$-monoid object in $\dStk_k^\perf$, then:
    \begin{enumerate}
        \item $*$-pullback along the structure maps and the diagonal of $\Mfr$ make $\QC(\Mfr)$ into an $(\EE_\infty, \EE_n)$-bialgebra in $\Pr^L_{\Sh(\pt)}$.
        \item $*$-pushforward along the structure maps and the diagonal of $\Mfr$ make $\QC(\Mfr)$ into an $(\EE_n, \EE_\infty)$-bialgebra in $\Pr^L_{\Sh(\pt)}$.
    \end{enumerate}
\end{cor}

\begin{proof}
    This follows by applying \cref{prop:6ff_bialgebra} to the presentable six-functor formalism of parts (2) and (3) of \cref{prop:3ff_qc}.
\end{proof}

\begin{rmk}
    Although part (2) of \cref{prop:3ff_qc} defines a presentable six-functor formalism, the right adjoint to $f_*$ behaves poorly in general, even for proper morphisms.
    For example, $f_*$ typically does not preserve compact objects, so its right adjoint does not preserve colimits.
    Dealing with these issues is one of the main reasons for the existence of the formalism of ind-coherent sheaves (as developed in \cite{gr} among others).
    We will not need to make direct use of the right adjoint of $f_*$, so we content ourselves with using $\QC$.
\end{rmk}

\section{Quiver presentations of derived categories} \label{sec:quiver}

\begin{notn} \label{notn:quiver_presentations}
    Throughout this section, fix a commutative reductive group $G$ (i.e.\ a product of a torus and a finite abelian group).
    Let $\phi: \Yfr \to BG$ be a morphism of derived stacks which is affine and almost of finite type. 
    Equivalently, there exists a unique affine derived scheme $Y$ almost of finite type over $\Spec k$ and a $G$-action $G \curvearrowright Y$ such that $\Yfr = [Y / G]$.
    Let $\Xfr$ be an open substack of $\Yfr$, and write $j_\Xfr: \Xfr \hookrightarrow \Yfr$ for the corresponding open immersion.
    Write $X = \Xfr \times_\Yfr Y$ so that $\fr X = [X/G]$.
\end{notn}

The derived category $\QC(\Yfr) = \QC([Y / G])$ may be understood as the derived category of $G$-equivariant quasicoherent sheaves on $Y$.
Because $G$ is reductive, for $\Fsc, \Gsc \in \QC(\Yfr)$, the $\Hom$-complexes $\Hom_\Yfr(\Fsc, \Gsc)$ may be computed by taking the $G$-invariants of the $\Hom$-complexes of the corresponding quasicoherent sheaves on $Y$.
An analogous description holds for $\Xfr$.
We write $\Osc_\Yfr(\chi)$ (resp.\ $\Osc_\Xfr(\chi)$) for the pullback of the line bundle (a.k.a.\ 1-dimensional representation) $\Osc_{BG}(\chi) \in \QC(BG)$ to $\Yfr$ (resp.\ $\Xfr$).

\subsection{Quivers and derived categories of quotient stacks}

We are interested in describing the symmetric monoidal category $(\QC(\Yfr), \otimes_\Osc)$ as a functor category out of some small $k$-linear category.

\begin{dfn} \label{dfn:weight_quiver}
    The \emph{total weight quiver} $\Qsf(\phi)$ is the essentially small full symmetric monoidal subcategory of $(\Perf(\Yfr), \otimes_\Osc)$ with objects given by $\{ \Osc_\Yfr(\chi) \}_{\chi \in \weight(G)}$.
    More generally, for any subset $S \subset \weight(G)$, the \emph{partial weight quiver} $\Qsf_S(\phi)$ corresponding to $S$ is the full (not necessarily monoidal) subcategory of $\Perf(\Yfr)$ with objects $\{ \Osc_\Yfr(\chi) \}_{\chi \in S}$.
\end{dfn}

\begin{rmk} \label{rmk:action_grading}
    Suppose $Y = \Spec R$ for a derived ring $R$.
    Then specifying the $G$-action on $Y$ is the same as specifying a $\weight(G)$-grading on the derived ring $R$, say $R = \oplus_{\chi \in \weight(G)} R_\chi$.
    From this perspective, quasicoherent sheaves on $\Yfr$ correspond to $\weight(G)$-graded $R$-modules.
    In this case, we have the following explicit (but less homotopy-coherent) description of the partial weight quivers $\Qsf_S(\phi)$:
    \begin{itemize}
        \item $\ob \Qsf_S(\phi) = S$.
        \item For $\chi_1, \chi_2 \in S$, we have $\Hom_{\Qsf_S(\phi)}(\chi_1, \chi_2) = R_{\chi_2 - \chi_1}$.
        \item Composition of morphisms in $\Qsf_S(\phi)$ is multiplication in $R$.
    \end{itemize}
    The symmetric monoidal structure on $\Qsf(\phi)$ is given on objects by addition in $\weight(G)$ and on morphisms\footnote{The behavior on morphisms is forced on us by the fact that the objects of $\Qsf(\phi)$ form a group under multiplication.} by multiplication in $R$.
\end{rmk}

\begin{ex} \label{ex:Beilinson_weight_quiver}
    If $G = \GG_m$ and $\Yfr = [\AA^{n+1} / \GG_m]$, the total weight quiver $\Qsf_\phi$ is the $k$-linearization of the infinite Beilinson quiver appearing in \cref{fig:infinite_beilinson_quiver}.
    The partial weight quiver for $\{ 0, \dots, n \} \subset \ZZ = \weight(\GG_m)$ is the $k$-linearization of the Beilinson quiver appearing in \cref{fig:beilinson_quiver}.
\end{ex}

Due to the above interpretations, we will often implicitly identify the objects of the category $\Qsf_S(\phi)$ with $S \subset \weight(G)$.
Following \cref{ex:Beilinson_weight_quiver} (and \cref{dfn:prelinear_quiver} below), we will generally think of $\Qsf_S(\phi)$ as a ``$k$-linear quiver with relations'' and think of $k$-linear functors from $\Qsf_S(\phi)$ (or $\Qsf_S(\phi)\op$) to $\Dsf(k)$ as derived quiver representations.

By \cref{prop:day_convolution_properties}, there exists a unique $k$-linear colimit-preserving symmetric monoidal structure $\odot_{\Qsf(\phi)}$ on $\Dsf(\Qsf(\phi)\op)$ extending the natural symmetric monoidal structure on $\Qsf(\phi)$ (viewed as a full subcategory of $\Dsf(\Qsf(\phi)\op)$ via the Yoneda embedding).
More concretely, for $V_1, V_2 \in \Dsf(\Qsf(\phi)\op) = \Fun_k(\Qsf(\phi)\op, \Dsf(k))$ and $\chi \in \weight(G)$, we have
\[
    (V_1 \odot_{\Qsf(\phi)} V_2)(\chi) = \bigoplus_{\chi_1 + \chi_2 = \chi} V_1(\chi_1) \otimes_k V_2(\chi_2).
\]

With the above definitions, it is not hard to give a description of $(\QC(\Yfr), \otimes_\Osc)$ in terms of $\Qsf(\phi)$.
Special cases and variants of this statement are well-known in the literature, especially in the context of equivariant mirror symmetry for affine toric varieties (see e.g.\ \cite[Proposition 2.1]{borisov2019equivariant}, as well as \cite[Theorem 1.1]{moulinos_geometry} and \cite[Proposition 3.3.1]{bai2025toric} for statements over the sphere spectrum).

\begin{prop} \label{prop:basic_underlying_equivalence}
    The inclusion $\Qsf(\phi) \hookrightarrow \QC(\Yfr)$ induces a symmetric monoidal equivalence 
    \[
       \Psi: (\Dsf(\Qsf(\phi)\op), \odot_{\Qsf(\phi)}) \simeq (\QC(\Yfr), \otimes_\Osc).
    \]
\end{prop}

\begin{proof}
    The functor $\Psi$ is symmetric monoidal by \cref{prop:day_convolution_properties}(1).
    It suffices to show that $\QC(\Yfr)$ is compactly generated by the line bundles $\{ \Osc_\Yfr(\chi) \}_{\chi \in \weight(G)}$ (as $\otimes_\Osc$ must then be the unique colimit-preserving $k$-linear extension of the symmetric monoidal structure on $\Qsf(\phi)$).
    This is standard, though we provide an ``intrinsically derived'' proof.
    
    Because $G$ is a reductive group, $\QC(BG)$ is compactly generated by the line bundles $\{\Osc_{BG}(\chi)\}_{\chi \in \weight(G)}$ (see e.g.\ \cite[Corollary 3.22]{bzfn_integral}).
    The natural map $\Yfr \to BG$ is affine, so we may write 
    \[
    \QC(\Yfr) = \Mod_{\Asc}(\QC(BG))
    \]
    for some $\Asc \in \CAlg(\QC^{\leq 0}(BG))$.
    Observe that $\Mod_{\Asc}(\QC(BG))$ is compactly generated by $\{\Asc(\chi)\}_{\chi \in \weight(G)}$.
    In fact, if $\Fsc \in \Mod_{\Asc}(\QC(BG))$ satisfies $\Hom_{\Mod_{\Asc}(\QC(BG))}(\Asc(\chi), \Fsc) = 0$ for all $\chi \in \weight(G)$, then by adjunction
    \[
        \Hom_{\QC(BG)}(\Osc_{BG}(\chi), \Fsc) = \Hom_{\Mod_{\Asc}(\QC(BG))}(\Asc(\chi), \Fsc) = 0
    \]
    for all $\chi \in \weight(G)$, so the underlying object of the $\Asc$-module $\Fsc$ is zero and we must have $\Fsc = 0$.
    It follows that $\QC(\Yfr)$ is compactly generated by $\{ \Osc_\Yfr(\chi) \}_{\chi \in \weight(G)}$.
\end{proof}

Combining \cref{prop:basic_underlying_equivalence} with Tannakian reconstruction theory gives a moduli interpretation of $\Yfr$ (which we will need in \cref{sec:basic}).

\begin{cor} \label{cor:tannakian_description}
    For all $R \in \dCAlg_k$, we have
    \[
        \Yfr(R) = \Fun^\otimes(\Qsf(\phi), \Perf^{\leq 0}(R))^\simeq.
    \]
\end{cor}

\begin{proof}
    By \cite[Theorem 1.3]{bhl_tannaka}, we may identify $\Yfr(R)$ with the subspace of $\Fun^{L,\otimes}_k(\QC(\Yfr), \Dsf(R))^\simeq$ consisting of functors which preserve connective objects.
    By \cref{prop:basic_underlying_equivalence} and \cref{prop:day_convolution_properties}(2), we may write
    \[
        \Fun^{L,\otimes}_k(\QC(\Yfr), \Dsf(R))^\simeq = \Fun^\otimes(\Qsf(\phi), \Dsf(R))^\simeq,
    \]
    and the subspace of $\Fun^{L,\otimes}(\QC(\Yfr), \Dsf(R))^\simeq$ consisting of functors which preserve connective objects is identified with $\Fun^\otimes(\Qsf(\phi), \Dsf^{\leq 0}(R))^\simeq$.
    Because every object of $\Qsf(\phi)$ is invertible with respect to $\otimes_\Qsf$, the image of every symmetric monoidal functor $\Qsf(\phi) \to \Dsf^{\leq 0}(R)$ must consist of line bundles, hence must lie in $\Perf^{\leq 0}(R)$.
    Thus $\Fun^\otimes(\Qsf(\phi), \Dsf^{\leq 0}(R))^\simeq = \Fun^\otimes(\Qsf(\phi), \Perf^{\leq 0}(R))^\simeq$.
\end{proof}

\subsection{Windows and quivers}

There is a deep relationship between the derived category of $\Yfr$ and that of an open substack $\Xfr \subset \Yfr$ (see e.g.\ \cite{hl_derived}).
In particular, in many cases, there is a nice embedding $\QC(\Xfr) \hookrightarrow \QC(\Yfr)$.
We can formalize one notion of ``niceness'' as follows:

\begin{dfn}
    A \emph{window} is a functor $W: \QC(\Xfr) \hookrightarrow \QC(\Yfr)$ such that:
    \begin{itemize}
        \item $j_\Xfr^* W = \id_{\QC(\Xfr)}$,
        \item $W$ preserves colimits, and
        \item $W$ preserves compact objects (i.e.\ perfect complexes).
    \end{itemize}
    The adjoint functor theorem \HTT{5.5.2.9} implies that $W$ admits a right adjoint $H: \QC(\Yfr) \to \QC(\Xfr)$, which we call the \emph{Hitchcock functor} corresponding to $W$.
\end{dfn}

\begin{rmk}
    Our hypotheses on $W$ guarantee that $W$ can be recovered from the restriction $W|_{\Perf(\Xfr)}: \Perf(\Xfr) \to \Perf(\Yfr)$.
    In the literature, the term ``window'' more frequently refers to this restriction (or to closely related notions, e.g.\ with $\Perf$ replaced by $\Coh$).
    The windows that interest us necessarily send perfect complexes to perfect complexes, so we do not gain anything by replacing $\QC$ by $\Indsf \Coh$.
\end{rmk}

We are particularly interested in windows coming from certain collections of line bundles on $\Xfr$.
For such a collection to determine a window, the line bundles should be induced by weights of $G$ satisfying the following conditions.

\begin{dfn} \label{dfn:transparent}
    A collection of weights $S \subset \weight(G)$ is \emph{transparent} for $\Xfr \subset \Yfr$ if:
    \begin{itemize}
        \item $\{ \Osc_\Xfr(\chi) \}_{\chi \in S}$ compactly generates $\QC(\Xfr)$, and
        \item $\Hom_\Xfr\big(\Osc_\Xfr(\chi_1), \Osc_\Xfr(\chi_2)\big) = \Hom_{\Yfr}\big(\Osc_\Yfr(\chi_1), \Osc_\Yfr(\chi_2)\big)$ for all $\chi_1, \chi_2 \in S$.
    \end{itemize}
\end{dfn}

Because $Y$ is affine and $G$ is reductive, for all $\chi_1, \chi_2 \in \weight(G)$, we have
\[
    \Hom_{\Yfr}^i\big(\Osc_\Yfr(\chi_1), \Osc_\Yfr(\chi_2)\big) = 0 \textrm{ for all } i > 0.
\]
In particular, if $S$ is transparent, then for all $\chi_1, \chi_2 \in S$, we have
\begin{equation} \label{eq:transparent_hom}
    \Hom_\Xfr\big(\Osc_\Xfr(\chi_1), \Osc_\Xfr(\chi_2)\big) = \Hom_{\Yfr}\big(\Osc_\Yfr(\chi_1), \Osc_\Yfr(\chi_2)\big) = R_{\chi_2 - \chi_1}
\end{equation}
using the notation of \cref{rmk:action_grading}. The following is one of the primary examples we should keep in mind.

The notions of ``transparent collection of weights'' and ``full strong exceptional collection of line bundles'' are related but distinct.
The key distinction for our purposes is that transparent collections are defined relative to an embedding $\Xfr \subset \Yfr$, while full strong exceptional collections are defined solely in terms of $\Xfr$.
Transparent collections also exist more frequently than full strong exceptional collections when $\Xfr$ is not proper.
(For example, the entire collection $\weight(G)$ is always transparent for $\Yfr \subset \Yfr$.)
Nevertheless, in many cases, the two notions coincide:

\begin{prop} \label{prop:transparent_fsec}
    Assume that $Y$ is an (affine) underived algebraic variety and $\Hom_\Yfr(\Osc_\Yfr, \Osc_\Yfr) = k[0]$, i.e.\ $\Osc_\Yfr$ is an exceptional object of $\QC(\Yfr)$.
    \begin{enumerate}
        \item If $S \subset \weight(G)$ is transparent for $\Xfr \subset \Yfr$, then the line bundles $\{ \Osc_\Xfr(\chi) \}_{\chi \in S}$ form a full strong exceptional collection in $\Perf(\Xfr)$.
        \item Suppose $Y$ is normal and $\codim_Y (Y \setminus X) \leq 2$.
        If $S \subset \weight(G)$ is such that the line bundles $\{ \Osc_\Xfr(\chi) \}_{\chi \in S}$ form a full strong exceptional collection in $\Perf(\Xfr)$, then $S$ is transparent for $\Xfr \subset \Yfr$.
        \item Suppose $Y = \Spec R$ where $R$ is a unique factorization domain, and suppose $\Pic(\Xfr) = \weight(G)$.
        Then transparent collections of weights for $\Xfr \subset \Yfr$ correspond precisely to full strong exceptional collections in $\Perf(\Xfr)$ via the correspondence $S \mapsto \{ \Osc_\Xfr(\chi) \}_{\chi \in S}$.
    \end{enumerate}
\end{prop}

\begin{proof}
    Throughout we use the equivalence between line bundles on a stack quotient $[Z / G]$ and $G$-equivariant line bundles on $Z$.
    
    (1). By supposition, $\Osc_\Xfr(\chi)$ is exceptional in $\Perf(\Xfr)$ for every $\chi \in \weight(G)$.
    Also, for $\chi \in \weight(G)$, note that only one of $H^0(\fr \Xfr,\Osc_\Xfr(\chi))$ and $H^0(\Xfr, \Osc_\Xfr(-\chi))$ can be non-zero.
    This gives a partial order on $S$.
    Thus, if $S \subset \weight(G)$ is transparent, then $\{\Osc_\Xfr(\chi)\}_{\chi \in S}$ can be ordered to be a full strong exceptional collection. 
    
    (2). It suffices to show that $\Hom_{\Xfr}^0\big( \Osc_\Xfr(\chi_1), \Osc_\Xfr(\chi_2) \big) = \Hom_{\Yfr}^0\big( \Osc_\Yfr(\chi_1), \Osc_\Yfr(\chi_2) \big)$ for all $\chi_1, \chi_2 \in S$.
    In fact, this holds for all $\chi_1, \chi_2 \in \weight(G)$, as we now show.
    
    We may identify $\Hom_\Xfr^0\big( \Osc_\Xfr(\chi_1), \Osc_\Xfr(\chi_2) \big)$ with the subspace of $H^0(X, \Osc_X)$ on which $G$ acts with weight $\chi_2 - \chi_1$.
    Any section of $H^0(X, \Osc_X)$ extends uniquely to a section of $H^0(Y, \Osc_Y)$ by the algebraic Hartogs' lemma (see e.g.\ \stacks{031T}), and the uniqueness implies that $G$ acts on this new section with the same weight.
    Thus we may identify $\Hom_\Xfr^0\big( \Osc_\Xfr(\chi_1),\Osc_\Xfr(\chi_2) \big)$ with the subspace of $H^0(Y, \Osc_Y)$ on which $G$ acts with weight $\chi_2 - \chi_1$.
    This last subspace may itself be identified with $\Hom_{\Yfr}^0\big( \Osc_\Yfr(\chi_1), \Osc_\Yfr(\chi_2) \big)$.

    (3). By hypothesis, every line bundle on $Y$ is (non-equivariantly) trivial.
    It follows that a $G$-equivariant line bundle on $Y$ is determined by the character by which $G$ acts on a (non-equivariant) non-vanishing section.
    In other words, there is a natural surjection $\weight(G) \twoheadrightarrow \Pic(\Yfr)$.
    Applying the restriction map $\Pic(\Yfr) \to \Pic(\Xfr) = \weight(G)$, we see that the identity map $\weight(G) \to \weight(G)$ factors through the surjection $\weight(G) \twoheadrightarrow \Pic(\Yfr)$, and thus $\Pic(\Yfr) = \weight(G)$.
    In particular, we see that the restriction $\Pic(\Yfr) \to \Pic(\Xfr)$ is an isomorphism, i.e.\ every $G$-equivariant line bundle on $X$ extends to a unique $G$-equivariant line bundle on $Y$.
    
    We claim that this implies $\codim_Y (Y \setminus X) \geq 2$.
    Indeed, if this is not true, we can choose a nonempty divisor $D \subset Y$ which is entirely contained in $Y \setminus X$.
    Replacing $D$ by $G \cdot D$, we may assume without loss of generality that $D$ is $G$-invariant.
    Thus $\Osc_Y(D)$ is naturally $G$-equivariant.
    Note that $\Osc_Y(D)$ restricts to the trivial $G$-equivariant line bundle on $X$, so $\Osc_Y(D)$ must be $G$-equivariantly trivial on $Y$.
    That is, there exists $f \in H^0(Y, \Osc_Y)^G$ such that the effective divisor $D$ is linearly equivalent to $V(f)$.
    But the hypothesis $H^0(Y, \Osc_Y)^G = \Hom_\Yfr(\Osc_\Yfr, \Osc_\Yfr) = k$ then implies $D$ must be empty, a contradiction.
    Thus $\codim_Y (Y \setminus X) \geq 2$ and we are in the situation of (2).
\end{proof}

We now introduce the windows arising from transparent collections of weights.
The following is known to experts,\footnote{We thank Daniel Halpern-Leistner for bringing a version of this statement to our attention.} and we state it partially to establish precise notation for later.

\begin{prop} \label{prop:transparent_induces_window}
    Suppose $S$ is transparent for $\Xfr \subset \Yfr$.
    Then:
    \begin{enumerate}
        \item There is a natural equivalence $\QC(\Xfr) = \Dsf(\Qsf_S(\phi)\op)$.
        \item There exists a window $W_S: \QC(\Xfr) \hookrightarrow \QC(\Yfr)$ such that $W_S(\Osc_\Xfr(\chi)) = \Osc_\Yfr(\chi)$ for all $\chi \in S$.
    \end{enumerate}
\end{prop}

\begin{proof}
    (1). Let $\chi_1, \chi_2 \in S$.
    By \cref{eq:transparent_hom}, we have $\Hom_{\Xfr}(\Osc_\Xfr(\chi_1), \Osc_\Xfr(\chi_2)) = \Hom_{\Qsf_S(\phi)}(\chi_1, \chi_2)$.
    Thus the full subcategory of $\QC(\Xfr)$ with objects $\{ \Osc_\Xfr(\chi) \}_{\chi \in S}$ is equivalent to $\Qsf_S(\phi)$.
    This full subcategory compactly generates $\QC(\Xfr)$, yielding the claim.

    (2). The inclusion $i_S: \Qsf_S(\phi) \hookrightarrow \Qsf(\phi)$ left Kan extends to a functor $(i_S)_!: \Dsf(\Qsf_S(\phi)\op) \to \Dsf(\Qsf(\phi)\op)$.
    Note that $(i_S)_!$ preserves colimits and compact objects by definition.
    Applying part (1) to $\weight(G)$ (which is transparent for $\Yfr \subset \Yfr$), we get an equivalence $\QC(\Yfr) = \Dsf(\Qsf(\phi)\op)$.
    Let $W_S: \QC(\Xfr) \to \QC(\Yfr)$ be the functor corresponding to $(i_S)_!$ via the above equivalences.
    Chasing the definitions shows that $W_S(\Osc_\Xfr(\chi)) = \Osc_\Yfr(\chi)$ for all $\chi \in S$.
    In particular, since $j_\Xfr^* W_S(\Osc_\Xfr(\chi)) = \Osc_\Yfr(\chi)$ and $\{ \Osc_\Xfr(\chi) \}_{\chi \in S}$ generates $\QC(\Xfr)$, we see that $j_\Xfr^* W_S = \id_{\QC(\Xfr)}$, so $W_S$ is a window.
\end{proof}

\begin{rmk} \label{rmk:compact_representations}
    If $S \subset \weight(G)$ is any collection of weights, note that the compact objects of $\Dsf(\Qsf_S(\phi)\op)$ are those which can be written as a finite colimit (possibly with shifts) of the generators $\big\{ \Osc_\Yfr(\chi) \big\}_{\chi \in S}$.
    In general, we cannot expect $\Dsf(\Qsf_S(\phi)\op)^\omega = \Fun_k(\Qsf_S(\phi)\op, \Perf(k))$, though this is true if $S$ arises from a finite full strong exceptional collection as in \cref{prop:transparent_fsec}.
\end{rmk}

\begin{notn}
    If $S$ is a transparent collection of weights for $\Xfr \subset \Yfr$, we fix the following notation:
    \begin{itemize}
        \item $W_S$ is the window of \cref{prop:transparent_induces_window}.
        \item $H_S$ is the corresponding Hitchcock functor (i.e.\ right adjoint of $W_S$).
        \item $\Ksc_S$ is the Fourier-Mukai kernel of $W_S$.
    \end{itemize}
\end{notn}

By the construction in \cref{prop:transparent_induces_window}, we see that the window $W_S$ may be identified with the left Kan extension functor $i_{S!}: \Dsf(\Qsf_S(\phi)\op) \to \Dsf(\Qsf(\phi)\op)$.
The right adjoint $H_S$ may therefore be identified with the ``restriction of quiver representations'' functor $i_S^*: \Dsf(\Qsf(\phi)\op) \to \Dsf(\Qsf_S(\phi)\op)$. 

\begin{ex}
    Let $\Xfr = \PP^n \subset \Yfr = [\AA^{n+1} / \GG_m]$, and let $S = \{ 0, \dots, n \} \subset \ZZ = \weight(\GG_m)$ (so $S$ is transparent for $\Xfr \subset \Yfr$).
    Direct computations using the ``restriction of quiver representations'' description of the Hitchcock functor $H_S$ show that
    \[
    H_S\big(\Osc_{[\AA^{n+1}/\GG_m]}(\ell)\big) = \begin{cases}
        \Osc_{[\AA^{n+1}/\GG_m]}(\ell) & \ell \geq 0 \\
        0 & \ell < 0
    \end{cases}
    \]
    for all $\ell \in \ZZ$.
    In particular, $H_S$ is distinct from the geometric pullback functor $j^*$ as functors from $\Dsf(\Qsf(\phi)\op)$ to $\Dsf(\Qsf_S(\phi)\op)$, although the two functors agree on the image of $W_S$.
\end{ex}

In general, it is difficult to give a geometric description of the windows $W_S$ or the Hitchcock functors $H_S$ in terms of Fourier-Mukai transforms.
However, in certain favorable cases, the windows $W_S$ arise via push-pull along compactifications of the diagonal of $\Xfr$ as in \cite{ballard2017kernels} (see \cref{rmk:toric_dm}).
We show this claim for a class of toric varieties in \cref{prop:toric_variety_geometric}.
\cref{prop:hitchcock_via_opposite} allows us to leverage this to provide geometric descriptions of the Hitchcock functors for the opposite collections of weights.
We use this to understand certain windows and Hitchcock functors for the inclusion $\PP^n \subset [\AA^{n+1} / \GG_m]$ in \cref{ex:Pn_window_push_pull} and \cref{ex:Pn_hitchcock}.

\subsection{Windows and resolutions of the diagonal} \label{sub:window_resolution}

A standard method (introduced in \cite{beilinson1978coherent}) for showing an exceptional collection in $\Perf(\Xfr)$ is full is to show that the collection can be used to produce a resolution of the diagonal of $\Xfr$.
When the exceptional collection in question corresponds to a collection of weights $S \subset \weight(G)$ which is transparent for $\Xfr \subset \Yfr$, this resolution of the diagonal may be used to concretely understand the window $W_S$.

\begin{prop} \label{prop:diagonal_window}
    Suppose $S \subset \weight(G)$ is such that $\Hom_\Xfr(\Osc_\Xfr(\chi_1), \Osc_\Xfr(\chi_2)) = \Hom_\Yfr(\Osc_\Yfr(\chi_1), \Osc_\Yfr(\chi_2))$ for all $\chi_1, \chi_2 \in S$.
    Then $S$ is transparent for $\Xfr \subset \Yfr$ if and only if there is an expression 
    \[
        \Delta_{\Xfr*} \Osc_\Xfr = \colim_{i \in I} \Asc_i \boxtimes \Osc_\Xfr(\chi_i)
    \]
    with all $\chi_i \in S$ and all $\Asc_i \in \QC(\Xfr)$.
    In this case, $\Ksc_S = \colim_{i \in I} \Asc_i \boxtimes \Osc_\Yfr(\chi_i) \in \QC(\Xfr \times \Yfr)$.
\end{prop}

\begin{proof}
    With the given hypotheses, transparency of $S$ is equivalent to the claim that $\{ \Osc_\Xfr(\chi) \}_{\chi \in S}$ generates $\QC(\Xfr)$.
    (Compactness of the sheaves $\{ \Osc_\Xfr(\chi) \}_{\chi \in S}$ is obvious because $\Xfr$ is perfect.)
    The equivalence in the statement of the Proposition now follows from \cref{lem:identity_resolution_diagonal}.
    The computation of $\Ksc_S$ is \cref{prop:fm_resolution_diagonal}.
\end{proof}

\begin{ex} \label{ex:Pn_window_complex}
    Let $\Xfr = \PP^n \subset [\AA^{n+1} / \GG_m]$.
    The Beilinson complex
    \[
    \begin{tikzcd}
        \Omega^n_{\PP^n}(n) \boxtimes \Osc_{\PP^n}(-n) \rar & \Omega^{n-1}_{\PP^n}(n-1) \boxtimes \Osc_{\PP^n}(-n+1) \rar & \dots \rar & \Osc_{\PP^n} \boxtimes \Osc_{\PP^n}
    \end{tikzcd}
    \]
    may be understood as a colimit over all of its terms (building the complex up one term at a time by a sequence of iterated mapping cones / cofibers).
    This complex is quasi-isomorphic to $\Delta_{\PP^n *} \Osc_{\PP^n}$, so by \cref{prop:diagonal_window}, we see that $S = \{ -n, \dots, 0 \}$ is transparent for $\PP^n \subset [\AA^{n+1} / \GG_m]$.
    The complex $\Ksc_S$ is given by
    \[
    \begin{tikzcd}
        \Omega^n_{\PP^n}(n) \boxtimes \Osc_{[\AA^{n+1}/\GG_m]}(-n) \rar & \dots \rar & \Osc_{\PP^n} \boxtimes \Osc_{[\AA^{n+1}/\GG_m]}.
    \end{tikzcd}
    \]
\end{ex}

When the resolution of the diagonal of $\Xfr$ extends to a resolution of the pushforward of the structure sheaf of a perfect stack over $\Xfr \times \Yfr$, we may obtain a more geometric description of $W_S$:

\begin{prop} \label{prop:check_geometric}
    Suppose that there exist a collection of weights $S \subset \weight(G)$, a perfect stack $\Wc_S$, and a morphism $q = (q_1, q_2): \Wc_S \to \Xfr \times \Yfr$ such that:
    \begin{itemize}
        \item $\Hom_\Xfr(\Osc_\Xfr(\chi_1), \Osc_\Xfr(\chi_2)) = \Hom_\Yfr(\Osc_\Yfr(\chi_1), \Osc_\Yfr(\chi_2))$ for all $\chi_1, \chi_2 \in S$.
        \item There is an expression $q_* \Osc_{\Wc_S} = \colim_{i \in I} \Asc_i \boxtimes \Osc_\Yfr(\chi_i)$ with all $\chi_i \in S$.
        \item There is a Cartesian square
        \[
            \begin{tikzcd}
                \Xfr \rar \dar["\Delta_{\Xfr}"] & \Wc_S \dar["q"] \\
                \Xfr \times \Xfr \rar["{(\id_\Xfr, j)}"] & \Xfr \times \Yfr.
            \end{tikzcd}
        \]
    \end{itemize}
    Then $S$ is transparent for $\Xfr \subset \Yfr$, and $W_S = q_{2*} q_1^*$.
\end{prop}

\begin{proof}
    By base change, we have
    \[
        \Delta_{\Xfr *} \Osc_\Xfr = (\id_\Xfr, j)^* q_* \Osc_{\Wc_S} = \colim_{i \in I} \Asc_i \boxtimes \Osc_\Xfr(\chi_i) \in \QC(\Xfr \times \Xfr).
    \]
    Thus the hypotheses of \cref{prop:diagonal_window} are satisfied, $S$ is transparent for $\Xfr \subset \Yfr$, and $W_S$ is given by the Fourier-Mukai transform with kernel $q_* \Osc_{\Wc_S}$.
    This Fourier-Mukai transform is exactly $q_{2*} q_1^*$.
\end{proof}

\begin{ex} \label{ex:Pn_window_push_pull}
    In the situation of \cref{ex:Pn_window_complex}, the complex $\Ksc_S$ has cohomology sheaves $\Hsc^i(\Ksc_S)$ concentrated in degree $0$.
    Furthermore, the map $\Osc_{\PP^n} \boxtimes \Osc_{[\AA^{n+1} / \GG_m]} \to \Ksc_S$ (induced from the brutal truncation) is a surjection on $\Hsc^0$.
    That is, $\Ksc_S$ is a quotient of the structure sheaf in the abelian category of quasicoherent sheaves.
    Thus $\Ksc_S$ is naturally a commutative algebra.
    Let $\Wc_S = \rSpec_{\PP^n \times [\AA^{n+1} / \GG_m]} \Ksc_S$, and let $q = (q_1, q_2): \Wc \to \PP^n \times [\AA^{n+1} / \GG_m]$ be the natural map.
    The conditions of \cref{prop:check_geometric} hold automatically, so we must have $W_S = q_{2*} q_1^*$.
\end{ex}

\subsection{Hitchcock functors via opposite collections of weights}

We may use the results of \cref{sub:window_resolution} to obtain geometric descriptions of the windows $W_S$ in some nice cases.
For our applications, we will also want to understand the Hitchcock functors $H_S$ and their relationship to the aforementioned windows.
It turns out that the Fourier-Mukai kernel of $H_S$ agrees with the Fourier-Mukai kernel of a window corresponding to a \emph{different} transparent collection of weights!

Let $S \subset \weight(G)$ be an arbitrary collection of weights.
The \emph{opposite collection} to $S$ is $-S = \bset{-\chi}{\chi \in S}$.
By taking duals of line bundles, we see that $\Qsf_S(\phi)\op \simeq \Qsf_{-S}(\phi)$.

\begin{prop}
    Let $S$ be transparent for $\Xfr \subset \Yfr$.
    Then $-S$ is also transparent for $\Xfr \subset \Yfr$.
    Furthermore, for all $\Fsc \in \Perf(\Xfr)$, there is a natural equivalence $W_{-S}(\Fsc) = W_S(\Fsc^\vee)^\vee$.
\end{prop}

\begin{proof}
    Let $\Fsc \in \Perf(\Xfr)$.
    By hypothesis, $\Fsc^\vee$ may be obtained from the objects $\Osc_\Xfr(\chi)$ for $\chi \in S$ via a finite sequence of taking shifts, (co)fibers, and direct summands.
    Dualizing, we see that $\Fsc$ can be obtained from the objects $\Osc_\Xfr(-\chi)$ for $\chi \in S$ via a finite sequence of taking shifts, (co)fibers, and direct summands.
    Because $\Perf(\Xfr)$ compactly generates $\QC(\Xfr)$, the collection $\{ \Osc_\Xfr(-\chi) \}_{\chi \in S}$ also compactly generates $\QC(\Xfr)$.

    For all $-\chi_1, -\chi_2 \in S$ and all $i > 0$, we have 
    \[
        \Hom_\Xfr^i\big(\Osc_\Xfr(-\chi_1), \Osc_\Xfr(-\chi_2)\big) = \Hom_\Xfr^i\big(\Osc_\Xfr(\chi_2), \Osc_\Xfr(\chi_1)\big) = 0
    \]
    by transparency of $S$.
    Thus $-S$ is transparent.

    To see the claim about $W_{-S}(\Fsc)$, by writing $\Fsc$ as a canonical colimit of objects $\Osc_\Xfr(-\chi)$, we reduce to showing that $W_{-S}\big(\Osc_\Xfr(-\chi)\big) = W_S(\Osc_\Xfr(-\chi)^\vee)^\vee$ for all $\chi \in S$.
    But this is obvious, as both $W_{-S}\big(\Osc_\Xfr(-\chi)\big)$ and $W_S(\Osc_\Xfr(-\chi)^\vee)^\vee$ are naturally equivalent to $\Osc_\Yfr(-\chi)$ by definition.
\end{proof}

We may use knowledge of $W_{-S}$ to compute the Hitchcock functor $H_S$:

\begin{prop} \label{prop:hitchcock_via_opposite}
    Let $S$ be transparent for $\Xfr \subset \Yfr$.
    Then $\Ksc_{-S}$ is the Fourier-Mukai kernel of the Hitchcock functor $H_S: \QC(\Yfr) \to \QC(\Xfr)$.
\end{prop}

\begin{proof}
    Let $\pi_\Yfr: \Yfr \times \Xfr \to \Yfr$ and $\pi_\Xfr: \Yfr \times \Xfr \to \Xfr$ be the natural projections.
    We need to show that
    \[
        \Hom_\Xfr(\Fsc, H_S \Gsc) = \Hom_\Xfr\big(\Fsc, \pi_{\Xfr*}(\pi_\Yfr^* \Gsc \otimes \Ksc_{-S})\big)
    \]
    for all $\Fsc \in \QC(\Xfr)$ and all $\Gsc \in \QC(\Yfr)$.
    Because $\Perf(\Xfr)$ compactly generates $\QC(\Xfr)$, $\Perf(\Yfr)$ compactly generates $\QC(\Yfr)$, and $H_S$ preserves colimits, it suffices to take $\Fsc \in \Perf(\Xfr)$ and $\Gsc \in \Perf(\Yfr)$.
    Then we may compute
    \begin{align*}
        \Hom_\Xfr(\Fsc, H_S \Gsc) &= \Hom_\Yfr(W_S \Fsc, \Gsc) \\
        &= \Hom_\Yfr\big((W_{-S}(\Fsc^\vee))^\vee, \Gsc\big) \\
        &= \Hom_\Yfr\big(\Gsc^\vee, W_{-S} \Fsc^\vee\big) \\
        &= \Hom_\Yfr\big(\Gsc^\vee, \pi_{\Yfr*}(\pi_\Xfr^* \Fsc^\vee \otimes \Ksc_{-S})\big) \\
        &= \Hom_{\Xfr \times \Yfr}\big(\pi_\Yfr^* \Gsc^\vee, \pi_\Xfr^* \Fsc^\vee \otimes \Ksc_{-S}\big) \\
        &= \Hom_{\Xfr \times \Yfr}\big(\pi_\Xfr^* \Fsc, \pi_\Yfr^* \Gsc \otimes \Ksc_{-S}\big) \\
        &= \Hom_\Xfr\big(\Fsc, \pi_{\Xfr*}(\pi_\Yfr^* \Gsc \otimes \Ksc_{-S})\big). \qedhere
    \end{align*}
\end{proof}

\begin{cor}
    Let $S$ be transparent for $\Xfr \subset \Yfr$, and suppose there exists a perfect stack $\Wc_{-S}$ and a morphism $q = (q_1, q_2): \Wc_{-S} \to \Xfr \times \Yfr$ such that $W_{-S} = q_{2*} q_1^*$.
    Then $H_S = q_{1*} q_2^*$.
\end{cor}

\begin{proof}
    Take $\Ksc_{-S} = q_* \Osc_{\Wc_{-S}}$ in \cref{prop:hitchcock_via_opposite}.
\end{proof}
\begin{ex} \label{ex:Pn_hitchcock}
    Let $\Xfr = \PP^n \subset [\AA^n / \GG_m]$, and let $q: \Wc \to \PP^n \times [\AA^n / \GG_m]$ be as in \cref{ex:Pn_window_push_pull}.
    For $S = \{ 0, \dots, n \}$, we see that $H_S = q_{1*} q_2^*$.
\end{ex}

\subsection{The toric case} \label{sub:toric_examples}

Our key examples of transparent collections of line bundles will arise from smooth toric stacks.
For our purposes, a \emph{smooth toric stack} is an open substack $\Xfr = [X / G]$ of a stack quotient $[\AA^n / G]$, where the action $G \curvearrowright \AA^n$ is induced by a homomorphism $G \to \GG_m^n$ and the diagonal action $\GG_m^n \curvearrowright \AA^n$.
As above, we write $\phi: [\AA^n / G] \to BG$ for the structure map.

We begin by introducing a combinatorial description of the total weight quiver $\Qsf(\phi)$.

\begin{dfn} \label{dfn:prelinear_quiver}
    Let $\phi: G \to \GG_m^n$ be a homomorphism.
    We may define a (non-enriched) symmetric monoidal discrete category $\Qsf^\pre(\phi)$, the \emph{prelinear total weight quiver}, where:
    \begin{itemize}
        \item $\ob \Qsf^\pre(\phi) = \weight(G)$.
        \item $\Map_{\Qsf^\pre(\phi)}(\chi_1, \chi_2)$ consists of monomials of degree $\chi_2 - \chi_1$ in the variables $x_1, \dots, x_n$.
        \item Composition is given by multiplication of monomials.
        \item Multiplication of objects is the group operation in $\weight(G)$.
        \item Multiplication of morphisms is multiplication of monomials.
    \end{itemize}
    For $S \subset \weight(G)$, we define the \emph{prelinear partial weight quiver} $\Qsf_S^\pre(\phi)$ as the full subcategory of $\Qsf^\pre(\phi)$ with objects $S$.
\end{dfn}

For $S \subset \weight(G)$, the weight quiver $\Qsf_S(\phi)$ agrees with the $k$-linearization of the prelinear weight quiver $\Qsf_S^\pre(\phi)$.

\begin{ex}
    For $\Xfr = \PP^n \subset [\AA^{n+1} / \GG_m]$, the category $\Qsf^\pre(\phi)$ is the infinite Beilinson quiver $Q_{n,\infty}$ of the Introduction (depicted in \cref{fig:infinite_beilinson_quiver}).
    The collection $S = \{ 0, \dots, n \}$ is transparent for $\PP^n \subset [\AA^{n+1} / \GG_m]$, and the category $\Qsf_S^\pre(\phi)$ is the Beilinson quiver $\Qsf_n$ depicted in \cref{fig:beilinson_quiver}.
\end{ex}

\begin{ex}
    Let $\Xfr = F_2$, the Hirzebruch surface of type $2$.
    We may realize $\Xfr$ as an open subset of the quotient $[\AA^4 / \GG_m^2]$, where $\GG_m^2$ acts on $\AA^4$ by
    \[
        (g_1, g_2) \cdot (y_1, y_2, y_3, y_4) = (g_1 y_1, g_1^{-2} g_2 y_2, g_1 y_3, g_2 y_4).
    \]
    (See e.g.\ \cite[Examples 14.2.17 and 14.2.20]{cls_toric}.)
    The category $\Qsf^\pre(\phi)$ is the infinite quiver
    
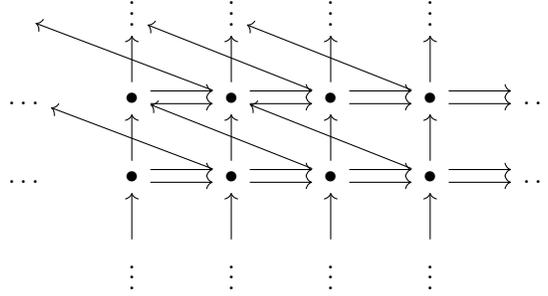
\[
    \begin{tikzcd}
        {} & \vdots & \vdots & \vdots & \vdots & {} \\
        \hdots & \bullet \ar[r, shift left] \ar[r, shift right] \ar[u] & \bullet \ar[ull] \ar[r, shift left] \ar[r, shift right] \ar[u] & \bullet \ar[ull] \ar[r, shift left] \ar[r, shift right] \ar[u] & \bullet \ar[ull] \ar[r, shift left] \ar[r, shift right] \ar[u] & \hdots \\
        \hdots & \bullet \ar[r, shift left] \ar[r, shift right] \ar[u] & \bullet \ar[ull] \ar[r, shift left] \ar[r, shift right] \ar[u] & \bullet \ar[ull] \ar[r, shift left] \ar[r, shift right] \ar[u] & \bullet \ar[ull] \ar[r, shift left] \ar[r, shift right] \ar[u] & \hdots \\
        {} & \vdots \uar & \vdots \uar & \vdots \uar & \vdots \uar & {}
    \end{tikzcd}
    \]
    where:
    \begin{itemize}
        \item the horizontal arrows are given by the monomials $y_1$ and $y_3$,
        \item the vertical arrows are given by the monomials $y_4$,
        \item the slanted diagonal arrows are given by the monomials $y_3$, and
        \item there are (implicit) relations $y_i y_j = y_j y_i$.
    \end{itemize}
    Identify $\weight(\GG_m^2) = \ZZ^2$ via the isomorphism corresponding to the $x$ and $y$ axes in our depiction of $\Qsf^\pre(\phi)$.
    The collection $S = \{ (0, 0), (1, 0), (-2, 1), (-1, 1)$ is transparent for $\Xfr \subset [\AA^4 / \GG_m^2]$ by the results of \cite[\S 6]{king_conj}.
\end{ex}

Let us say that the open immersion $j_\Xfr: \Xfr \subset [\AA^n / G]$ is \emph{decent} if $\weight(G) = \Pic(\Xfr)$.
If $\Xfr$ is a smooth toric stack and $\Xfr \hookrightarrow [\AA^n / G]$ is a decent open immersion, then any full strong exceptional collection of line bundles on $\Xfr$ corresponds to a transparent collection of weights $S$ for $\Xfr \subset [\AA^n / G]$ by \cref{prop:transparent_fsec}.
Thus we obtain a combinatorial description of $\QC(\Xfr)$:
\[
    \QC(\Xfr) = \Fun\big(\Qsf_S^\pre(\phi), \Dsf(k)\big).
\]
This description applies to many key examples.

\begin{ex} \label{ex:toric_cox}
    If $\Xfr$ is a smooth toric variety without torus factors, the Cox construction (see \cite[Chapter 5]{cls_toric}) gives a decent open immersion $\Xfr \hookrightarrow [\AA^n / G]$, where $n$ is the number of rays in the fan of $\Xfr$.
\end{ex}

\begin{ex}
    Let $d_0, \dots, d_n$ be positive integers, and let $\PP(d_0, \dots, d_n)$ be the weighted projective stack with weights $d_0, \dots, d_n$.
    Assume $n \geq 1$.\footnote{When $n = 0$, we have $\PP(d_0) = B (\ZZ / {d_0})$, and one can check the existence of a transparent collection of line bundles ``by hand.''}
    Then the standard open immersion $\PP(d_0, \dots, d_n) \hookrightarrow [\AA^{n+1} / \GG_m]$ is decent.
\end{ex}

In the situation of \cref{ex:toric_cox}, the window corresponding to the \emph{Bondal-Thomsen collection} of weights has a particularly nice description whenever it exists.
We first recall some relevant definitions.

\begin{dfn}[{\cite{bondal2006derived}}] \label{dfn:bondal_thomsen}
    Let $\Xfr$ be a smooth complete toric variety, and let $j: \Xfr \hookrightarrow \Yfr = [\AA^n / G]$ be the decent open immersion of \cref{ex:toric_cox}.
    Let $\phi^*: \ZZ^n \to \weight(G)$ be the natural map induced by the action $G \curvearrowright \AA^n$.
    Let $M = \ker \phi^*$, and let $M_\RR = M \otimes_\ZZ \RR$.
    Define a (discontinuous) \emph{anti-Bondal-Ruan map} $F_\phi: M_\RR / M \to \weight(G)$ (where $M_\RR / M$ is viewed as a real torus) by
    \[
        F_\phi\big(\sum_i a_i e_i + M\big) = \phi^*\big(\sum_i \lfloor a_i \rfloor \mathbf{e}_i\big)
    \]
    The \emph{anti-Bondal-Thomsen collection} (or \emph{anti-BT collection}) is $-\Theta_\Xfr = \im F_\phi \subset \weight(G)$.
    The \emph{Bondal-Thomsen collection} is $\Theta_\Xfr = -(-\Theta_\Xfr)$.
    Following \cite[Definition 5.3]{favero2025homotopy}, we say that $\Xfr$ is of \emph{Bondal-Ruan type} if $\Theta_\Xfr$ (equivalently, $-\Theta_\Xfr$) is transparent for $\Xfr \subset [\AA^n / G]$.
\end{dfn}

\begin{rmk} \label{rmk:dual_conventions}
    We use the term ``anti'' above to stress that we are prioritizing different conventions from those of \cite{bondal2006derived} and subsequent works.
    This is motivated by \cref{prop:hitchcock_via_opposite}: the window corresponding to the Bondal-Thomsen collection has a nice geometric description, so the Hitchcock functor corresponding to the anti-BT collection has a corresponding description.
    Because of the relevance of Hitchcock functors to the results of \cref{sec:ec}, we prefer to use the anti-BT collection.
\end{rmk}

\begin{rmk}
    The generalization of \cref{dfn:bondal_thomsen} to toric Deligne-Mumford stacks is straightforward and profitable; see e.g.\ \cite[\S 2.3]{hanlon2024resolutions}.
    We avoid doing so here due to the difficulties sketched in \cref{rmk:toric_dm} below.
\end{rmk}

Intersection-theoretic sufficient conditions for a smooth complete toric variety to be of Bondal-Ruan type are given in \cite{bondal2006derived}.
(See also \cite[Proposition 5.18]{favero2025homotopy} for a more detailed proof.)
These conditions hold for many examples of interest, e.g.\ all but two smooth toric Fano threefolds.

\begin{ex}
    For $\Xfr = \PP^n$, the Bondal-Thomsen collection of weights is $\{-n, \dots, 0 \} \subset \ZZ = \weight(\GG_m)$.
    This is transparent for $\PP^n \subset [\AA^{n+1} / \GG_m]$, so $\PP^n$ is of Bondal-Ruan type.
\end{ex}

One should also note that many smooth toric varieties (e.g.\ Hirzebruch surfaces) are not of Bondal-Ruan type but still admit transparent collections of weights.
In such cases, we can construct corresponding windows using \cref{prop:transparent_induces_window}, but a geometric interpretation of such windows remains elusive.
In the Bondal-Ruan case, the situation is much nicer:

\begin{prop} \label{prop:toric_variety_geometric}
    Let $\Xfr$ be a smooth complete toric variety of Bondal-Ruan type, and fix notation as in \cref{dfn:bondal_thomsen}.
    Let $\Wc$ be the scheme-theoretic closure (taken in the category of relative schemes over $\Xfr  \times \Yfr$) of the diagonal of $\Xfr$ in $\Xfr \times \Yfr$.\footnote{Alternatively, note that the map $(\id_\Xfr, j): \Xfr \to \Xfr \times \Yfr$ is a locally closed immersion, so its (relative) affinization $\Wc = \rSpec_{\Xfr \times \Yfr} \tau^{\leq 0} (\id_\Xfr, j)_* \Osc_\Xfr$ is a closed substack of $\Xfr \times \Yfr$.}
    Let $q = (q_1, q_2): \Wc \to \Xfr \times \Yfr$ be the natural closed immersion.
    Then $W_{\Theta_\Xfr} = q_{2*} q_1^*$.
\end{prop}

\begin{proof}
    We need only check that $q$ satisfies the hypotheses of \cref{prop:check_geometric} for $S = \Theta_\Xfr$.
    The first condition follows by the definition of Bondal-Ruan type.
    The second condition follows from \cite[Theorem A]{hanlon2024resolutions}, as
    the structure sheaf $q_* \Osc_\Wc$ admits a resolution by direct sums of line bundles $\Osc_\Xfr(\chi_1) \boxtimes \Osc_\Yfr(\chi_2)$ with $\chi_1, \chi_2 \in \Theta_\Xfr$.
    The third condition holds by the definition of $\Wc$.
\end{proof}

In the case where $\Xfr$ is $\PP^n$ (or more generally a weighted projective stack), we may understand the stack $\Wc$ and the push-pull explicitly using blowups. 
This seems to be known to experts (see e.g.\ \cite{ballard2017kernels}), but we include details for completeness.

\begin{ex} \label{ex:wps}
    Let $\Xfr = \PP(d_0, \dots, d_n)$ be a weighted projective stack.
    Consider the weighted stacky blow-up $\sBl^w_0 \AA^{n+1} \to \AA^{n+1}$ of $\AA^{n+1}$ at the origin with weight $w = (d_0,\dots,d_n)$ in the sense of \cite{quek2021weighted}.
    Unpacking the definitions, the stack $\sBl^w_0 \AA^{n+1}$ is equivalent to the (stacky) relative affine spectrum $\rSpec_{\Xfr} \Sym^\bullet_{\Osc_\Xfr} \Osc_{\Xfr}(1)$ (i.e., the (stacky) total space of $\Osc_\Xfr(-1)$) since they are both equivalent to the quotient stack 
    \[
        [(\Spec k[x_0,\dots,x_n][u] \setminus V(x_0,\dots,x_n))/\GG_m]
    \]
    where the $\GG_m$-action on $\Spec k[x_0,\dots,x_n][u]$ corresponds to the $\ZZ$-grading on $k[x_0,\dots,x_n][u]$ given by $\deg x_i = d_i$ and $\deg u = -1$.
    Now, consider another $\GG_m$-action on $\Spec k[x_0,\dots,x_n][u]$ corresponding to the $\ZZ$-grading given by $\deg x_i = 0$ and $\deg u = 1$.
    Combining these two $\GG_m$-actions, we have a $\GG_m^2$-action on $\Spec k[x_0,\dots,x_n][u] \setminus V(x_0,\dots,x_n)$ and we set
    \[
    \Zfr = [(\Spec k[x_0,\dots,x_n][u] \setminus V(x_0,\dots,x_n))/\GG_m^2].
    \]
    Geometrically, the second $\GG_m$-action corresponds to the $\GG_m$-action on $\rSpec_{\Xfr} \Sym^\bullet_{\Osc_\Xfr} \Osc_{\Xfr}(1) \simeq \sBl^w_0 \AA^{n+1}$ given by scaling fibers and we have
    \[
    \Zfr \simeq [\rSpec_{\Xfr} \Sym^\bullet_{\Osc_\Xfr} \Osc_{\Xfr}(1)/\GG_m] \simeq [\sBl^w_0 \AA^{n+1}/\GG_m] 
    \]
    in the sense of \cite[Remark 2.4]{romagny2005group}.
    Now, note that if we consider the diagonal $\GG_m$-action of weight $w$ on $\AA^{n+1}$, then the stacky weighted blow-up $\sBl^w_0 \AA^{n+1} \to \AA^{n+1}$ is $\GG_m$-equivariant. Thus, we obtain a morphism 
    \[
    b: \Zfr \simeq [\sBl^w_0 \AA^{n+1}/\GG_m] \to [\AA^{n+1}/\GG_m] = \Yfr. 
    \]
    Algebraically, this reduces to saying that the ring homomorphism 
    \begin{align*}
        k[x_0,\dots,x_n] &\to k[x_0,\dots,x_n][u] \\
        x_i &\mapsto x_i u^{d_i}
    \end{align*}
    respects $\ZZ^2$-gradings given by $\deg x_i = (d_i,0)$ on $k[x_0,\dots,x_n]$ and by $\deg x_i = (0,d_i)$ and $\deg u = (1,-1)$ on $k[x_0,\dots,x_n][u]$, respectively (cf. \cite[Section 5.2]{quek2021weighted} and \cite{ballard2017kernels} for details).
    
    Similarly, the projection $\rSpec_{\Xfr} \Sym^\bullet_{\Osc_\Xfr} \Osc_{\Xfr}(1) \to \Xfr$ is $\GG_m$-equivariant with respect to the trivial $\GG_m$-action on $\Xfr$, so we obtain a morphism 
    \[
    \pi: \Zfr \simeq [\rSpec_{\Xfr} \Sym^\bullet_{\Osc_\Xfr} \Osc_{\Xfr}(1)/\GG_m] \to \Xfr \times B\GG_m \to \Xfr.
    \] 
    In this case, the Bondal-Thomsen collection $\Theta_\Xfr = \big\{ 1 - \sum_i d_i, \dots, 0 \big\}$ forms a full strong exceptional collection in $\PP(d_0, \dots, d_n)$.
    By comparing with the proof of \cite[Proposition 4.1.5]{ballard2017kernels} (or by direct verification), we can see $W_{\Theta_\Xfr} \simeq b_*\pi^*$.
\end{ex}

\begin{rmk} \label{rmk:toric_dm}
    It is not clear how one would obtain the description of $W_{\Theta_\Xfr}$ in \cref{ex:wps} from the results of \cite{hanlon2024resolutions} (unless e.g. $(d_0,\dots,d_n) = (1,\dots,1)$), as the diagonal morphism of $\PP(d_0, \dots, d_n)$ is finite but not a closed immersion in general.
    One would hope that a refined approach can be used to extend the virtual resolutions of \cite{hanlon2024resolutions} to virtual resolutions of diagonals of separated smooth toric Deligne-Mumford stacks.
    In this case, our results for smooth complete toric varieties of Bondal-Ruan type may be extended to separated smooth toric Deligne-Mumford stacks of Bondal-Ruan type.
\end{rmk}

\subsection{Quiver-theoretic examples}

We may use the approach of \cite{bergman2008moduli} to view a full strong exceptional collection of line bundles on a general smooth projective toric variety as transparent for a certain open immersion.
This open immersion is constructed in terms of moduli spaces of quiver representations.

Let $\Qsf$ be a $k$-linear category such that:
\begin{itemize}
    \item $\Qsf$ has finitely many objects.
    \item All $\Hom$ objects in $\Qsf$ are finite-dimensional vector spaces concentrated in degree zero.
    \item There is a total order on $\ob \Qsf$ such that $\Hom_\Qsf(q_1, q_2) = 0$ only if $q_1 \leq q_2$ 
    \item $\Hom_\Qsf(q, q) = k$ for all $q \in \Qsf$.
\end{itemize}

Let $Y_\Qsf$ be the parameter space of representations of $\Qsf$ with dimension vector $(1, \dots, 1)$, i.e.\ 
\[
    Y \subset \prod_{q_1, q_2 \in \ob \Qsf} \AA\Big(\Hom_k\big(\Hom_\Qsf(q_1, q_2), \Hom_k(k, k)\big)\Big)
\]
is the subscheme defined by the composition relations in $\Qsf$.
Let $G_\Qsf = \Big(\prod_{q \in \ob \Qsf} \GG_m\Big) \Big/ \GG_m^{\pi_0(\Qsf)}$ act on $Y$ by conjugation, where $\pi_0(\Qsf)$ is the set of connected components of $\Qsf$.
Then $\Yfr_\Qsf = [Y_\Qsf / G_\Qsf]$ is the \emph{rigidified moduli stack of vertexwise one-dimensional representations} of $\Qsf$.

If $\Xfr$ is a smooth projective variety and $\{ \Lsc_1, \dots, \Lsc_n \}$ is a full strong exceptional collection of line bundles on $\Xfr$, we may define a small $k$-linear category $\Qsf_\Lsc$ as the opposite of the full subcategory of $\QC(\Xfr)$ consisting of the objects $\{ \Lsc_1, \dots, \Lsc_n \}$.
Let us identify the objects of $\Qsf_\Lsc$ with $\{ 1, \dots, n \}$.
There is a tautological map $j: \Xfr \to \Yfr_{\Qsf_\Lsc}$ defined by sending a point $u \in \Xfr$ to the representation $i \mapsto \Hom_\Xfr(\Lsc_i, \Osc)$.
By \cite[proof of Theorem 2.4]{bergman2008moduli}, the map $j$ is an open immersion.
The vector bundles $\Lsc_i$ are all pullbacks $j^* \Osc_{\Yfr_{\Qsf_\Lsc}}(\chi_i)$ for some $\chi_i \in \weight(G)$, and the collection $S = \{ \chi_1, \dots, \chi_n \} \subset \weight(G)$ is transparent for $\Xfr \subset \Yfr$.

\begin{ex}
    By \cite[Theorem 5.14]{hille2011exceptional}, if $\Xfr$ is a del Pezzo surface with $\rk \Pic \Xfr \leq 7$, then $\Xfr$ admits a full strong exceptional collection of line bundles.
    We do not know when the corresponding $k$-linear categories $\Qsf_\Lsc$ admit ``quiver tensor products.''
\end{ex}

The quiver-theoretic presentations of smooth projective varieties are generally larger and harder to control than those arising in our toric examples.
In particular, the varieties $Y_\Qsf$ are typically singular.
We do not know when it is possible to describe the corresponding windows $W_S$ in terms of push-pull operations as in \cref{prop:toric_variety_geometric}.

\begin{ex}
    For $\Xfr = \PP^n$, if $\Qsf$ is the quiver corresponding to (the opposite of) the Beilinson collection $\{ \Osc_{\PP^n}(-n+1), \dots, \Osc_{\PP^n}$, then $Y_\Qsf$ is a closed subscheme of $\AA^{n(n+1)}$, and $G_\Qsf \cong \GG_m^n$.
    By contrast, the toric quotient presentation allows us to view $\Xfr$ as an open substack of $[\AA^n / \GG_m]$.
\end{ex}

\section{Quiver tensor products: the basic case} \label{sec:basic}

\begin{notn} \label{notn:homomorphism}
    Let $G$ be a commutative reductive group.
    Let $\phi: \Mfr \to BG$ be a homomorphism of $\EE_n$-monoid derived stacks which is (as a morphism of derived stacks) affine and almost of finite type, so $\Mfr$ is a perfect stack.
    In particular, taking $M = \Mfr \times_{BG} \Spec k$, the affine derived scheme $M$ is an $\EE_n$-monoid derived scheme, and the natural morphism $M \to \Mfr$ is a morphism of $\EE_n$-monoids.
    Write $\mu: \Mfr \times \Mfr \to \Mfr$ for the binary multiplication map.
\end{notn}

\begin{rmk}
    In general, we expect that such an $\EE_n$-homomorphism $\Mfr \to BG$ is obtained from a ``normal'' $\EE_n$-homomorphism $G \to M$ (defined analogously to the inclusion of a normal subgroup).
    However, we are not aware of a workable characterization of normal $\EE_n$-homomorphisms in the $\infty$-categorical context.
    In our classical applications, the $\EE_n$-structure on $\Mfr \to BG$ can be checked ``by hand,'' so we leave development of a theory of normal $\EE_n$-homomorphisms to future work.
\end{rmk}

By \cref{prop:basic_underlying_equivalence}, we know that $\QC(\Mfr) \simeq \Dsf(\Qsf(\phi)\op)$ and that the usual tensor product $\otimes_\Osc$ on $\Qsf(\phi)$ corresponds to the Day convolution product on $\Dsf(\Qsf(\phi)\op)$.
Under the above hypotheses, $\QC(\Mfr)$ also admits an $\EE_n$-monoidal convolution product $\star_\Mfr$ by the results of \cref{sub:3ff_convolution}.
Our goal in this section is to show (\cref{thm:basic_comparison}) that $\star_\Mfr$ corresponds to a ``quiver tensor product'' $\otimes_\Qsf$ on $\Dsf(\Qsf(\phi)\op)$.
Compared with \cref{prop:basic_underlying_equivalence}, establishing \cref{thm:basic_comparison} as a homotopy-coherent $\EE_n$-monoidal equivalence is nontrivial, and we are not aware of versions of this result in the literature (even at the level of homotopy categories). 

\subsection{The quiver tensor product} \label{sub:quiver_tensor_definition}

To rigorously construct the quiver tensor product on $\Dsf(\Qsf(\phi)\op)$, first recall that we may view $\QC(\Mfr)$ as an $(\EE_n, \EE_\infty)$-bialgebra in $\Pr_k^L$ by \cref{cor:QC_convolution_bialgebra}.
The multiplication on $\QC(\Mfr)$ is the usual tensor product of quasicoherent sheaves, while the (binary) comultiplication is the pullback functor
\[
    \mu^*: \QC(\Mfr) \to \QC(\Mfr \times \Mfr) = \QC(\Mfr) \otimes_k \QC(\Mfr),
\]
where $\mu: \Mfr \times \Mfr \to \Mfr$ is the binary multiplication associated with the monoid structure on $\Mfr$.
In particular, for $\chi \in \weight(G)$, we have
\begin{equation} \label{eq:pullback_computation}
    \mu^* \Osc_\Mfr(\chi) = \Osc_\Mfr(\chi) \boxtimes \Osc_\Mfr(\chi),
\end{equation}
so the $\EE_n$-coalgebra structure on $\QC(\Mfr) \simeq \Dsf(\Qsf(\phi)\op)$ restricts to an $\EE_n$-coalgebra structure on $\Qsf(\phi)$.
Because $\Qsf(\phi)$ is closed under the usual tensor product on $\QC(\Mfr)$, we actually get a stronger statement: the $(\EE_\infty, \EE_n)$-bialgebra structure on $\QC(\Mfr)$ (defined by pullbacks) restricts to an $(\EE_\infty, \EE_n)$-bialgebra structure on $\Qsf(\phi)$.

\begin{ex} \label{ex:toric_coalgebra}
    Consider the setting of \cref{dfn:prelinear_quiver}.
    Because the category $\Cat$ of small discrete categories is Cartesian monoidal, every small category $\Cc$ is (uniquely) an $\EE_\infty$-coalgebra in $\Cat$ by \HA{2.4.3.10}.\footnote{This does not require discreteness -- the only reason we mention ``discreteness'' here is that the categories $\Qsf^\pre(\phi)$ of \cref{dfn:prelinear_quiver} are discrete.}
    The comultiplication on $\Cc$ is given by the diagonal $\Delta_\Cc: \Cc \to \Cc \times \Cc$.
    In particular, the category $\Qsf^\pre(\phi)$ obtains an $\EE_\infty$-coalgebra structure in this way, and the above $\EE_\infty$-coalgebra structure on $\Qsf(\phi)$ is the $k$-linearization of that on $\Qsf^\pre(\phi)$.
\end{ex}

By \cref{cor:maps_form_algebra}, the $\EE_n$-coalgebra structure on $\Qsf(\phi)$ combines with the symmetric monoidal structure $\otimes_k$ on $\Dsf(k)$ to produce a \emph{quiver tensor product} $\otimes_\Qsf := \otimes_{\Qsf(\phi)}$ on $\Dsf(\Qsf(\phi)\op) = \Fun_k(\Qsf(\phi)\op, \Dsf(k))$.
More concretely, this is defined for $V_1, V_2: \Qsf(\phi)\op \to \Dsf(k)$ and $\chi \in \weight(G)$ by
\[
    (V_1 \otimes_{\Qsf} V_2)(\chi) = V_1(\chi) \otimes_k V_2(\chi).
\]
The $\EE_n$-coalgebra structure on $\Qsf(\phi)$ is used to define the behavior of $\otimes_\Qsf$ on morphisms.

\subsection{Moduli of representations and comparison functors}

To establish \cref{thm:basic_comparison}, we will use the moduli stack of representations of a small $k$-linear $\infty$-category, as developed in \cite{tv_moduli} and subsequent papers (see in particular \cite[\S 5]{ag_brauer} for a discussion in the language of $\infty$-categories).
We review the theory of these moduli stacks here.

\begin{dfn}[{\cite[Definition 3.2]{tv_moduli}}]
    Let $\Repfr: \Cat_k\op \to \dStk_k$ be the functor defined by
    \[
		\Repfr(\Cc)(R) = \Fun_k(\Cc, \Perf(R))^\simeq
	\]
    for $\Cc \in \Cat_k$ and $R \in \dCAlg_k$.
    (This defines a derived stack for the \'etale topology by \cite[Lemma 3.1]{tv_moduli} -- see also \cite[Lemma 5.4]{ag_brauer} for an $\infty$-categorical treatment.)
	We call $\Repfr(\Cc)$ the \emph{moduli stack of representations of $\Cc$}.
\end{dfn}

\begin{rmk}
    We use conventions opposite to those of \cite{tv_moduli}, which refers to what we call $\Repfr(\Cc)$ as the \emph{moduli stack of (pseudo-perfect) $\Cc\op$-modules}. 
\end{rmk}

\begin{rmk}
    The moduli stacks $\Repfr(\Cc)$ are typically very ``large'' and poorly behaved.
    For example, $\Repfr(k)$ is the moduli stack of all perfect complexes.
    This will not present difficulties for us: we are not interested in the stacks $\Repfr(\Cc)$ themselves.
    Instead, we want to study maps $\Xfr \to \Repfr(\Cc)$ where $\Xfr$ is ``small'' and well-behaved (e.g.\ perfect).
\end{rmk}

Recall that:
\begin{itemize}
    \item $\Indsf(-) = \Fun((-)\op, \Dsf(k))$ sends a small $k$-linear $\infty$-category $\Cc$ to the $\infty$-category of $k$-linear presheaves on $\Cc$.
    \item $(-)^\omega$ sends a presentable $k$-linear $\infty$-category to its (small) full subcategory of compact objects.
\end{itemize}
The left adjoint to the inclusion $\Cat_k^\perf \hookrightarrow \Cat_k$ is given by $\Cc \mapsto (\Indsf \Cc)^\omega$.
Thus $\Repfr(\Cc) \simeq \Repfr(\Indsf(\Cc)^\omega)$ for all $\Cc \in \Cat_k$.
We will implicitly identify these moduli stacks.
As a consequence, we lose no generality by restricting the domain of $\Repfr$ to $(\Cat_k^\perf)\op$.

Viewing $\Repfr$ as a functor $(\Cat_k^\perf)\op \to \dStk_k$ allows us to find an explicit left adjoint to $\Repfr$.
In fact, we can do even better: $\Repfr$ fits into a lax symmetric monoidal adjunction (cf.\ \cref{dfn:lax_monoidal_adjunction}).
Here subscripts $(-)_*$, resp.\ superscripts $(-)^*$, are used to indicate that a functor sends morphisms to the corresponding pushforward functors, resp.\ pullback functors.

\begin{prop} \label{prop:moduli_adjunction}
    There are lax symmetric monoidal adjunctions
    \[
        (\Perf^*)\op: (\dStk_k, \times) \rightleftarrows \big(\big(\Cat_k^\perf\big)\op, \otimes_k\big) :\Repfr
    \]
    and
    \[
        \Indsf\Perf_*: (\dStk_k, \times) \rightleftarrows \big(\Pr_{k,\omega}^R, \otimes\big) :\Repfr((-)^\omega).
    \]
\end{prop}

\begin{proof}
    By using the symmetric monoidal equivalence 
    \[
        \Indsf: \big(\big(\Cat_k^\perf\big)\op, \otimes_k\big) \xleftrightarrow{\sim} \big(\Pr_{k,\omega}^R, \otimes_k\big) : (-)^\omega
    \]
    we see that it suffices to construct the former adjunction.
    Existence of this adjunction is standard (cf.\ \cite[Proposition 3.4]{toen_homotopy}), though we recall the argument in $\infty$-categorical language for completeness.
    
    Let $\Cc \in \Cat_k^\perf$.
    Any derived stack $\Xfr: \dAff_k\op \to \Sc$ may be written as a colimit of representables, i.e.\ $\Xfr = \colim_i \Spec R_i$ for some $R_i$.
    We may then compute
    \begin{align*}
        \Map_{\dStk_k}\big(\Xfr, \Repfr(\Cc)\big) &= \Map_{\dStk_k}\big(\colim_i \Spec R_i, \Repfr(\Cc)\big) \\
        &= \lim_i \Map_{\dStk_k}\big(\Spec R_i, \Repfr(\Cc)\big) \\
        &= \lim_i \Fun_k\big(\Cc, \Perf(R_i)\big) \textrm{ by definition of $\Repfr(\Cc)$} \\
        &= \Fun_k\big(\Cc, \lim_i \Perf(R_i)\big) \\
        &= \Fun_k\big(\Cc, \Perf(\Xfr)\big) \textrm{ by definition of $\Perf(\Xfr)$.} 
    \end{align*}
    This is the needed adjunction.

    To get the lax symmetric monoidal structure on the adjunction, we just need to find a natural oplax symmetric monoidal structure on the left adjoint $(\Perf^*)\op$.
    The lax symmetric monoidal structure on $\QC^*: \dStk_k\op \to \Pr_k^L$ restricts to a lax symmetric monoidal structure on $\Perf^*: \dStk_k\op \to \Cat_k^\perf$.
    Reversing the direction of arrows, we obtain an oplax symmetric monoidal structure on $(\Perf^*)\op: \dStk_k \to \big(\Cat_k^\perf\big)\op$ as needed.
\end{proof}

In particular, for a derived stack $\Xfr \in \dStk_k$ and a small $k$-linear $\infty$-category $\Cc \in \Cat_k$, the following data are equivalent:
\begin{itemize}
    \item A morphism of derived stacks $f_{\Esc_\bullet}: \Xfr \to \Repfr(\Cc) \simeq \Repfr(\Indsf(\Cc)^\omega)$.
    \item A $k$-linear functor $\Esc_\bullet: \Cc \to \Perf(\Xfr)$.
    \item A $k$-linear functor of the form $\Hom(\Esc_\bullet, -): \Indsf\Perf(\Xfr) \to \Indsf(\Cc)$.
\end{itemize}
Here $\Hom(\Esc_\bullet, -)$ refers to the functor satisfying
\[
    \Hom(\Esc_\bullet, \Fsc)(c) = \Hom_{\Indsf\Perf(\Xfr)}(\Esc_c, \Fsc)
\]
for $\Fsc \in \Indsf\Perf(\Xfr)$ and $c \in \Cc$.
The functors of this form are precisely the functors from $\Indsf\Perf(\Xfr)$ to $\Indsf(\Cc)$ in $\Pr_{k,\omega}^R$ by a version of Yoneda's lemma.

Suppose that $\Cc \in \Alg_{\EE_n}(\Cat_k\op)$, i.e.\ $\Cc$ is an $\EE_n$-coalgebra in small $k$-linear $\infty$-categories.
Then $\Indsf(\Cc)$ is an $\EE_n$-coalgebra in $\Pr_{k,\omega}^L$.
Applying the equivalence $(\Pr_{k,\omega}^L)\op \simeq \Pr_{k,\omega}^R$, we see that $\Indsf(\Cc)$ is an $\EE_n$-algebra in $\Pr_{k,\omega}^R$, i.e.\ an $\EE_n$-monoidal category such that the monoidal structure preserves colimits and limits (formed in the powers $\Indsf(\Cc)^{\boxtimes k}$).

Furthermore, for $R \in \dCAlg_k$, \cref{cor:maps_form_algebra} lets us combine the symmetric monoidal structure $\otimes_\Osc$ on $\Perf(R)$ with the $\EE_n$-coalgebra structure on $\Cc$ to produce an $\EE_n$-monoid structure on the space $\Repfr(\Cc)(R) = \Fun_k(\Cc, \Perf(R))^\simeq$.
Thus $\Repfr(\Cc)$ is an $\EE_n$-monoid object of $\dStk_k$.
This $\EE_n$-monoid structure gives a useful source of symmetric monoidal functors into $(\Indsf(\Cc), \otimes_\Cc)$:

\begin{prop} \label{prop:general_comparison}
    Let $\Xfr \in \Alg_{\EE_n}(\dStk_k^\perf)$ be a perfect $\EE_n$-monoid derived stack, and let $\Cc \in \Alg_{\EE_n}(\Cat_k\op)$.
    There is a natural equivalence
    \[
        \Map_{\Alg_{\EE_n}(\dStk_k)}(\Xfr, \Repfr(\Cc)) \simeq \Fun_k^{\EE_n,R}\big((\QC(\Xfr), \star_\Xfr), (\Indsf(\Cc), \otimes_\Cc)\big)^\simeq
    \]
    enhancing the adjunction equivalence of \cref{prop:moduli_adjunction}.
    In particular, if $f_{\Esc_\bullet}: \Xfr \to \Repfr(\Cc)$ is a homomorphism of commutative monoid derived stacks, then the functor
    \[
        \Hom(\Esc_\bullet, -): (\QC(\Xfr), \star_\Xfr) \to (\Indsf(\Cc), \otimes_\Cc)
    \]
    is $\EE_n$-monoidal.
\end{prop}

\begin{proof}
    The adjunction of \cref{prop:moduli_adjunction} is only lax symmetric monoidal in general.
    However, on the full subcategory $\dStk_k^\perf \subset \dStk_k$, the functor $\Indsf\Perf_*$ is symmetric monoidal and agrees with $\QC_*$.
    The result follows from \cref{prop:partial_monoidal_adjunction_algebra} (keeping in mind \cref{rmk:O_algebra_in_SMC}), using the fact that $\Fun_k^{\EE_n,R}(-, -)^\simeq$ is the space of morphisms in $\Alg_{\EE_n}\big(\Pr_k^R\big)$ by definition.
\end{proof}

\subsection{Proof of the $\EE_n$-monoidal equivalence}

We have now assembled all of the ingredients necessary to upgrade the underlying equivalence \cref{prop:basic_underlying_equivalence} to a symmetric monoidal equivalence

\begin{thm} \label{thm:basic_comparison}
    The underlying equivalence of categories of \cref{prop:basic_underlying_equivalence} upgrades to an $\EE_n$-monoidal equivalence
    \[
        \big(\QC(\Mfr), \star_\Mfr\big) \simeq \big(\Dsf(\Qsf(\phi)\op), \otimes_{\Qsf(\phi)\op}\big).
    \]
\end{thm}

\begin{proof}
    Let $R \in \dCAlg_k$.
    There is a natural $\EE_n$-monoid structure on $\Fun^\otimes(\Qsf(\phi), \Perf^{\leq 0}(R))^\simeq$ constructed as follows.
    Note that $\Qsf(\phi)$ is an $(\EE_\infty, \EE_n)$-bialgebra in $(\Cat_k^\perf, \otimes_k)$, or equivalently an $\EE_n$-coalgebra in $\big(\CAlg(\Cat_k^\perf), \otimes_k\big)$.
    Furthermore, $\Perf^{\leq 0}(R)$ is a symmetric monoidal $\infty$-category, i.e.\ an object of $\CAlg(\Cat_k^\perf) = \CAlg(\CAlg(\Cat_k^\perf))$ (by \HA{3.2.4.5}).
    Thus the mapping space 
    \[
        \Fun^\otimes(\Qsf(\phi), \Perf^{\leq 0}(R))^\simeq = \Map_{\CAlg(\Cat_k^\perf)}(\Qsf(\phi), \Perf^{\leq 0}(R))
    \]
    obtains an $\EE_n$-algebra structure from \cref{cor:maps_form_algebra}. 
    
    By \cref{cor:tannakian_description}, we may identify $\Mfr(R) \simeq \Fun^\otimes(\Qsf(\phi), \Perf^{\leq 0}(R))^\simeq$ as $\EE_n$-monoids.
    Forgetting the symmetric monoidal structure on functors and applying the inclusion $\Perf^{\leq 0}(R) \hookrightarrow \Perf(R)$, we obtain a map
    \[
        f_{\Osc(-),R}: \Mfr(R) \simeq \Fun^\otimes(\Qsf(\phi), \Perf^{\leq 0}(R))^\simeq \to \Fun(\Qsf(\phi), \Perf(R))^\simeq =: \Repfr(\Qsf(\phi))(R).
    \]
    Each map $f_{\Osc(-),R}$ is a homomorphism of $\EE_n$-monoid spaces by \cref{cor:maps_form_algebra} and \cref{cor:maps_form_algebra_functorial}.
    The maps $f_{\Osc(-),R}$ are natural in $R$ and thus assemble into a map of $\EE_n$-monoid derived stacks $f_{\Osc(-)}: \Mfr \to \Repfr(\Qsf(\phi))$.
    
    By \cref{prop:general_comparison}, $f_{\Osc(-)}$ gives an $\EE_n$-monoidal functor
    \begin{align*}
        \Hom(\Osc(-), -): \big(\QC(\Mfr), \star_\Mfr\big) &\to \big(\Dsf(\Qsf(\phi)\op), \otimes_\Qsf\big) \\
        \Fsc &\mapsto \big(\chi \mapsto \Hom(\Osc_\Mfr(\chi), \Fsc)\big)
    \end{align*}
    The underlying functor of this $\EE_n$-monoidal functor is the inverse of the underlying equivalence of \cref{prop:basic_underlying_equivalence}.
    Thus $\Hom(\Osc(-), -)$ gives the desired $\EE_n$-monoidal equivalence.
\end{proof}

\begin{rmk} \label{rmk:moulinos}
    The proof of \cref{thm:basic_comparison} works just as well over the sphere spectrum $\SS$ provided that we know $\QC(BG)$ is compactly generated by invertible objects.
    In particular, taking $G = \GG_m$ (the flat multiplicative group over $\SS$) and $M = \AA^1$ (the flat affine line over $\SS$), we extend \cite[Theorem 1.1]{moulinos_geometry} to an equivalence
    \[
        \big(\QC([\AA^1 / \GG_m]), \star\big) \simeq \big(\Fun((\ZZ, \leq)\op, \Sp), \otimes_\SS\big).
    \]
    It would be interesting to understand the topological implications of this statement.
\end{rmk}

\section{Extended convolution products} \label{sec:ec}

Suppose $\Mfr$ is a perfect $\EE_n$-monoid stack with multiplication $\mu: \Mfr \times \Mfr \to \Mfr$.
If $\Xfr$ is an open substack of $\Mfr$, then the multiplication on $\Mfr$ typically does not induce a multiplication on $\Xfr$: there is no reason for $\Xfr$ to be closed under $\mu$.

However, the situation is somewhat different if we categorify.
Sometimes the convolution product $\star_\Mfr$ on $\QC(\Mfr)$ (together with a choice of window $W_S$) induces an ``extended convolution'' product $\star'_{\Mfr,S}$ on $\QC(\Xfr)$ \emph{even when $\Xfr$ is not closed under $\mu$}.
Our goal in this section is to demonstrate this claim and understand the behavior of $\star'_{\Mfr,S}$ in the case where $\Mfr$ is a (nice) global quotient $[M / G]$.

\subsection{Definitions and examples}

We begin by formalizing the data we will use to construct extended convolution products.

\begin{dfn} \label{dfn:ec}
    For $1 \leq n \leq \infty$, an $\EE_n$-\emph{extended convolution setup} (or $\EE_n$-\emph{EC setup}) is a triple 
    \[
        (\phi: \Mfr \to BG, j: \Xfr \hookrightarrow \Mfr, S)
    \]
    where:
    \begin{itemize}
        \item $G$ is a commutative reductive group over $k$.
        \item $\phi: \Mfr \to BG$ is a homomorphism of $\EE_n$-monoid derived stacks which is (as a morphism of derived stacks) affine and almost of finite type.
        \item $j: \Xfr \hookrightarrow \Mfr$ is an open immersion.
        \item $S \subset \weight(G)$ is a collection of weights which is transparent for $\Xfr \subset \Mfr$.
    \end{itemize}
    A \emph{geometric} $\EE_n$-EC setup is an $\EE_n$-EC setup as above together with:
    \begin{itemize}
        \item a morphism of perfect stacks $q = (q_1, q_2): \Wc_{-S} \to \Xfr \times \Mfr$ such that $W_{-S} = q_{2*} q_1^*$.
    \end{itemize}
    Given an $\EE_n$-EC setup as above, we write $\mu: \Mfr \times \Mfr \to \Mfr$ for the binary multiplication.
\end{dfn}

Before developing the theory of extended convolution, we introduce a few key examples of $\EE_n$-EC setups.

\begin{ex} \label{ex:algebra_ec_setup}
    Let $A$ be a finite dimensional $k$-algebra.
    The natural map $k^\times \subset k \to A$ corresponds to a multiplicative homomorphism $\GG_m \to \AA(A)$, and the quotient stack $[\AA(A) / \GG_m]$ inherits a monoid structure.
    Let $\Xfr = \PP(A)$, and identify $\weight(\GG_m) \simeq \ZZ$ in such a way that $\Oc_{\PP(A)}(1) \in \Pic \PP(A) \simeq \weight(\GG_m)$ corresponds to $1 \in \ZZ$.
    Write $S = \{ 0, \dots, \dim_k A \} \subset \ZZ$.
    Then 
    \[
        (\phi: [\AA(A) / \GG_m] \to B\GG_m, \PP(A) \hookrightarrow [\AA(A) / \GG_m], S)
    \]
    is a geometric $\EE_1$-EC setup, where $q: \Wc_{-S} \to \PP(A) \times [\AA(A) / \GG_m]$ is the morphism of \cref{prop:toric_variety_geometric}.
    If $A$ is commutative, the above $\EE_1$-EC setup upgrades to an $\EE_\infty$-EC setup.
\end{ex}

\begin{ex} \label{ex:toric_ec_setup}
    Let $\Xfr$ be a smooth toric stack with a decent open immersion $j: \Xfr \hookrightarrow [\AA^n / G]$.
    We may view $[\AA^n / G]$ as an $\EE_\infty$-monoid stack, where $\AA^n$ is equipped with its coordinatewise multiplication.
    Let $S$ be any collection of weights corresponding to a full strong exceptional collection of line bundles in $\Perf(\Xfr)$.
    Then 
    \[
        (\phi: [\AA^n / G] \to BG, j: \Xfr \to [\AA^n / G], S)
    \]
    is an $\EE_\infty$-EC setup (though there is no reason \emph{a priori} for it to be geometric).
\end{ex}

We may do even better in the Bondal-Ruan case:

\begin{ex} \label{ex:toric_geometric_ec_setup}
    Suppose that $\Xfr$ is a smooth complete toric variety of Bondal-Ruan type, and let $j: \Xfr \hookrightarrow [\AA^n / G]$ be the open immersion arising from the Cox construction.
    Then
    \[
        (\phi: [\AA^n / G] \to BG, j: \Xfr \to [\AA^n / G], -\Theta_\Xfr)
    \]
    is a geometric $\EE_\infty$-EC setup by \cref{prop:toric_variety_geometric}.
\end{ex}

Given an $\EE_n$-EC setup $(\phi, j, S)$ with notation as in \cref{dfn:ec}, the category $\QC(\Xfr)$ inherits a convolution product from that of $\QC(\Mfr)$. 

\begin{prop} \label{prop:ec_exists}
    Let $(\phi, j, S)$ be an $\EE_n$-EC setup with notation as in \cref{dfn:ec}.
    There exists a unique $\EE_n$-monoidal structure $\star'_{\Mfr,S}$ on $\QC(\Xfr)$ such that the Hitchcock functor 
    \[
        H_S: (\QC(\Mfr), \star_\Mfr) \to (\QC(\Xfr), \star_{\mu,S}')
    \]
    is $\EE_n$-monoidal.
\end{prop}

\begin{proof}
    We check that the adjunction 
    \[
        W_S: \QC(\Xfr) \rightleftarrows \QC(\Mfr) :H_S
    \]
    satisfies the conditions of \cref{prop:transport_convolution}.
    Conditions (1) and (2) are satisfied by the definition of $W_S$.
    Condition (3) is satisfied by the definition of the convolution product: $- \star_\Mfr - = \mu_*(- \boxtimes -)$ and $\mu_*$ has left adjoint $\mu^*$.
    To check condition (4), note that $\QC(\Xfr)$ is compactly generated by $\{ \Osc_\Xfr(\chi) \}_{\chi \in S}$ and $W_S(\Osc_\Xfr(\chi)) = \Osc_{\Mfr}(\chi)$ for $\chi \in S$.
    Since
    \[
        \mu^* \Osc_\Mfr(\chi) \simeq \Osc_\Mfr(\chi) \boxtimes \Osc_\Mfr(\chi)
    \]
    by $G$-equivariance of the multiplication on $\Mfr$, condition (4) holds.
\end{proof}

We call $\star'_{\Mfr,S}$ the \emph{extended convolution product} (or \emph{EC product}) associated with the $\EE_n$-EC setup $(\phi, j, S)$.
The name is justified by \cref{cor:ec_extends_convolution}, which implies in particular that $\star'_{\Mfr,S}$ extends the convolution product on any open submonoid stack of $\Mfr$ contained in $\Xfr$.

By \cref{prop:right_localization_properties}, the $\EE_n$-monoidal functor $H_S$ satisfies the following universal property: if $F: (\QC(\Mfr), \star_\Mfr) \to (\Cc, \otimes_\Cc)$ is an $\EE_n$-monoidal functor and the underlying functor of $F$ satisfies $F = F' \circ H_S$, then there is a unique $\EE_n$-monoidal structure on $F'$ such that $F = F' \circ H_S$ as $\EE_n$-monoidal functors.

\subsection{Functoriality of extended convolution}

To understand the behavior of EC products geometrically, it is useful (though perhaps not strictly necessary for later developments) to define morphisms of $\EE_n$-EC setups. 
We do so here and show (\cref{prop:ec_functoriality}) that morphisms of $\EE_n$-EC setups induce $\EE_n$-monoidal functors.

\begin{dfn} \label{dfn:ec_morphism}
    Let $(\phi_1, j_1, S_1)$ and $(\phi_2, j_2, S_2)$ be $\EE_n$-EC setups with notation as in \cref{dfn:ec}.
    For notational simplicity, write $j_i = j_{\Xfr_i}$, $W_i = W_{S_i}$, $H_i = H_{S_i}$, and $\star'_i = \star'_{\Mfr_i,S_i}$ for $i = 1, 2$.
    A \emph{morphism of $\EE_n$-EC setups} $(\alpha, \beta): (\phi_1, j_1, S_1) \to (\phi_2, j_2, S_2)$ consists of:
    \begin{itemize}
        \item an $\EE_n$-homomorphism $\alpha: \Mfr_1 \to \Mfr_2$ and
        \item a group homomorphism $\beta: G_1 \to G_2$,
    \end{itemize}
    such that the diagram
    \[
        \begin{tikzcd}
            \Mfr_1 \rar["\alpha"] \dar["\phi_1"] & \Mfr_2 \dar["\phi_2"] \\
            BG_1 \rar["B\beta"] & BG_2
        \end{tikzcd}
    \]
    commutes and $\im \alpha^* W_2 \subset \im W_1$. 
\end{dfn}

\begin{ex}
    If $G_1 = G_2 = G$ and $S_1 = \weight(G)$ (so $\Xfr_1 = \Mfr_1$), then any morphism $\alpha: \Mfr_1 \to \Mfr_2$ over $BG$ gives a morphism of $\EE_n$-EC setups $(\alpha, \id_G): (\phi_1, \id_{\Xfr_1}, \weight(G)) \to (\phi_2, j_2, S_2)$.
\end{ex}

\begin{ex}\label{ex:ec_morph_fdalg}
    Let $f: A \to B$ be a homomorphism of finite-dimensional $k$-algebras such that $\dim A \geq \dim B$.
    Then the natural map $\alpha_f: [\AA(A) / \GG_m] \to [\AA(B) / \GG_m]$ gives a morphism of $\EE_1$-EC setups
    \begin{align*}
        (\alpha_f, \id_{\GG_m}): &\big(\phi_A: [\AA(A) / \GG_m] \to B\GG_m, \PP(A) \hookrightarrow [\AA(A) / \GG_m], \{0, \dots, \dim_k A\}\big) \\
        &\to \big(\phi_B: [\AA(B) / \GG_m] \to B\GG_m, \PP(B) \hookrightarrow [\AA(B) / \GG_m], \{0, \dots, \dim_k B\}\big).
    \end{align*}
    If $A$ and $B$ are both commutative, then $\alpha_f$ is in fact a morphism of $\EE_\infty$-EC setups.
\end{ex}

Morphisms of $\EE_n$-EC setups induce functors that preserve the corresponding EC products.
To prove this, we first need a lemma allowing us to simplify certain composite functors involving windows and Hitchcock functors.

\begin{lem} \label{lem:window_functoriality}
    For $i \in \{1, 2\}$, fix the following data:
    \begin{itemize}
        \item A commutative reductive group $G_i$,
        \item A morphism of derived stacks $\phi_i: \Yfr_i \to BG_i$ which is affine and almost of finite type,
        \item An open immersion $j_i: \Xfr_i \hookrightarrow \Yfr_i$, and
        \item A transparent collection of weights $S_i \subset \weight(G_i)$ for $\Xfr_i \subset \Yfr_i$.
    \end{itemize}
    Write $W_i: \QC(\Xfr_i) \to \QC(\Yfr_i)$ and $H_i: \QC(\Yfr_i) \to \QC(\Xfr_i)$ for the corresponding window and Hitchcock functor (respectively).
    Let $f: \Yfr_1 \to \Yfr_2$ be a morphism such that $f^*(\im W_2) \subset \im W_1$.\footnote{The functor $W_1$ is fully faithful, so it suffices to check this inclusion on objects.}
    Then:
    \begin{enumerate}
        \item $f^* W_2 = W_1 j_1^* f^* W_2$.
        \item $H_2 f_* = H_2 f_* j_{1*} H_1$.
        \item $H_2 f_* W_1$ = $H_2 f_*j_{1*}$.
    \end{enumerate}
\end{lem}

\begin{proof}
    (1). Suppose $\Fsc \in \QC(\Xfr_2)$.
    Then $f^* W_2 \Fsc = W_1 \Gsc$ for some $\Gsc \in \QC(\Xfr_1)$.
    Applying $j_1^*$ gives 
    \[
        j_1^* f^* W_2 \Fsc = j_1^* W_1 \Gsc = \Gsc,
    \]
    so 
    \[
        f^* W_2 \Fsc = W_1 \Gsc = W_1 j_1^* f^* W_2 \Fsc.
    \]
    Naturality of this isomorphism is clear because $W_1$ is fully faithful.

    (2). Take right adjoints of all functors involved in (1).

    (3). We compute
    \begin{align*}
        H_2 f_* W_1 &= H_2 f_* j_{1*} H_1 W_1 \textrm{ by (2)} \\
        &= H_2 f_* j_{1*} \id_{\Xfr_1} \textrm{ because $W_1$ is fully faithful} \\
        &= H_2 f_* j_{1*}. \qedhere
    \end{align*}
\end{proof}

\begin{prop} \label{prop:ec_functoriality}
    Let $(\alpha, \beta): (\phi_1, j_1, S_1) \to (\phi_2, j_2, S_2)$ be a morphism of $\EE_n$-EC setups with notation as above.
    Then the functor $H_2 \alpha_* j_{1*}: \big(\QC(\Xfr_1), \star'_1\big) \to \big(\QC(\Xfr_2), \star'_2\big)$ is symmetric monoidal. 
\end{prop}

\begin{proof}
    By \cref{lem:window_functoriality}(2), we have $H_2 \alpha_* = H_2 \alpha_* j_{1*} H_1$.
    Using the universal property of $H_1$ mentioned in \cref{prop:ec_exists}, it suffices to show that $H_2 \alpha_*$ is $\EE_n$-monoidal.
    But this is clear as both $H_2$ and $\alpha_*$ are $\EE_n$-monoidal.
\end{proof}

\begin{cor} \label{cor:ec_extends_convolution}
    Let $(\phi, j, S)$ be an $\EE_n$-EC setup with notation as in \cref{dfn:ec}.
    Let $\alpha: \Nfr \to \Mfr$ be a homomorphism of $\EE_n$-monoid derived stacks over $BG$, and assume $\alpha$ factors through the inclusion $j: \Xfr \hookrightarrow \Yfr$, say $\alpha = j \circ a$.
    Then the pushforward functor $a_*: \big(\QC(\Nfr), \star_{\Nfr}\big) \to \big(\QC(\Xfr), \star'_{\Mfr,S}\big)$ is $\EE_n$-monoidal.
\end{cor}

\begin{proof}
    In \cref{prop:ec_functoriality}, take $(\phi_1, j_1, S_1) = \big(\phi, \id_{\Xfr}, \weight(G)\big)$, $(\phi_2, j_2, S_2) = (\phi, j, S)$, $\alpha = j \circ a$, and $\beta = \id_G$.
    Then $H_2 \alpha_* j_{1*} H_1 = H_S j_* a_* = a_*$.
\end{proof}

\subsection{Fourier-Mukai kernels and geometric descriptions of extended convolution}

We would like to understand the operations $\star'_{\Mfr,S}$ using the geometry of $\Xfr$.
This is difficult in general for the simple reason that the window associated with a transparent collection is hard to understand geometrically.
However, we shall show that giving a geometric description of $\star'_{\Mfr,S}$ is \emph{no more difficult} than giving a geometric description of the Fourier-Mukai kernel $\Ksc_{-S}$ of $W_{-S}$.
In particular, for geometric $\EE_n$-EC setups, we obtain a simple geometric description of $\star'_{\Mfr,S}$.

For simplicity, our claims here will be made for the binary product $\star'_{\Mfr,S}: \QC(\Xfr) \boxtimes \QC(\Xfr) \to \QC(\Xfr)$.
The analogous claims for $n$-ary products can be established using the same arguments.

\begin{prop} \label{prop:ec_fm_computation}
    Let $(\phi, j, S)$ be an $\EE_n$-EC setup with notation as in \cref{dfn:ec}.
    View $\Ksc_{-S}$ as an object of $\QC(\Mfr \times \Xfr)$.
    Then $\star'_{\Mfr,S}$ is given by the Fourier-Mukai transform with kernel $(j \times j \times \id_\Xfr)^* (\mu \times \id_\Xfr)^* \Ksc_{-S} \in \QC(\Xfr^3)$.
\end{prop}

\begin{proof}
    Let $\Fsc, \Gsc \in \QC(\Xfr)$.
    By \cref{lem:window_functoriality}(3), we have
    \[
        \Fsc \star'_{\Mfr,S} \Gsc = H_S(j_* \Fsc \star_\Mfr j_* \Gsc) = H_S \mu_* (j \times j)_* (\Fsc \boxtimes \Gsc).
    \]
    Writing $H_S = \Phi_{\Ksc_{-S}}$ and applying the general formula $\Phi_\Ksc \circ f_* = \Phi_{(f \times \id)^* \Ksc}$, we see the claim.
\end{proof}

If we have an explicit resolution of the diagonal, we may obtain an explicit algebraic description of the corresponding EC product.

\begin{cor}
    Suppose $(\phi, j, S)$ is an $\EE_n$-EC setup with notation as in \cref{dfn:ec}.
    Write $\Delta_{\Xfr*} \Osc_\Xfr$ as a complex
    \[
        \begin{tikzcd}
            \dots \rar & \displaystyle\bigoplus_{\chi \in S} \Osc_\Xfr(-\chi) \boxtimes \Asc_{i,\chi} \rar["d_i"] & \displaystyle\bigoplus_{\chi \in S} \Osc_\Xfr(-\chi) \boxtimes \Asc_{i+1,\chi} \rar & \dots
        \end{tikzcd}
    \]
    Then the Fourier-Mukai kernel of $\star'_{\Mfr,S}$ is
    \[
        \begin{tikzcd}
            \dots \rar & \displaystyle\bigoplus_{\chi \in S} \Osc_\Xfr(-\chi) \boxtimes \Osc_\Xfr(-\chi) \boxtimes \Asc_{i,\chi} \ar[rrr, "(j \times j \times \id_\Xfr)^* (\mu \times \id_\Xfr)^* d_i"] & & & \displaystyle\bigoplus_{\chi \in S} \Osc_\Xfr(-\chi) \boxtimes \Osc_\Xfr(-\chi) \boxtimes \Asc_{i+1,\chi} \rar & \dots
        \end{tikzcd}
    \]
\end{cor}

\begin{proof}
    By \cref{prop:diagonal_window}, the kernel $\Ksc_{-S} \in \QC(\Mfr \times \Xfr)$ is given by
    \[
        \begin{tikzcd}
            \dots \rar & \displaystyle\bigoplus_{\chi \in S} \Osc_\Mfr(-\chi) \boxtimes \Asc_{i,\chi} \rar["d_i"] & \displaystyle\bigoplus_{\chi \in S} \Osc_\Mfr(-\chi) \boxtimes \Asc_{i+1,\chi} \rar & \dots
        \end{tikzcd}
    \]
    Applying \cref{prop:ec_fm_computation} gives the result.
\end{proof}

\cref{prop:ec_fm_computation} also yields a geometric description of the EC products associated with geometric $\EE_n$-EC setups.

\begin{prop} \label{prop:ec_geometric_computation}
    Suppose $(\phi, j, S, q)$ is a geometric $\EE_n$-EC setup with notation as in \cref{dfn:ec}.
    Let $\Zfr_{\mu,S}$ be defined by the Cartesian square (in $\dStk_k$)
    \[
        \begin{tikzcd}
            \Zfr_{\mu,S} \ar[rr, "p_3"] \dar["p_1 \times p_2"] & & \Wc_{-S} \dar["q_2"] \\
            \Xfr \times \Xfr \rar["j \times j"] & \Mfr \times \Mfr \rar["\mu"] & \Mfr
        \end{tikzcd}
    \]
    Then, for $\Fsc, \Gsc \in \QC(\Xfr)$, there is a natural isomorphism $\Fsc \star'_{\Mfr,S} \Gsc = q_{1*} p_{3*} (p_1^*\Fsc \otimes p_2^* \Gsc)$.
\end{prop}

\begin{proof}
    By hypothesis, $\Ksc_{-S} = q_* \Osc_{\Wc_{-S}}$, so the Fourier-Mukai kernel of $\star'_{\Mfr,S}$ is
    \[
        (j \times j \times \id_\Xfr)^* (\mu \times \id_\Xfr)^* q_* \Osc_{\Wc_{-S}}
    \]
    by \cref{prop:ec_fm_computation}.
    Applying base change for the commutative square
    \[
        \begin{tikzcd}
            \Zfr_{\mu,S} \ar[rr, "p_3"] \dar["p_1  \times p_2 \times (q_1 \circ p_3)", swap] & & \Wc_{-S} \dar["q"] \\
            \Xfr^3 \rar["j \times j \times \id_\Xfr"] & \Mfr \times \Mfr \times \Xfr \rar["\mu \times \id_\Xfr"] & \Mfr \times \Xfr
        \end{tikzcd}
    \]
    shows that
    \[
        (j \times j \times \id_\Xfr)^* (\mu \times \id_\Xfr)^* q_* \Osc_{\Wc_{-S}} = (p_{1,2} \times (q_1 \circ p_3))_* \Osc_{\Zfr_{\mu,S}}.
    \]
    For $\Fsc, \Gsc \in \QC(\Xfr)$, we may now compute (letting $\pi_i: \Xfr^3 \to \Xfr$ be projection onto the $i$th coordinate):
    \begin{align*}
        \Fsc \star'_{\Mfr,S} \Gsc &= \pi_{3*}\big(\pi_1^* \Fsc \otimes \pi_2^* \Gsc \otimes (p_1 \times p_2 \times (q_1 \circ p_3))_* \Osc_{\Zfr_{\mu,S}}\big) \\
        &= \pi_{3*} \big(p_1 \times p_2 \times (q_1 \circ p_3)\big)_* \big(p_1^* \Fsc \otimes p_2^* \Gsc\big) \textrm{ by the projection formula} \\
        &= q_{1*} p_{3*}(p_1^* \Fsc \otimes p_2^* \Gsc). \qedhere
    \end{align*}
\end{proof}

\begin{ex} \label{ex:ec_geometric_diagonal}
    Suppose $\Xfr$ is a separated scheme and $\Wc_{-S} \subset \Xfr \times \Mfr$ is the restriction of the diagonal closed substack of $\Mfr$ to $\Xfr \times \Mfr$ (this is true for \cref{ex:algebra_ec_setup} and \cref{ex:toric_geometric_ec_setup}).
    Then $\Zfr_{\mu,S} \subset \Xfr^3$ is (the restriction to $\Xfr^3$ of) the closure of the graph of $\mu \circ (j \times j): \Xfr \times \Xfr \to \Mfr$.
    The legs of the correspondence
    \[
    \begin{tikzcd}
        {} & \Zfr_{\mu,S} \ar[dl, "p_1 \times p_2", swap] \ar[dr, "q_1 \circ p_3"] & \\
        \Xfr \times \Xfr & & \Xfr
    \end{tikzcd}
    \]
    are just the projections to the factors when $\Zfr_{\mu,S}$ is viewed as a closed substack of $\Xfr^3$.
\end{ex}

\begin{rmk}
    If $\alpha: (\phi_1, j_1, S_1) \to (\phi_2, j_2, S_2)$ is a morphism of $\EE_n$-EC setups, we may use similar methods to the above to describe the functor $H_2 \alpha_* j_{1*}$.
\end{rmk}

\subsection{Quiver tensor products via extended convolution}

Let $(\phi, j, S)$ be an $\EE_n$-EC setup with notation as in \cref{dfn:ec}.
Because $\mu^* \Osc_\Mfr(\chi) = \Osc(\chi) \boxtimes \Osc(\chi)$ for all $\chi \in \weight(G)$, the $\EE_n$-coalgebra structure on $\Qsf(\phi)$ restricts to an $\EE_n$-coalgebra structure on $\Qsf_S(\phi)$.
As in \cref{sub:quiver_tensor_definition}, this produces an $\EE_n$-monoidal \emph{quiver tensor product} $\otimes_{\Qsf_S(\phi)}$ on $\Dsf(\Qsf_S(\phi)\op)$.

\begin{thm} \label{thm:ec_comparison}
    Let $(\phi, j, S)$ be an $\EE_n$-EC setup with notation as in \cref{dfn:ec}.
    Then there is an $\EE_n$-monoidal equivalence
    \begin{align*}
        \big(\QC(\Xfr), \star'_{\Mfr,S}\big) &\xrightarrow{\sim} \big(\Dsf(\Qsf_S(\phi)\op), \otimes_{\Qsf_S(\phi)}\big) \\
        \Fsc &\mapsto \big(\chi \mapsto \Hom_\Xfr(\Osc_\Xfr(\chi), \Fsc)\big).
    \end{align*}
\end{thm}

\begin{proof}
    The functor in question is an equivalence by \cref{prop:transparent_induces_window}, so it suffices to show that said functor is $\EE_n$-monoidal.
    For this, let $i_S: \Qsf_S(\phi) \to \Qsf(\phi)$ be the inclusion, and recall that the diagram
    \[
        \begin{tikzcd}
            \QC(\Mfr) \rar["H_S"] \dar["\sim"] & \QC(\Xfr) \dar["\sim"] \\
            \Dsf(\Qsf(\phi)\op) \rar["i_S^*"] & \Dsf(\Qsf_S(\phi)\op)
        \end{tikzcd}
    \]
    commutes by the definition of $H_S$.
    Thus, by the universal property of \cref{prop:ec_exists}, it suffices to show that the composite $\QC(\Mfr) \to \Dsf(\Qsf_S(\phi)\op)$ is $\EE_n$-monoidal.
    But the equivalence $\big(\QC(\Mfr), \star_\Mfr\big) \xrightarrow{\sim} \big(\Dsf(\Qsf(\phi)\op), \otimes_\Qsf\big)$ is $\EE_n$-monoidal by \cref{thm:basic_comparison}, and functoriality of \cref{cor:maps_form_algebra} implies $i_S^*$ is $\EE_n$-monoidal, so the same must be true of the composite.
\end{proof}

\begin{rmk} \label{rmk:ec_perfect}
    As a consequence of \cref{thm:ec_comparison} and \cref{rmk:compact_representations}, we see that, if $S$ arises from a finite full strong exceptional collection of line bundles on $\Xfr$, then the EC product $\star'_{\Mfr,S}$ preserves perfect complexes.
    (One may also deduce this from the definition of $\star'_{\Mfr,S}$.)
    However, even in this case, we are not aware of a way to define $\star'_{\Mfr,S}$ geometrically without using the full categories $\QC$.
    Indeed, convolution on $\Mfr$ uses the functor $\mu_*$, which typically does not preserve perfect complexes.
\end{rmk}

\section{Application: (symmetric) monoidal structures on \texorpdfstring{$\Perf \PP^d$}{Perf P\^n}} \label{sec:Pn}

Using the methods of \cref{sec:ec}, we can construct new (symmetric) monoidal structures on $\Perf \PP^d$ from finite-dimensional algebras.
Recall that we use the term ``$\EE_n$-monoidal structures'' for the sole purpose of stating results about usual monoidal structures ($n =1$) and symmetric monoidal structures ($n = \infty$) concisely, and no other value of $n$ will be considered in this section (cf. \cref{rmk:en-algebra-in-1-cat} and \cref{rmk:en-algebra-in-2-cat}). 

Throughout the section, let $A$ be a finite-dimensional $k$-algebra with $\dim A = d + 1$, and let $\AA(A)$ be the corresponding affine monoid scheme.
By standard results on algebraic monoids (see e.g.\ \cite[page 1]{vinberg1995reductive}), the group of units $\AA(A)^\times$ is open and Zariski dense in $\AA(A)$.
We may use this to obtain the following:

\begin{prop} \label{mainthm_algebra}
    Let $A$ be a nonzero\footnote{When $A = 0$, we have $\PP(A) = \varnothing$, so the corresponding EC-setup exists but the description of the geometric structure fails.} finite-dimensional $k$-algebra.
    Write $\Zfr_A$ for the closure in $\PP(A)^3$ of the graph of the binary multiplication on $[\AA(A)^\times / \GG_m]$.
    Then:
    \begin{enumerate}
        \item Push-pull along the correspondence
        \[
            \begin{tikzcd}
                & \Zfr_A \ar[dl] \ar[dr] & \\
                \PP(A) \times \PP(A) & & \PP(A)
            \end{tikzcd}
        \]
        defines a monoidal structure $\star'_A := \star'_{[\AA(A) / \GG_m],-\Theta_{\PP(A)}}$ on $\Perf \PP(A)$.
        \item If $j': [\AA(A)^\times / \GG_m] \hookrightarrow \PP(A)$ is the inclusion, then the pushforward functor 
        \[
            j'_*: \big(\QC([\AA(A)^\times / \GG_m]), \star_{[\AA(A)^\times / \GG_m]}\big) \to \big(\QC(\PP(A)), \star'_A\big)
        \]
        is monoidal.
        \item The construction of $\star'_A$ is functorial in surjections of finite-dimensional $k$-algebras.
    \end{enumerate}
    When $A$ is commutative, ``monoidal'' may be upgraded to ``symmetric monoidal'' throughout.
\end{prop}

\begin{proof}
    We use the geometric $\EE_1$-EC setup $(\phi, j, \{ 0, \dots, \dim_k A \}, q)$ of \cref{ex:algebra_ec_setup}.
    When $A$ is commutative, we upgrade this to an $\EE_\infty$-EC setup.
    
    (1). The scheme $\Zfr_A$ agrees with the restriction to $\PP(A)^3$ of the graph of the multiplication morphism
    \[
        \PP(A) \times \PP(A) \to [\AA(A) / \GG_m].
    \]
    Thus \cref{prop:ec_geometric_computation} and \cref{ex:ec_geometric_diagonal} imply the claim.

    (2). This is a direct consequence of \cref{cor:ec_extends_convolution}.
    
    (3). This follows from \cref{ex:ec_morph_fdalg} and \cref{prop:ec_functoriality}.
\end{proof}

\subsection{Computations}

One can compute EC products $\Fsc \star'_A \Gsc$ using the equivalence of \cref{thm:ec_comparison}.
More precisely, using \cref{rmk:action_grading} with $R = \Sym A^\vee$, the category $\Qsf_S(\phi)\op$ is given by
\[
    \begin{tikzcd}
        q_0 & \lar["A^\vee", swap] q_1 & \lar["A^\vee", swap] \dots &\lar["A^\vee", swap] q_n
    \end{tikzcd}
\]
where we have relations $\alpha_1 \alpha_2 = \alpha_2 \alpha_1 \in \Sym A^\vee$ for $\alpha_i \in A^\vee$.
Via the equivalence of \cref{thm:ec_comparison}, $\Fsc \in \QC(\Perf \PP(A))$ corresponds to the derived $\Qsf_S(\phi)\op$-representation\footnote{Recall that our convention is that all the functors are implicitly derived.}
\[
    \begin{tikzcd}
        \Gamma(\PP(A), \Fsc) & \lar  \Gamma(\PP(A), \Fsc(-1)) & \lar \dots &\lar \Gamma(\PP(A), \Fsc(-n))
    \end{tikzcd}
\]
and similarly for $\Gsc$.
The EC product $\Fsc \star'_A \Gsc$ then corresponds to the quiver tensor product
\[
    \begin{tikzcd}
        \Gamma(\PP(A), \Fsc) \otimes_k \Gamma(\PP(A), \Gsc) & \lar \dots &\lar \Gamma(\PP(A), \Fsc(-n)) \otimes_k \Gamma(\PP(A), \Gsc(-n)),
    \end{tikzcd}
\]
where the coalgebra structure on $A^\vee$ is used to construct the tensor product of morphisms.

This is easiest to understand for specific classes of sheaves:

\begin{ex} \label{ex:Pn_simple_representation}
    For $i \in \{0, \dots, \dim A - 1\}$, the sheaf $\Omega^i_{\PP(A)}(i)[i]$ corresponds to the simple $\Qsf_S(\phi)\op$-representation
    \[
        \begin{tikzcd}
            0 & \lar \dots &\lar 0 & \lar k &\lar 0 & \lar \dots &\lar 0
        \end{tikzcd}
    \]
    sending $q_i$ to $k$ and $q_j$ to $0$ for $j \neq i$ by Bott's formula.
    In particular, we see that
    \[
        \Omega^i_{\PP(A)}(i)[i] \star'_A \Omega^j_{\PP(A)}(j)[j] = \begin{cases}
            0 & i \neq j \\
            \Omega^i_{\PP(A)}(i)[i] & i = j 
        \end{cases}
    \]
\end{ex}

\begin{ex} \label{ex:Pn_skyscraper}
    For $[a] \in \PP(A)$, the skyscraper sheaf $k([a])$ corresponds to the quiver representation
    \[
        \begin{tikzcd}
            k & \lar["a", swap] k & \lar["a", swap] \dots & \lar["a", swap] k
        \end{tikzcd}
    \]
    where we view $a$ as an element of $A^{\vee \vee} = \Hom_k(A^\vee, \Hom_k(k, k))$.
\end{ex}

\begin{prop} \label{prop:algebra_ec_skyscraper}
    Let $a_1, a_2 \in A \setminus \{ 0\}$, so $k([a_1]), k([a_2]) \in \Perf \PP(A)$.
    Then
    \[
        k([a_1]) \star'_A k([a_2]) = \begin{cases}
            k\big([\mu_A(a_1, a_2)]\big) & \mu_A(a_1, a_2) \neq 0. \\
            \oplus_{i=0}^{\dim A-1} \Omega^i_{\PP(A)}(i)[i] & \mu_A(a_1, a_2) = 0.
        \end{cases}
    \]
\end{prop}

\begin{proof}
    Using \cref{ex:Pn_skyscraper} and the fact that the coalgebra structure on $A^\vee$ is the dual of the algebra structure on $A$, we see that $k([a_1]) \star'_A k([a_2])$ corresponds to the quiver representation
    \[
        \begin{tikzcd}[column sep=large]
            k & \lar["{\mu_A(a_1, a_2)}", swap] k & \lar["{\mu_A(a_1, a_2)}", swap] \dots & \lar["{\mu_A(a_1, a_2)}", swap] k
        \end{tikzcd}
    \]
    If $\mu_A(a_1, a_2) \neq 0$, this representation corresponds to $k\big([\mu_A(a_1, a_2)]\big)$, establishing the first case of the claim.
    
    Otherwise $\mu_A(a_1, a_2) = 0$, so $k([a_1]) \star'_A k([a_2])$ corresponds to the quiver representation 
    \[
        \begin{tikzcd}
            k & \lar["0", swap] k & \lar["0", swap] \dots & \lar["0", swap] k
        \end{tikzcd}
    \]
    This decomposes as a direct sum of the simple representations of $\Qsf_S(\phi)\op$.
    \cref{ex:Pn_simple_representation} lets us convert these quiver representations back into perfect complexes, giving
    \[
        k([a_1]) \star'_A k([a_2]) = \bigoplus_{i=0}^{\dim A-1} \Omega^i_{\PP(A)}(i)[i]. \qedhere
    \]
\end{proof}

\subsection{Categorical compactifications} \label{sub:compactifications}

The monoidal $\infty$-categories $\big(\Perf \PP(A), \star'_A\big)$ give categorified ``compactifications'' of many well-known algebras and groups.
We formalize this notion as follows.

\begin{dfn}
    Let $\Nfr$ be a perfect $\EE_n$-monoid stack over $k$.
    A \emph{categorical compactification} of $\Nfr$ is a fully faithful, $\EE_n$-monoidal functor
    \[
        \iota: (\QC(\Nfr), \star_\Nfr) \hookrightarrow (\Indsf \Cc, \otimes_\Cc)
    \]
    where $(\Cc, \otimes_\Cc) \in \CAlg\big(\Cat_k^{\perf}\big)$.
    We will often abuse notation and refer to $(\Cc, \otimes_\Cc)$ as a categorical compactification of $\Nfr$.
\end{dfn}

By \cref{cor:ec_extends_convolution}, if $\Nfr$ is an open submonoid of $[\AA(A) / G]$ such that $\Nfr \subset \PP(A)$, then pushforward along the inclusion $j': \Nfr \subset \Xfr$ exhibits $(\Perf \PP(A), \star'_{\Mfr,S})$ as a categorical compactification of $\Nfr$.
We obtain the following examples in this way.

\begin{ex} \label{ex:compactification_algebra}
   Any finite-dimensional $k$-algebra $A'$ admits a categorical compactification.
   In fact, if we let $A = A' \times k$, then the inclusion $\AA(A') \hookrightarrow \PP(A)$ exhibits $(\Perf \PP(A), \star'_{\Mfr,S})$ as a categorical compactification of $\AA(A)$.
\end{ex}

\begin{ex} \label{ex:compactification_groups}
    We may also construct categorical compactifications of many linear algebraic groups.
    We list some such groups $G$ together with algebras $A$ such that $(\Perf \PP(A), \star'_{\Mfr,S})$ is a categorical compactification of $G$:
    \begin{itemize}
        \item For $G = \GG_m^n$, we can take $A = k^{n+1}$.
        \item For $G = \GL_n$, we can take $A = \End_k(k^n) \times k$.
        \item For $G = \mathrm{PGL}_n$, we can take $A = \End_k(k^n)$. 
        \item For $G = \GG_a$, we can take $A = k[\epsilon] / \epsilon^2$.
        \item For $G = B_n$ the group of invertible upper triangular matrices, let $\ol{B}_n$ be the $k$-algebra of all upper triangular matrices.
        Then we can take $A = \ol{B}_n \times k$.
        \item For $G = B_n / \GG_m$ (a Borel subgroup of $\mathrm{PGL}_n$), we can take $A = \ol{B}_n$.
    \end{itemize}
\end{ex}

\begin{rmk}
    The question of which algebraic groups admit categorical compactifications remains open in general.
    Categorical compactifications are typically far from unique, even if we impose ``connectivity'' hypotheses on $\Cc$.\footnote{For example, any full strong exceptional collection on a smooth complete toric variety gives a categorical compactification of the dense torus by \cref{ex:toric_geometric_ec_setup}.}
    It would be interesting to understand whether one can construct a \emph{canonical} categorical compactification of a nice (e.g.\ semisimple) algebraic group.
\end{rmk}

Although the EC products in this subsection can be defined geometrically, they may exhibit surprising behavior at infinity.

\begin{ex} \label{ex:multiply_zero_by_infinity}
    Consider $\PP^1 = \PP(k^2)$, so that $\Perf \PP^1$ is a categorical compactification of $\GG_m$.
    By \cref{prop:algebra_ec_skyscraper}, we see that $k([1 : 0]) \star'_{k^2} k([0 : 1])$ is the (globally supported) complex $\Osc_{\PP^1} \oplus \Osc_{\PP^1}(-1)[1]$.
    Loosely: we have extended the multiplication of $\GG_m$ to $\PP^1$, but multiplying zero by infinity produces a perfect complex rather than a number!
\end{ex}

Let us mention some additional subtleties that may arise when attempting to apply these methods over non-algebraically closed fields.
 
\begin{rmk}
    Suppose for the context of this remark that $k$ is not algebraically closed, and let $k \subset F$ be a finite extension of fields.
    In this case, we still obtain a symmetric monoidal EC product on $\Perf \PP(F)$.
    One might na\"ively expect that this EC product agrees with a genuine convolution product, with the multiplication on $\PP(F)$ constructed by restricting the multiplication on $\AA(F)$ to $\AA(F) \setminus \{ 0 \}$.
    
    However, such a construction is typically not possible: even though the set of $k$-points $\big(\AA(F) \setminus \{ 0 \}\big)(k) = F^\times$ is closed under multiplication, this is no longer true for the set of $F$-points $\big(\AA(F) \setminus \{ 0 \}\big)(F)$.
    In particular, the construction of the EC product on $\Perf \PP(F)$ does not contradict the fact that all complete connected algebraic groups are abelian varieties.
\end{rmk}

\subsection{Recovering $A$ from $(\Perf \PP(A), \star'_A)$}

We have the following result allowing us to reconstruct a finite-dimensional algebra $A$ from the corresponding EC product on $\Perf \PP(A)$:

\begin{prop} \label{prop:algebra_reconstruction}
    Let $A$ and $A'$ be finite-dimensional $k$-algebras.
    Then $A$ and $A'$ are isomorphic as $k$-algebras if and only if there is a monoidal equivalence $(\Perf(\PP(A)),\star_A') \simeq (\Perf(\PP(A')),\star_{A'}')$.
\end{prop}

To prove \cref{prop:algebra_reconstruction}, we will need a linear-algebraic lemma.

\begin{notn}
    Let $V$ and $V'$ be finite-dimensional vector spaces.
    Given an element $v \in V$ (possibly zero), let $[v]$ denote the corresponding class in $V / k^\times$. 
    Given a linear map $\phi: V \to V'$, let $[\phi]$ denote the corresponding class in $V/k^\times \to V/k^\times$. 
\end{notn}

\begin{lem}\label{lemma-findimalgebra is determined by projective monoid}
    Let $A$ and $A'$ be finite-dimensional $k$-algebras.
    Suppose $\phi: A \to A'$ is a linear map such that $[\phi]: A / k^\times \to {A'} / k^\times$ is a homomorphism of monoids (in $\Set$).
    Then there exists $c \in k^\times$ such that $c\phi: A \to A'$ is a monoid homomorphism.
\end{lem}

\begin{proof}
    Invertible elements of $A$ are dense in $\AA(A)$ (see e.g.\ \cite[page 1]{vinberg1995reductive}), so we may fix a basis $\{ e_0, \dots, e_n \}$ of $A$ such that $e_i$ is invertible in $A$ for all $i$.
    Because $[\phi]$ is a monoid homomorphism, the images $[\phi]([e_i])$ are invertible in $A' / k^\times$.
    Thus the images $\phi(e_i)$ must also be invertible in $A'$.
    In particular, for fixed $i$, the sets $\{ \mu_{A'}(\phi(e_i), \phi(e_j)) \}_{j=0}^n$ are linearly independent in $A'$.
    
    For $i, j, \ell \in \{ 0, \dots, n \}$, we may write
    \[
        \phi\big(\mu_A(e_i,e_j)\big) = c_{i,j}\mu_{A'}\big(\phi(e_i),\phi(e_j)\big)
    \]
    and
    \[
        \phi\big(\mu_A(e_i, e_j + e_\ell)\big) = c_{i,j\ell}\mu_{A'}\big(\phi(e_i),\phi(e_j + e_\ell)\big)
    \]
    for some $c_{i,j}, c_{i,j\ell} \in k^\times$.
    For any $i, j, \ell \in \{ 0 , \dots, n \}$, we get
    \begin{align*}
        c_{i,j}\mu_{A'}\big(\phi(e_i),\phi(e_j)\big) + c_{i,\ell}\mu_{A'}\big(\phi(e_i),\phi(e_\ell)\big) &= \phi\big(\mu_A(e_i, e_j+e_\ell)\big) \\
        &= c_{i,j\ell}\mu_{A'}\big(\phi(e_i),\phi(e_j+e_\ell)\big) \\
        &= c_{i,j\ell}\mu_{A'}\big(\phi(e_i),\phi(e_j)\big) + c_{i,j\ell}\mu_{A'}\big(\phi(e_i),\phi(e_\ell)\big)
    \end{align*}
    so that $c_{i,j} = c_{i,j\ell} = c_{i,\ell}$ by the aforementioned linear independence.
    In particular, for all $i, j$, we have $c_{i,j} = c_{i,0}$.
    Repeating the argument with the order of inputs to $\mu_{A'}$ reversed shows that $c_{i,j} = c_{0,j}$.
    Thus $c_{i,j} = c_{0,0}$ for all $i, j$, i.e.\
    \[
        \phi\big(\mu_A(e_i, e_j)\big) = c_{0,0} \mu_{A'}\big(\phi(e_i), \phi(e_j)\big).
    \]
    Multiplying both sides of this equation by $c_{0,0}$ shows that $(c_{0,0} \phi) \circ \mu_A = \mu_{A'} \circ (c_{0,0} \phi \otimes c_{0,0} \phi)$, i.e.\ $c_{0,0} \phi$ is a monoid homomorphism.
\end{proof}

\begin{proof}[Proof of \cref{prop:algebra_reconstruction}]
    The ``only if'' direction follows by the functoriality of the construction of $\star'$ on surjective maps (\cref{prop:ec_functoriality}).
    
    For the ``if'' direction, we may assume $A$ and $A'$ have the same underlying vector space $V$.
    Every monoidal equivalence $(\Perf(\PP(A)),\star_A') \simeq (\Perf(\PP(A')),\star_{A'}')$ is induced by an autoequivalence $\tau$ of $\Perf (\PP(V))$.
    By \cite[Theorem 3.1]{bondal2001reconstruction}, we can write $\tau = [\phi]_* (-\otimes \Lsc [n])$ for some $\phi \in \GL(V)$ (so $[\phi] \in \PGL(V)$), $\Lsc \in \Pic(\PP(V))$, and $n \in \ZZ$.
    As $\tau(k([1_A])) = k([1_{A'}])$, we must have $n = 0$, i.e.\ $\tau = [\phi]_* (- \otimes \Lsc)$.\footnote{One can also show that in this situation we must have $\Lsc \simeq \Osc$, though this does not simplify matters for us.}

    By \cref{lemma-findimalgebra is determined by projective monoid}, it suffices to show that $[\phi]: A / k^\times \to A' / k^\times$ is a monoid homomorphism. 
    To this end, let $x, y \in A$.
    We may assume $\dim A \geq 2$, allowing us to argue by cases:
    \begin{itemize}
        \item If $x = 0$ or $y = 0$, we must have 
        \begin{equation} \label{eq:monoid_hom_zero}
            [\phi](\mu_A(x,y)) = [0] = \big[\mu_{A'}\big(\phi(x), \phi(y)\big)\big].
        \end{equation}
        \item If $x$ and $y$ are both nonzero and $\mu_A(x,y) = 0$, then $\mu_{A'}(\phi(x),\phi(y)) = 0$, so \cref{eq:monoid_hom_zero} still holds.
        In fact, if we had $\mu_{A'}(\phi(x),\phi(y))  \neq 0$, then 
        \[
            k([\phi(x)]) \star'_{A'} k([\phi(y)]) = k\big([\mu_{A'}(x, y)]\big)
        \]
        is indecomposable while \cref{prop:algebra_ec_skyscraper} implies $k([x]) \star'_A k([y])$ is decomposable.
        But then
        \[
            \tau\big(k([x]) \star'_A k([y])\big) = k([\phi(x)]) \star'_{A'} k([\phi(y)])
        \]
        contradicts the assumption that $\tau$ is an equivalence.
        
        \item Otherwise, $\mu_A(x,y) \neq 0$, so 
        \[
            k\big(\big[ \phi(\mu_A(x,y)) \big]\big) \simeq \tau\big(k([\mu_A(x,y)])\big) \simeq \tau\big(k([x]) \star'_A k([y])\big) \simeq k([\phi(x)]) \star'_{A'}k([\phi(y)]).
        \]
        Hence, by \cref{prop:algebra_ec_skyscraper}, we have $[\phi(\mu_A(x,y))] = [\mu_{A'}(\phi(x),\phi(y))]$.
    \end{itemize}
    Thus $[\phi]:A/k^\times \to {A'}/k^\times$ is a monoid homomorphism.
\end{proof}

\subsection{Invariants}

Fix a finite-dimensional algebra $A$ and consider the monoidal category $\big(\Perf \PP(A), \star'_A\big)$.
In this section we compute:
\begin{itemize}
    \item The Balmer spectrum $\Spc_{\star'_A} \PP(A)$
    \item The Grothendieck ring $K^0\big(\Perf \PP(A), \star'_A\big)$.
    \item The Picard group $\Pic(\Perf \PP(A), \star'_A)$.
\end{itemize}
We observe that the Balmer spectrum and the Grothendieck ring depend only on $\dim A$.

We first compute $\Spc_{\star'_A} \PP(A)$.
Here we follow \cite{nakano2022noncommutative} for a definition of the Balmer spectrum of a stably monoidal $\infty$-category which need not be symmetric monoidal.

\begin{dfn}
    Let $(\Cc, \otimes_\Cc)$ be a small stably monoidal $\infty$-category.
    A two-sided thick $\otimes$-ideal $\Pc$ in $\Cc$ is \emph{prime} if $\Pc \neq \Cc$ and, whenever we have $\Ic \otimes_\Cc \Jc \subset \Pc$ for two-sided thick $\otimes$-ideals $\Ic,\Jc$ of $\Cc$, then either $\Ic \subset \Pc$ or $\Jc \subset \Pc$.
    The \emph{Balmer spectrum} of $(\Cc, \otimes_\Cc)$ is
    \[
        \Spc_{\otimes_\Cc} \Cc = \bset{\Pc}{\textrm{$\Pc$ is a prime two-sided thick $\otimes$-ideal of $\Cc$}}
    \]
    with a topology defined similarly to the (Zariski) topology of the usual Balmer spectrum (see \cite[Section 1.2]{nakano2022noncommutative} for details). 
\end{dfn}

\begin{lem} \label{lem:spectrum_computation} 
    Let $(\Cc, \otimes_\Cc)$ be a small stably monoidal $\infty$-category.
    Let $S_0, \dots, S_n$ be a collection of exceptional objects of $\Cc$ with $\langle S_i \rangle \neq \langle S_j \rangle$ for $i\neq j$, and for each $i$, let $\Pc_i = \bangle{S_j}{j \neq i}$.
    Suppose that:
    \begin{enumerate}
        \item $\Cc = \langle S_0,\dots,S_n\rangle$
        \item $\Pc_i = \ker (-\otimes S_i) = \ker (S_i \otimes -)$ for all $i$.
        \item $\langle S_i \rangle$ is a two-sided thick $\otimes$-ideal for all $i$. 
    \end{enumerate}
    Then
    \[
    \Spc_{\otimes_\Cc} \Cc = \bigsqcup_{i = 0}^n \Pc_i. 
    \]
\end{lem}

\begin{proof}
    First, $\Pc_i$ is clearly a two-sided thick $\otimes$-ideal since $\Pc_i = \ker (-\otimes S_i)= \ker (S_i \otimes -)$. 
    Because $\Cc/\Pc_i = \langle S_i \rangle = \Perf(k)$, we see that $\Pc_i$ is maximal among thick subcategories of $\Cc$.
    Thus $\Pc_i$ is a prime two-sided thick $\otimes$-ideal of $\Cc$ by \cite[Theorem 3.2.3]{nakano2022noncommutative}.

    Conversely, let $\Pc$ be a prime thick $\otimes$-ideal of $\Cc$.
    Because $\Pc \neq \Cc$, there exists $i$ such that $S_i \not \in \Pc$.
    Then, since $\langle S_i \rangle \otimes \langle S_j \rangle = \langle 0 \rangle \subset \Pc$ for each $j \neq i$ and $\langle S_i \rangle \not \subset \Pc$, we have $S_j \in \langle S_j \rangle \subset \Pc$ for every $j \neq i$.
    Thus, $\Pc = \Pc_i$ as $\Pc_i$ is maximal. 
\end{proof}

\begin{prop}\label{prop:balmer_spectrum}
    Let $A$ be a finite dimensional algebra $A$ over $k$.
    Then 
    \[
        \Spc_{\star'_A} \PP(A) = \bigsqcup_{i = 0}^{\dim A-1} \bangle{\Omega^j(j)[j]}{j \neq i}.
    \]
    In particular, $\Spc_{\star'_A} \PP(A)$ depends only on the dimension of $A$. 
\end{prop}

\begin{proof}
    This follows from \cref{lem:spectrum_computation} by taking $S_i = \Omega^i(i)[i]$ for $i = 0, \dots, \dim A - 1$ and using the computation of \cref{ex:Pn_simple_representation}. 
\end{proof} 

\begin{rmk}
    Recent work on higher Zariski geometry (\cite{aoki2025higher}) proves that the $\infty$-category of $2$-rings (i.e., rigid stably symmetric monoidal $\infty$-categories) embeds into the $\infty$-category of $2$-ringed spaces via an enhanced version of the Balmer spectrum construction.
    If $A$ is commutative, the above computations can be used to show that the $2$-ringed space associated with $(\Perf(\PP(A)), \star'_A)$ depends only on $\dim A$, i.e.\ this $2$-ringed space does not distinguish between different commutative algebras of the same dimension.
    However, it is easy to show directly that $\star'_A$ is not rigid, so we do not contradict the results of \cite{aoki2025higher}. 
\end{rmk}

\cref{ex:Pn_simple_representation} also allows us to compute the Grothendieck ring of $(\Perf \PP(A), \star'_A)$:

\begin{prop}
    Let $A$ be a finite dimensional $k$-algebra.
    Then there is an isomorphism
    \[
        K^0\big(\Perf \PP(A), \star'_A\big) \simeq \ZZ^{\dim A}.
    \]
\end{prop}

\begin{proof}
    Consider the the basis of $K^0(\Perf \PP(A), \star'_A)$ given by $v_i: = \big[\Omega^i(i)[i]\big]$ for $i = 0,\dots, \dim A - 1$.
    By \cref{ex:Pn_simple_representation}, regardless of the algebra structure of $A$, we have $v_i^2 = v_i$ for all $i$ and $v_i \cdot v_j = 0$ for all $i\neq j$.
    Thus the basis $\{ v_i \}_{i=0}^{\dim A - 1}$ gives rise to the desired isomorphism.
\end{proof}

Unlike the Balmer spectrum and the Grothendieck ring, the Picard group of $\Pic\big(\Perf \PP(A), \star_A\big)$ can depend on the choice of $A$.

\begin{prop}
    Let $A$ be a finite dimensional $k$-algebra.
    Then
    \[
    \Pic\big(\Perf \PP(A), \star_A\big) = (A^\times)^{\dim A - 1} / k^\times \times \bb Z.  
    \]
\end{prop}

\begin{proof}
    Suppose $\Fsc$ is invertible in $\Pic\big(\Perf \PP(A), \star_A\big)$.
    Then the derived quiver representation corresponding to $\Fsc$ must take invertible values at each vertex, i.e.\ $\Fsc$ corresponds to a derived quiver representation of the form
    \[
        \begin{tikzcd}
            k[i_0] & \lar k[i_1] & \lar \dots & \lar k[i_{\dim A - 1}].
        \end{tikzcd}
    \]

    We first claim that $i_0 = i_1 = \dots = i_{\dim A - 1}$.
    Indeed, if this were not the case, then we would necessarily have $i_p \neq i_{p-1}$ for some $p$.
    This forces $H^0\big(\Hom_k(k[i_p], k[i_{p-1}])\big) = \Ext_k^{i_{p-1} - i_p}(k, k)= 0$.
    Thus, for any $\Gsc \in \Perf \PP(A)$, the derived quiver representation corresponding to $\Fsc \star'_A \Gsc$ would necessarily have $0$ as one of its morphisms.
    In particular, because $\Fsc$ is invertible, one of the morphisms in the quiver representation corresponding to the unit of $\star'_A$ is $0$.
    But this is impossible as the unit of $\star'_A$ is $k(1_A)$, corresponding to the quiver representation
    \[
        \begin{tikzcd}
            k & \lar["1_A", swap] k & \lar["1_A", swap] \dots & \lar["1_A", swap] k.
        \end{tikzcd}
    \]

    Thus every invertible object of $\Pic\big(\Perf \PP(A), \star_A\big)$ corresponds to a derived quiver representation of the form
    \[
        \begin{tikzcd}
            k[i] & \lar k[i] & \lar \dots & \lar k[i]
        \end{tikzcd}
    \]
    for some $i \in \ZZ$.
    We may define a homomorphism $\alpha: \Pic\big(\Perf \PP(A), \star_A\big) \to \ZZ$ sending a quiver representation of the above form to $i$.
    The homomorphism 
    \begin{align*}
        \ZZ &\to \Pic\big(\Perf \PP(A), \star_A\big) \\
        i &\mapsto k(1_A)[i]
    \end{align*}
    defines a section of $\alpha$.
    
    It remains to show that $\ker \alpha = (A^\times)^{\dim A - 1} / k^\times$.
    We may define a map $\beta: (A^\times)^{\dim A - 1} \to \ker \alpha$ sending $(a_1, \dots, a_{\dim A - 1})$ to the object of $\ker \alpha$ corresponding to the quiver representation
    \[
        \begin{tikzcd}
            k & \lar["a_1", swap] k & \lar["a_2", swap] \dots & \lar["a_{\dim A - 1}", swap] k
        \end{tikzcd}
    \]
    It is clear that every object of $\ker \alpha$ arises in this way.
    Furthermore, any isomorphism $\beta(a_1, \dots, a_{\dim A - 1}) \simeq \beta(a'_1, \dots, a'_{\dim A - 1})$ is witnessed by a commutative diagram
    \[
        \begin{tikzcd}
            k \dar["\sim"] & \lar["a_1", swap] k \dar["\sim"] & \lar["a_2", swap] \dots & \lar["a_{\dim A - 1}", swap] k \dar["\sim"] \\
            k & \lar["a'_1", swap] k & \lar["a'_2", swap] \dots & \lar["a'_{\dim A - 1}", swap] k,
        \end{tikzcd}
    \]
    where invertibility of each $a_p$ and each $a'_p$ ensures that the vertical arrows are all determined by the rightmost vertical arrow.
    Thus $\ker \beta = k^\times$ and $\ker \alpha = (A^\times)^{\dim A - 1} / \ker \beta = (A^\times)^{\dim A - 1} / k^\times$.
\end{proof}

\begin{rmk}
    The invariants we have discussed in this section reflect the differences in behavior between the EC products on $\Perf \PP^d$ and the usual tensor product on $\Perf \PP^d$.
    Indeed, for the usual tensor product $\otimes_\Osc$, we have:
    \begin{itemize}
        \item $\Spc_{\otimes_\Osc} \Perf \PP^d = \PP^d$.
        \item $K^0\big(\Perf \PP^d, \otimes_\Osc\big) \simeq \ZZ[\eta] / \eta^{n+1}$ where $\eta = [\Osc_{\PP^d}(1)] - 1$.
        \item $\Pic(\Perf \PP(A), \otimes_\Osc) \simeq \ZZ \cdot \Osc_{\PP^d}(1) \times \ZZ$ (where the second factor arises from the shift functor $[1]$).
    \end{itemize}
\end{rmk}

\begin{ex}
    We may use the above results to construct two distinct symmetric monoidal structures on $\Perf \PP^1$ with the same Balmer spectrum and the same ring structure on $K^0$.
    Namely, take $A = k^2$ and $A' = k[\epsilon] / \epsilon^2$.
    The Balmer spectra and Grothendieck rings of $\big(\Perf \PP(A), \star'_{A}\big)$ and $\big(\Perf \PP(A'), \star'_{A'}\big)$ agree because $\dim A = \dim A'$.
    However, $\Pic \big(\Perf \PP(A), \star'_{A}\big) = k^\times \times k^\times \times \ZZ$, while $\Pic \big(\Perf \PP(A'), \star'_{A'}\big) = k \times k^\times \times \ZZ$, so the symmetric monoidal structures are distinct.
\end{ex}

Note that the Picard group is not a complete invariant.

\begin{ex}
    Let $A$ be any finite dimensional $k$-algebra.
    In general, $A$ is not isomorphic to its opposite algebra $A\op$.
    However, we have
    \[
        \Pic (\Perf \PP(A), \star_A) = (A^\times)^{\dim A - 1} / k^\times \times \ZZ \cong ((A\op)^\times)^{\dim A - 1} / k^\times \times \ZZ = \Pic (\Perf \PP(A\op), \star_{A\op})
    \]
    where we use the isomorphism $A^\times \cong (A\op)^\times = (A^\times)\op$ defined by $a \mapsto a\inv$.
    Note that this isomorphism does not extend to a linear map $A \to A\op$, so the argument of \cref{prop:algebra_reconstruction} does not apply.
\end{ex}

\begin{ex}
    For an example where all algebras involved are commutative, let $A = k[\epsilon]/\epsilon^3$ and $A' = k[x,y]/\langle x^2,xy,y^2 \rangle$, then $A^\times \iso k^\times \times k^2 \iso {A'}^\times$, where the first isomorphism is given by 
    \begin{align*}
    A^\times = k^\times \times (1 + \fr m_A) &\to k^\times \times k^2 \\
    (c, 1 + a\epsilon + b \epsilon^2) &\mapsto (c, a, b - a^2/2).
    \end{align*}
    As before, this isomorphism is not linear, so this does not contradict the argument of \cref{prop:algebra_reconstruction}.
\end{ex}

\section{Application: Tensor products in toric mirror symmetry} \label{sec:hms}

Suppose $\Xfr = \Xfr_\Sigma$ is the smooth complete toric variety associated with a fan $\Sigma \subset N_\RR$ (where $N$ is a lattice, $M$ is its dual lattice, and $(-)_\RR = (-) \otimes_\ZZ \RR$).
Homological mirror symmetry for toric varieties (in this context also called the ``coherent-constructible correspondence'') gives a symmetric monoidal equivalence
\begin{equation} \label{eq:toric_hms}
    \big(\QC(\Xfr_\Sigma), \otimes_\Osc \big) \simeq \big(\Sh_{\Lambda_\Sigma}(M_\RR / M), \star_{M_\RR / M}\big)
\end{equation}
where:
\begin{itemize}
    \item $M_\RR / M$ is a real torus,
    \item $\Lambda_\Sigma \subset T^*(M_\RR / M)$ is a certain Lagrangian defined from the combinatorics of $\Sigma$,
    \item $\Sh_{\Lambda_\Sigma}(M_\RR / M)$ is the category of constructible \footnote{Here we use constructible in a weak sense, i.e.\ we do not require any finiteness conditions on stalks.
    See \cite{kuwagaki2017nonequivariant} for a more comprehensive discussion of finiteness conditions in this context.} sheaves on $M_\RR / M$ with coefficients in $k$ and singular support in $\Lambda_\Sigma$, and
    \item $\star_{M_\RR / M}$ is the convolution product on $T^*(M_\RR / M)$.
\end{itemize}
There is a vast literature on this subject, and we will content ourselves by noting that a fairly complete (and significantly more general) discussion of \cref{eq:toric_hms} may be found in \cite{kuwagaki2017nonequivariant}, which builds on \cite{bondal2006derived, fang2011categorification} and many other sources.

In this section we will discuss the appearance of quiver tensor products in homological mirror symmetry, with a particular focus on smooth complete toric varieties of Bondal-Ruan type (cf.\ \cref{dfn:bondal_thomsen}).
We begin by noting the following explicit description of quiver tensor products in the Bondal-Ruan case.

\begin{prop} \label{mainthm_toric}
    Let $\Xfr$ be a smooth complete toric variety of Bondal-Ruan type, and let $\Xfr \subset \Mfr$ be the Cox presentation of $\Xfr$.
    Write $\Zfr_{\mu,-\Theta_\Xfr}$ for the closure in $\Xfr^3$ of the graph of the binary multiplication on the dense torus in $\Xfr$.
    Then push-pull along the correspondence
    \[
        \begin{tikzcd}
              & \Zfr_{\mu,-\Theta_\Xfr} \ar[dl] \ar[dr] & \\
            \Xfr \times \Xfr & & \Xfr
        \end{tikzcd}
    \]
    defines a symmetric monoidal structure $\star'_{\Mfr,-\Theta_\Xfr}$ on $\Perf(\Xfr)$.
    Furthermore, there is a symmetric monoidal equivalence
    \[
        \big(\Perf(\Xfr), \star'_{\Mfr,-\Theta_\Xfr}\big) \simeq \Big(\Fun\big(\Qsf_{\Theta_\Xfr}^\pre(\phi), \Perf(k)\big), \otimes_\Qsf\Big).
    \]
\end{prop}

\begin{proof}
    We use the geometric EC setup $(\phi, j, -\Theta_\Xfr, q)$ of \cref{ex:toric_geometric_ec_setup}.
    The first claim follows by combining \cref{prop:ec_geometric_computation} and \cref{ex:ec_geometric_diagonal}.
    The second claim is \cref{thm:ec_comparison}, noting that $\Qsf_{-\Theta_\Xfr}^\pre(\phi)\op \simeq \Qsf_{\Theta_\Xfr}^\pre(\phi)$  by taking duals of line bundles.
\end{proof}

\subsection{The mirror of the constructible tensor product}

Given the underlying equivalence of categories of \cref{eq:toric_hms}, it is natural to ask when and whether one can describe the tensor product $\otimes_k$ of constructible sheaves on $M_\RR / M$ in terms of the algebraic geometry of $\Xfr_\Sigma$.
This question does not always make sense -- for $\Fsc, \Gsc \in \Sh_{\Lambda_\Sigma}(M_\RR / M)$, the singular support of the tensor product $\Fsc \otimes_k \Gsc$ need not lie in $\Lambda_\Sigma$.
However, when the Lagrangian $\Lambda_\Sigma$ arises from a stratification $Z = \{ Z_i \}_{i \in I}$ of $M_\RR / M$, the category $\Sh_{\Lambda_\Sigma}(M_\RR / M)$ is closed under $\otimes_k$, and the question does make sense.
In this case, we shall write $\Sh_Z(M_\RR / M) := \Sh_{\Lambda_\Sigma}(M_\RR / M)$.

By \cite{bondal2006derived} (see also \cite[\S 5]{favero2025homotopy} for a more detailed presentation and generalization), we know that $\Lambda_\Sigma$ arises from a stratification when $\Xfr$ is of Bondal-Ruan type.
More precisely, let $F_\phi: M_\RR / M \to \weight(G)$ be the anti-Bondal-Ruan map of \cref{dfn:bondal_thomsen}.
We may decompose $M_\RR / M$ as a disjoint union of the level sets $\{ Z_\chi := F_\phi\inv(\chi) \}_{\chi \in -\Theta_\Xfr}$.\footnote{In \cite{favero2025homotopy}, $Z_\chi$ is denoted by $S_{\chi}$ for $\chi \in \hat G = \weight(G)$. We use the letter $Z$ only to avoid conflict with our notation for transparent collections.}

There is a natural order on $-\Theta_\Xfr$ given by the transitive closure of the following rule: $\chi_1 \leq \chi_2$ if $Z_{\chi_1} \subset \ol{Z_{\chi_2}}$ (cf.\ \cite[Corollary 5.7 and Definitions 4.8, 4.13, and 4.28]{favero2025homotopy}.
The order on $-\Theta_\Xfr$ may also be understood algebraically: by the discussion at the beginning of \cite[\S 5.2]{favero2025homotopy}, we have $\chi_1 \leq \chi_2$ if and only if $H^0(\Xfr, \Osc(\chi_2 - \chi_1)) \neq 0$.
In particular, $Z = \{ Z_\chi \}_{\chi \in -\Theta_\Xfr}$ is a (non-conical) stratification, which we call the \emph{anti-Bondal-Ruan stratification}, of $M_\RR / M$ by the poset $-\Theta_\Xfr$ (cf. \cite[Definition 4.2]{favero2025homotopy}).

\begin{prop} \label{prop:hms_tensor}
    Let $\Xfr$ be a smooth complete toric variety of Bondal-Ruan type, and let $Z$ be the corresponding anti-Bondal-Ruan stratification of $M_\RR / M$.
    There is a symmetric monoidal equivalence
    \[
        \big(\QC(\Xfr), \star'_{[\AA^n/G],-\Theta_\Xfr} \big) \simeq \big(\Sh_Z(M_\RR / M), \otimes_k\big)
    \]
\end{prop}

\begin{proof}
   Let $\Ex_Z(M_\RR / M)$ be the $\infty$-category of \emph{exit paths} of the stratification $Z$.
   We refer to \cite[Proposition 2.2.10]{haine2024exodromy} for a precise definition of $\Ex_Z(M_\RR / M)$, but we note as a heuristic (cf.\ \cite[Observation 5.1.6]{haine2024exodromy}) that:
   \begin{itemize}
        \item Objects of $\Ex_Z(M_\RR / M)$ are points of $Z$.
        \item Morphisms are \emph{exit paths} in $M_\RR / M$, i.e.\ paths $\gamma: [0, 1] \to M_\RR / M$ (from the domain to the codomain) such that if $t_1 \leq t_2$ and $\gamma(t_i) \in Z_{\chi_i}$ for $i = 1, 2$, then $\chi_1 \leq \chi_2$.
        \item Composition is concatenation of paths.
    \end{itemize}
    By \cite[Theorem 0.4.2 and Example 5.3.10]{haine2024exodromy}, there is a symmetric monoidal equivalence
    \[
        \big(\Sh_Z(M_\RR / M), \otimes_k\big) = \big(\Fun(\Ex_Z(M_\RR / M), \Dsf(k)), \otimes_k\big)
    \]
    given by sending a constructible sheaf to the collection of its stalks and specialization maps.\footnote{The symmetric monoidal structure on $\Fun(\Ex_Z(M_\RR / M), \Dsf(k))$ is just the objectwise tensor product, so saying that this equivalence is symmetric monoidal is just saying that stalks and specialization maps commute with tensor products.}
    By \cite[Proposition 5.5 and Proposition 5.6]{favero2025homotopy}, there is a natural equivalence
    \[
        \Qsf_{-\Theta_\Xfr}^\pre(\phi)\op \xrightarrow{\sim} \Ex_Z(M_\RR / M)
    \]
    given by sending an object $\Osc_\Xfr(\chi) \in \Qsf_{-\Theta_\Xfr}^\pre(\phi)\op$ to the corresponding stratum $Z_\chi$.
    Thus we obtain a chain of equivalences
    \begin{align*}
        \big(\Sh_Z(M_\RR / M), \otimes_k\big) &= \big(\Fun(\Ex_Z(M_\RR / M), \Dsf(k)), \otimes_k\big) \\
        &= \big(\Fun(\Qsf_{-\Theta_\Xfr}^\pre(\phi)\op, \Dsf(k)), \otimes_k\big) \\
        &= \big(\Dsf(\Qsf_{-\Theta_\Xfr}(\phi)\op), \otimes_{\Qsf}\big) \\
        &= \big(\QC(\Xfr), \star'_{[\AA^n/G],\Theta_\Xfr} \big). \qedhere
    \end{align*} 
\end{proof}

\begin{rmk}
    \cref{prop:hms_tensor} may also be stated for the usual Bondal-Ruan stratification, though in this case the geometric interpretation of the EC product $\star'_{[\AA^n/G],-\Theta_\Xfr}$ on $\QC(\Xfr)$ is less clear.
    Going between the Bondal-Ruan and anti-Bondal-Ruan stratifications corresponds to taking the negative of the corresponding Lagrangians.
    This does not affect the existence of the homological mirror symmetry equivalence, as both $\QC(\Xfr)$ and $\Sh_{\Lambda_\Sigma}(M_\RR / M)$ are self-dual.
\end{rmk}

When $\Xfr$ is not of Bondal-Ruan type, the situation can be more subtle:

\begin{ex} \label{ex:ms_hirzebruch}
    Let $\Xfr = \Xfr_\Sigma = F_n$ be a Hirzebruch surface of type $n \geq 2$.
    By \cite[Proposition 6.1]{king_conj}, $\Xfr$ admits a full strong exceptional collection of line bundles.
    This corresponds (by \cref{ex:toric_cox}) to a collection of weights $S$ which is transparent for an embedding $\Xfr \subset [\AA^4 / \GG_m^2]$.
    Let $\Qsf_S^\pre(\phi)$ be the discrete category of \cref{dfn:prelinear_quiver}.
    By \cref{ex:toric_coalgebra}, we obtain a symmetric monoidal structure $\otimes_{\Qsf_S(\phi)}$ on $\Dsf(\Qsf_S(\phi)\op) \simeq \Sh_{\Lambda_\Sigma}(M_\RR / M)$.
    
    However, this symmetric monoidal structure \emph{does not agree} with the tensor product of constructible sheaves on $M_\RR / M$.
    In fact, the category $\Sh_{\Lambda_\Sigma}(M_\RR / M)$ is not closed under $\otimes_k$!
    (This relates to the fact that $\Lambda_\Sigma$ is not the Lagrangian of conormals to a stratification of $M_\RR / M$.)
    The discrepancy arises because, even though $\Xfr$ has a full strong exceptional collection of line bundles, the equivalence of \cref{eq:toric_hms} is not induced by said collection.
\end{ex}

\subsection{A conjecture about Cox categories}

The recent paper \cite{ballard2024king} introduces a ``Cox category'' $\QC_{\mathrm{Cox}}(\Xfr)$ (for an arbitrary semiprojective toric variety $\Xfr$) in which the Bondal-Thomsen collection $\Theta_\Xfr$ behaves as if it were ``transparent'' in a suitable sense.
In general, the $\QC_{\mathrm{Cox}}(\Xfr)$ is not the derived category of any genuine stack $\Xfr$, though it has a geometric interpretation in terms of ``gluing birational models of $\Xfr$.''
Via homological mirror symmetry, $\QC_{\mathrm{Cox}}(\Xfr)$ is expected to correspond to a category $\Sh_{\Lambda_\mathrm{Cox}}(M_\RR / M)$ where the singular support Lagrangian $\Lambda_{\mathrm{Cox}}$ is obtained by taking the unions of the Lagrangians of these birational models.

\begin{conj} \label{conj:hms_tensor_cox}
    Let $\Xfr$ be a smooth semiprojective toric variety.
    The category $\Sh_{\Lambda_\mathrm{Cox}}(M_\RR / M)$ is closed under $\otimes_k$, and there is a symmetric monoidal equivalence
    \[
        \big(\QC_{\mathrm{Cox}}(\Xfr), \star'\big) \simeq \big(\Sh_{\Lambda_\mathrm{Cox}}(M_\RR / M), \otimes_k\big)
    \]
    where $\star'$ is a ``birationally glued EC product'' extending the convolution on the maximal torus of $\Xfr$.
\end{conj}

\appendix

\section{\texorpdfstring{$\Oc$}{O}-monoidal structures and adjoints} \label{app:monoidal_adjunctions}

In this appendix, we collect various useful results about $\Oc$-algebras, $\Oc$-monoidal structures, and $\Oc$-monoidal adjunctions.
We expect that much of this appendix is well-known to the experts and have tried to give references to the literature where possible; any failure of attribution is due to the authors' ignorance.

The results here are stated in the ``non-enriched'' (i.e., enriched in spaces) context in which they originally appeared in the literature.
In the body of the paper, we will largely use versions of these results enriched in spectra or $D(k)$.
One may use the results of \cite[Appendix A]{mgs_universal} (especially \cite[Theorem A.3.8]{mgs_universal}) to upgrade the results here to their corresponding enriched versions by regarding presentable $\infty$-categories enriched in a presentably symmetric monoidal $\infty$-category $\Vc$ as presentable $\infty$-categories with $\Vc$-actions.

Fix an $\infty$-operad $\Oc^\tens$ and write $\Oc$ for $\Oc^\otimes_{\langle 1 \rangle}$. Suppose that $\Oc^\tens$ is single-colored, i.e., $\Oc$ is contractible (e.g.\ this is true for $\Oc = \EE_n$ for some $n \in \{1, \dots, \infty\}$).\footnote{This restriction is not logically necessary, but it makes the statements cleaner.}
We begin by recalling some notation from \HA{\S 2}.

\begin{notn}\label{notn:infinity_operads}
    Let $\Cc$ and $\Dc$ be $\Oc$-monoidal categories, i.e. fix coCartesian fibrations $p_\Cc: \Cc^\tens \to \Oc^\tens$ and $p_\Dc: \Dc^\tens \to \Oc^\tens$ of $\infty$-operads with $\Cc^\otimes \times_{\Oc^\otimes} \pt = \Cc$ and $\Dc^\otimes \times_{\Oc^\otimes} \pt = \Dc$.
    Let
    \[
    \Alg_{\Cc/\Oc}(\Dc) = \Fun^{\Oc,\lax}(\Cc,\Dc)
    \]
    denote the full subcategory of $\Fun_{\Oc^\otimes}(\Cc^\otimes,\Dc^\otimes)$ spanned by the maps of $\infty$-operads.
    By abuse of notation, a lax $\Oc$-monoidal functor $F: \Cc \to \Dc$ is a functor $F: \Cc \to \Dc$ together with $F^\tens \in \Fun_{\Oc}^\lax(\Cc,\Dc)$ inducing $F$ on the underlying $\infty$-categories.
    When $\Cc = \Oc$ and $p_\Cc = \id_\Oc$, we write $\Alg_{/\Oc}(\Dc) := \Alg_{\Oc/\Oc}(\Dc)$.
    
    When $\Oc^\otimes = \EE_\infty^\otimes \simeq N(\Fin_*)$ is the commutative operad, we write $\Alg_{\Cc}(\Dc) := \Alg_{\Cc/\EE_\infty}(\Dc)$.
    We set
    \[
    \CAlg(\Cc):= \Alg_{\EE_\infty/\EE_\infty}(\Cc) = \Alg_{\EE_\infty}(\Cc) = \Alg_{/\EE_\infty}(\Cc) 
    \]
    and set
    \[
    \Fun^{\otimes,\lax}(\Cc,\Dc) := \Fun^{\Oc,\lax}(\Cc,\Dc). 
    \]
    
    By abuse of notation, if $\Cc$ is a symmetric monoidal $\infty$-category, we say $c \in \Cc$ is an $\Oc$-algebra object if we may write $c = F(*) \in \Cc$ for some $F \in \Alg_{\Oc}(\Cc) = \Fun^{\otimes,\lax}(\Oc,\Cc)$. 
    
\end{notn}

\begin{rmk}
    Note that the above is not an extreme abuse of notation in the case of classical operads.
    Let $\Cc$ be a symmetric monoidal (ordinary) category, and equip the nerve $N(\Cc)$ with the symmetric monoidal $\infty$-category structure of \HA{2.1.2.21}.
    Then $\CAlg(N(\Cc))$ can be identified with the nerve of the category of commutative algebra objects on $\Cc$ in the classical sense.
\end{rmk}

\begin{rmk} \label{rmk:O_algebra_in_SMC}
    If $(\Cc, \otimes)$ is a symmetric monoidal category and $\Oc^\otimes$ is any $\infty$-operad, we may view $\Cc$ as an $\Oc$-monoidal category (by pulling back the structure map $\Cc^\otimes \to N(\Fin_*)$ along $\Oc^\otimes \to N(\Fin_*)$.)
    By the universal property of fiber products, we see that $\Alg_{/\Oc}(\Cc^\otimes \times_{N(\Fin_*)} \Oc^\otimes) = \Alg_{\Oc}(\Cc^\otimes)$.
    That is, $\Oc$-algebras in $\Cc$ (viewed as an $\Oc$-monoidal $\infty$-category) are the same as $\Oc$-algebras in $\Cc$ viewed as a symmetric monoidal $\infty$-category.
\end{rmk}

\subsection{Day convolution}

If $\Cc$ is a small $\Oc$-monoidal $\infty$-category, then the presheaf $\infty$-category $\PSh(\Cc) := \Fun(\Cc\op, \Sc)$ can be equipped with a \emph{Day convolution} $\Oc$-monoidal structure: for an $n$-ary operation $f$ in $\Oc$, we define
\begin{equation} \label{eq:day_convolution}
    \big(\otimes_f \{ \Fsc_i \}_{i=1}^n\big)(c) = \colim_{\otimes_f \{ c_i \}_{i=1}^n \to c} \prod_{i=1}^n \Fsc_i(c_i).
\end{equation}
We recall some useful properties of Day convolution as follows -- proofs of these (or very similar statements) can be found in \cite[\S 2.2.6]{ha} and \cite[\S 2 and \S 3]{torii_perfect}.

\begin{prop} \label{prop:day_convolution_properties}
    Let $\Cc$ be a small $\Oc$-monoidal $\infty$-category.
    Then \cref{eq:day_convolution} defines the unique colimit-preserving $\Oc$-monoidal structure on $\PSh(\Cc)$ such that the Yoneda embedding $y_\Cc: \Cc \to \PSh(\Cc)$ is $\Oc$-monoidal.
    Furthermore:
    \begin{enumerate}
        \item If $\Dc$ is a small $\Oc$-monoidal $\infty$-category and $F: \Cc \to \Dc$ is an oplax $\Oc$-monoidal functor, then the pullback functor $F^*: \PSh(\Dc) \to \PSh(\Cc)$ is naturally lax $\Oc$-monoidal.
        \item If $\Ec$ is a cocomplete $\infty$-category with a colimit-preserving $\Oc$-monoidal structure, then left Kan extension along the Yoneda embedding $y_\Cc: \Cc \hookrightarrow \PSh(\Cc)$ induces an equivalence 
        \[
            \Fun^{\Oc,L}(\PSh(\Cc), \Ec) \xrightarrow{\sim}{\to} \Fun^\Oc(\Cc, \Ec).
        \]
        The same is true with $\Fun^\Oc$ replaced by $\Fun^{\Oc,\lax}$ on both sides.
        \item In particular, if $\Dc$ is a small $\Oc$-monoidal $\infty$-category and $F: \Cc \to \Dc$ is an $\Oc$-monoidal functor (resp.\ lax $\Oc$-monoidal functor), then the left Kan extension functor $F_!: \PSh(\Cc) \to \PSh(\Dc)$ is naturally $\Oc$-monoidal (resp.\ lax $\Oc$-monoidal).
        \item There is a natural equivalence $\Alg_{/\Oc}(\PSh(\Cc)) \simeq \Fun^{\Oc,\lax}(\Cc\op, \Sc)$.
    \end{enumerate}
\end{prop}
\begin{proof}
    The fact that Eq. \eqref{eq:day_convolution} defines an $\Oc$-monoidal structure on the Yoneda embedding is well-known (e.g. \cite[Lemma 2.4]{torii_perfect}) and the uniqueness is clear as $\PSh(\Cc)$ is generated by $\Cc$ via colimits.
    
    (1) is \cite[Corollary 2.12]{torii_perfect}.
    The lax monoidal case of (2) is \cite[Proposition 2.7]{torii_perfect}, from which the lax monoidal case of (3) immediately follows (cf. \cite[Definition 2.14]{torii_perfect}).
    The strong monoidal cases of both statements follow by direct computation (using the fact that Day convolution preserves colimits to reduce to looking at representable functors).
    Finally, (4) is \HA{2.2.6.8}.
\end{proof}

In particular, we may use Day convolution to construct $\Oc$-algebra structures on mapping spaces.
The following is essentially contained in the proof of \cite[Corollary 6.8]{nikolaus_stable}:

\begin{cor} \label{cor:maps_form_algebra}
    Let $\Cc$ be a symmetric monoidal $\infty$-category.
    Let $c_1 \in \Alg_\Oc(\Cc\op)$ and let $c_2 \in \CAlg(\Cc)$.
    Then $\Map_\Cc(c_1, c_2)$ is naturally an $\Oc$-algebra in $(\Sc, \times)$,\footnote{As in \cref{notn:infinity_operads}, this $\Oc$-algebra structure is encoded formally by a functor $\Oc^\otimes \to \Sc^\times$.} and this $\Oc$-algebra structure is functorial in $c_1$ and $c_2$.
\end{cor}

\begin{proof}
    We may assume $\Cc$ is small.
    As $y_\Cc$ is symmetric monoidal, \cref{prop:day_convolution_properties} (4) gives
    \[
        y_\Cc(c_2) \in \CAlg(\PSh(\Cc)) \simeq \Fun^{\otimes,\lax}(\Cc\op, \Sc).
    \]
    Thus $y_\Cc(c_2)$ defines a functor $\Alg_\Oc(\Cc\op) \to \Alg_\Oc(\Sc)$ (by composition). 
    In particular, $y_\Cc(c_2) \circ c_1$ gives an $\Oc$-algebra structure to $\Map_\Cc(c_1, c_2) = y_\Cc(c_2)(c_1)$.
    The functoriality follows as $c_1$ also defines a functor $\CAlg(\PSh(\Cc)) \to \Alg_\Oc(\Sc)$ by precomposition.
\end{proof}

\subsection{$\Oc$-monoidal adjunctions and $\Oc$-algebra objects}

We will implicitly use the following fundamental theorem on adjunctions between $\Oc$-monoidal $\infty$-categories:

\begin{prop}[{\cite[Proposition A]{hhln_lax}, \cite[Theorem 1.1]{torii_perfect}}]\label{prop:doctrinal_adjunction}
    Let $\Cc$ and $\Dc$ be $\Oc$-monoidal $\infty$-categories, and let $L: \Cc \rightleftarrows \Dc :R$ be an adjunction of the underlying $\infty$-categories.
    Then the data of an oplax $\Oc$-monoidal structure on $L$ is equivalent to the data of a lax $\Oc$-monoidal structure on $R$.
\end{prop}

With this in mind, we define:

\begin{dfn} \label{dfn:lax_monoidal_adjunction}
    A \emph{lax $\Oc$-monoidal adjunction} $L: \Cc \rightleftarrows \Dc :R$ is an adjunction between $\Oc$-monoidal $\infty$-categories together with an $\Oc$-monoidal structure on $L$.
    An \emph{$\Oc$-monoidal adjunction} is a lax $\Oc$-monoidal adjunction as above in which $L$ is (strong) $\Oc$-monoidal. 
\end{dfn}

\begin{ex}
    If $L: \Cc \rightleftarrows \Dc :R$ is an $\Oc$-monoidal adjunction, the functor $R$ need not be $\Oc$-monoidal.
    For example, if $f: X \to Y$ is any morphism of schemes, then the adjunction $f^*: (\QC(X), \otimes_{\Osc_X}) \rightleftarrows (\QC(Y), \otimes_{\Osc_Y}) : f_*$ is symmetric monoidal, but the functor $f_*$ typically is not symmetric monoidal.
\end{ex}

Lax $\Oc$-monoidal functors induce functors on $\infty$-categories of $\Oc$-algebra objects.
In particular, there is the following standard result:

\begin{prop}[{\HA{7.3.2.13}}] \label{prop:monoidal_adjunction_algebra}
    Let $L: \Cc \rightleftarrows \Dc :R$ be an $\Oc$-monoidal adjunction.
    Then $L$ and $R$ induce an adjunction
    \[
        L: \Alg_{/\Oc}(\Cc) \rightleftarrows \Alg_{/\Oc}(\Dc) :R.
    \]
\end{prop}

We may use \cref{prop:monoidal_adjunction_algebra} to understand the behavior of \cref{cor:maps_form_algebra} with respect to symmetric monoidal functors:

\begin{cor} \label{cor:maps_form_algebra_functorial}
    Let $F: \Cc \to \Dc$ be a symmetric monoidal functor between symmetric monoidal $\infty$-categories.
    Let $c_1 \in \Alg_\Oc(\Cc\op)$ and let $c_2 \in \CAlg(\Cc)$.
    Then the natural map $F_{c_1,c_2}: \Map_\Cc(c_1, c_2) \to \Map_{\Dc}(F(c_1), F(c_2))$ is naturally an $\Oc$-algebra homomorphism.    
\end{cor}

\begin{proof}
    We may write
    \begin{align*}
        \Map_\Dc(F(c_1), F(c_2)) &= y_\Dc\big(F(c_2)\big)\big(F(c_1)\big) \\
        &= F_!\big(y_\Cc(c_2)\big)\big(F(c_1)\big) \\
        &= \big((F^* F_!)(y_\Cc(c_2))\big)(c_1) \textrm{ by the definition of $F^*$}.
    \end{align*}
    The map $F_{c_1,c_2}$ is equivalent to the map $y_\Cc(c_2)(c_1) \to \big((F^* F_!)(y_\Cc(c_2))\big)(c_1)$ induced by the unit map 
    \[
        \eta_{y_\Cc(c_2)}: y_\Cc(c_2) \to (F^* F_!)\big(y_\Cc(c_2)\big).
    \]
    By the functoriality of the construction of \cref{cor:maps_form_algebra}, to show that $F_{c_1,c_2}$ is a homomorphism of $\Oc$-algebras in $\Sc$, it suffices to show that $\eta_{y_\Cc(c_2)}$ upgrades to a homomorphism of commutative algebras in $\PSh(\Cc)$.
    
    The adjunction
    \[
        F_!: \PSh(\Cc) \rightleftarrows \PSh(\Dc) : F^*
    \]
    is symmetric monoidal by \cref{prop:day_convolution_properties}, so by \cref{prop:monoidal_adjunction_algebra} we get an adjunction
    \[
        F_!: \CAlg\big(\PSh(\Cc)\big) \rightleftarrows \CAlg\big(\PSh(\Dc)\big) : F^*.
    \]
    The unit map of this adjunction (at $y_\Cc(c_2) \in \CAlg\big(\PSh(\Cc)\big)$) gives the desired homomorphism of commutative algebras upgrading $\eta_{y_\Cc(c_2)}$.
\end{proof}

In the body of the paper, we will need the following modest strengthening of \cref{prop:monoidal_adjunction_algebra}.

\begin{prop} \label{prop:partial_monoidal_adjunction_algebra}
    Let $L: \Cc \rightleftarrows \Dc :R$ be a lax $\Oc$-monoidal adjunction, and
    suppose that $i: \Cc' \to \Cc$ is the inclusion of a full $\Oc$-monoidal subcategory of $\Cc$ such that the restriction $L \circ i$ is $\Oc$-monoidal.
    Then, for all $c' \in \Alg_{/\Oc}(\Cc')$ and $d \in \Alg_{/\Oc}(\Dc)$, there is a natural equivalence
    \[
        \Map_{\Alg_{/\Oc}(\Dc)}\big(L(i(c')), d\big) \simeq \Map_{\Alg_{/\Oc}(\Cc)}\big(i(c'), R(d)\big).
    \]
\end{prop}

\begin{proof}
    By \cite[Lemma 2.16]{torii_perfect} (passing to a larger universe as necessary to avoid set-theoretic issues) there is a lax $\Oc$-monoidal adjunction
    \[
        L_!: \PSh(\Cc) \rightleftarrows \PSh(\Dc) : R_!.
    \]
    Composing this with the $\Oc$-monoidal adjunction $i_!: \PSh(\Cc') \rightleftarrows \PSh(\Cc) : i^*$, we obtain an adjunction
    \[
        (L \circ i)_!: \PSh(\Cc') \rightleftarrows \PSh(\Dc) : i^* R_!
    \]
    which is $\Oc$-monoidal because $L \circ i$ is.
    By \cref{prop:monoidal_adjunction_algebra}, this induces an adjunction on the corresponding categories of $\Oc$-algebras.

    The Yoneda lemma allows us to embed $\Alg_{/\Oc}(\Cc') \subset \Alg_{/\Oc}(\PSh(\Cc'))$ and $\Alg_{/\Oc}(\Dc) \subset \Alg_{/\Oc}(\PSh(\Dc))$ by \cref{prop:day_convolution_properties}(4).
    Thus, for $c' \in \Alg_{/\Oc}(\Cc')$ and $d \in \Alg_{/\Oc}(\Dc)$, we have
    \begin{align*}
        \Map_{\Alg_{/\Oc}(\Dc)}\big(L(i(c')), d\big) &= \Map_{\Alg_{/\Oc}(\PSh(\Dc))}\big((L \circ i)_!(c'), d\big) \\
        &= \Map_{\Alg_{/\Oc}(\PSh(\Cc'))}\big(c', i^* R_!(d)\big) \\
        &= \Map_{\Alg_{/\Oc}(\PSh(\Cc))}\big(i_!(c'), R_!(d)\big) \textrm{ because $i_!$ is fully faithful} \\
        &= \Map_{\Alg_{/\Oc}(\Cc)}\big(i(c'), R(d)\big). \qedhere
    \end{align*}
\end{proof}

\subsection{Right localization of $\EE_n$-monoidal structures}

We can also ask: given an adjunction of $\infty$-categories in which one category is $\Oc$-monoidal, is there a natural $\Oc$-monoidal structure on the other such that the adjunction upgrades to a (lax) $\Oc$-monoidal adjunction?

When the left adjoint is a localization in the sense of \HTT{5.2.7.2} (i.e.\ when the right adjoint is fully faithful), necessary and sufficient conditions are well-known -- see \HA{2.2.1.9} for the general case as well as \HA{4.1.7.4} for the (symmetric) monoidal case.
In this subsection, we collect some results on the dual problem:
\begin{itemize}
    \item Let $\Dc$ be an $\Oc$-monoidal $\infty$-category, and let $L: \Cc \rightleftarrows \Dc :R$ be an adjunction.
    Assume that $L$ is fully faithful.
    What other conditions guarantee that there exists a natural $\Oc$-monoidal structure on $\Cc$ such that $R: \Dc \to \Cc$ upgrades to a (lax) $\Oc$-monoidal functor?
\end{itemize}

We will restrict to the case $\Oc = \EE_n$ (for some $n \in \{1, \dots, \infty\}$) for notational convenience.

We begin by collecting some well-known results on the behavior of such adjunctions:

\begin{prop}[{\HA{4.1.7.4}\footnote{The claims in \HA{4.1.7.4} are made for the (symmetric) monoidal case only, but the arguments work for $\EE_n$-monoidal structures for any $n$.}}] \label{prop:right_localization_properties}
    Let $L: (\Cc, \otimes_\Cc) \rightleftarrows (\Dc, \otimes_\Dc) :R$ be a lax $\EE_n$-monoidal adjunction.
    Suppose that $L$ is fully faithful and $R$ is strong $\EE_n$-monoidal.
    Then:
    \begin{enumerate}
        \item $\otimes_\Cc$ can be computed via the formula 
        \[
            c_1 \otimes_\Cc c_2 = R(L(c_1)) \otimes_\Cc R(L(c_2)) = R\big(L(c_1) \otimes_\Dc L(c_2)\big).
        \]
        for all $c_1, c_2 \in \Cc$.
        \item The $\EE_n$-monoidal functor $(\Dc, \otimes_\Dc) \to (\Cc, \otimes_\Cc)$ is universal among $\EE_n$-monoidal functors with source $(\Dc, \otimes_\Dc)$ whose underlying functors factor through $\Cc$.
        Symbolically:
        \[
            \Fun^\otimes\big((\Cc, \otimes_\Cc), (\Ec, \otimes_\Ec)\big) = \Fun^\otimes\big((\Dc, \otimes_\Dc), (\Ec, \otimes_\Ec)\big) \times_{\Fun(\Dc, \Ec)} \Fun(\Cc, \Ec).
        \]
        In particular, $\otimes_\Cc$ is unique whenever it exists.
    \end{enumerate}
\end{prop}

In the stable case, $\EE_n$-monoidal adjunctions can be constructed from quotients by (two-sided) thick $\otimes$-ideals.
This is well-known to the experts, though we include the details for completeness.

\begin{dfn} \label{dfn:tensor_ideal}
    Let $(\Cc, \otimes_\Cc)$ be a stably $\EE_n$-monoidal $\infty$-category.
    A \emph{two-sided thick $\otimes$-ideal of $\Cc$} is a stable full subcategory $\Ic \subset \Cc$ such that:
    \begin{itemize}
        \item If $c \oplus c' \in \Ic$, then $c \in \Ic$ and $c' \in \Ic$.
        \item If $i \in \Ic$ and $c \in \Cc$, then $i \otimes_\Cc  c$ and $c \otimes_\Cc i$ are both in $\Ic$.
    \end{itemize}
\end{dfn}

\begin{lem}[{Stable, dual version of \HA{2.2.1.9}}] \label{lem:transport_adjunction}
    Let $(\Dc, \otimes_\Dc)$ be a stably $\EE_n$-monoidal $\infty$-category, and let $\Cc$ be a stable $\infty$-category.
    Let $L: \Cc \rightleftarrows \Dc :R$ be an adjunction.
    Suppose:
    \begin{enumerate}
        \item $L$ is fully faithful, and
        \item The full subcategory 
        \[
            \ker R := \bset*{d \in \ob \Dc}{R(d) \simeq 0}
        \]
        is a two-sided thick $\otimes$-ideal of $\Dc$.{}
    \end{enumerate}
    Then there exists a unique $\EE_n$-monoidal structure on $\Cc$ such that $L: \Cc \rightleftarrows \Dc :R$ is $\EE_n$-monoidal.
\end{lem}

\begin{proof}
    By the dual of \HA{2.2.1.9}, it suffices to show that, if $f$ and $g$ are morphisms in $\Dc$ such that $R(f)$ and $R(g)$ are equivalences in $\Cc$, then $R(f \otimes_\Dc g)$ is also an equivalence in $\Dc$.
    Factoring $f \otimes_\Dc g$ as $(f \otimes_\Dc \id_{d_1}) \circ (\id_{d_2} \otimes_\Dc g)$, we see that it suffices to show that, if $f$ is a morphism in $\Dc$ such that $R(f)$ is an equivalence in $\Cc$, then $R(f \otimes_\Dc \id_d)$ is an equivalence for all $d \in \Dc$.\footnote{Of course, when $n = 1$, we must also show the same result for $R(\id_d \otimes_\Dc f)$, but the same argument will apply.}
    
    Since $R(f)$ is an equivalence, we have $R(\fib f) = \fib R(f) \simeq 0$, so $\fib f \in \ker R$.
    Because $\ker R$ is a $\otimes_\Dc$-ideal, we have $\fib (f \otimes_\Dc \id_d) = \fib f \otimes_\Dc d \in \ker R$.
    Thus $\fib R(f \otimes \id_d) = 0$ and $R(f \otimes_\Dc \id_d)$ is an equivalence.
\end{proof}

\begin{rmk}
    In the context of \cref{lem:transport_adjunction}, the category $\ker R$ is automatically stable and closed under taking direct summands.
    Thus it suffices to check that if $i \in \ker R$ and $d \in \Dc$, then $i \otimes_\Dc d$ and $d \otimes_\Dc i$ are both in $\Dc$.
\end{rmk}

When $\Cc$ and $\Dc$ are compactly generated, we may provide a criterion for the existence of a corresponding $\EE_n$-monoidal structure on $\Cc$ that requires us to understand only the behavior of $L$ (rather than $R$).
We use this result in the body of the paper (see \cref{prop:ec_exists}).

\begin{prop} \label{prop:transport_convolution}
    Fix $n \in \{ 1, \dots, \infty\}$.
    Suppose $(\Dc, \otimes_\Dc)$ is a compactly generated, stable, and presentably $\EE_n$-monoidal $\infty$-category.\footnote{We do \emph{not} assume that the monoidal structure $\otimes_\Dc$ preserves compact objects.}
    Let $G: \Dc \otimes \Dc \to \Dc$ be the unique colimit-preserving functor such that $- \otimes_\Dc - = G(- \boxtimes -)$.\footnote{$G$ exists by the universal property of the Lurie tensor product, since $- \otimes_\Dsc -: \Dsc \times \Dsc \to \Dsc$ preserves colimits in each variable separately.}
    Let $\Cc$ be a compactly generated stable $\infty$-category, and let $L: \Cc \rightleftarrows \Dc :R$ be an adjunction. 
    Suppose:
    \begin{enumerate}
        \item $L$ is fully faithful,
        \item $L$ preserves compact objects,
        \item $G$ has a left adjoint $F: \Dc \to \Dc \otimes \Dc$, and
        \item The image of $F \circ L$ lies in the full subcategory of $\Dc \otimes \Dc$ compactly generated by
        \[
            \big\{ L(c_1) \boxtimes L(c_2) \,\big|\, c_1, c_2 \in \Cc^\omega \big\}.
        \]
    \end{enumerate}
    Then there exists a unique $\EE_n$-monoidal structure $\otimes_\Cc$ on $\Cc$ such that $R: (\Dc, \otimes_\Dc) \to (\Cc, \otimes_\Cc)$ is $\EE_n$-monoidal.
\end{prop}

\begin{proof}
    By \cref{lem:transport_adjunction}, it suffices to check that, if $d_1 \in \ker R$ and $d_2 \in \Dc$, then $d_1 \otimes_\Dc d_2 \in \ker R$.\footnote{When $n = 1$, we must show $d_2 \otimes_\Dc d_1 \in \ker R$, but the same argument works.}
    By the Yoneda lemma, it suffices to show that $\Hom_\Cc\big(c, R(d_1 \otimes_\Dc d_2)\big) = 0$ for all $c \in \Cc$.
    Observe that
    \[
        \Hom_\Cc\big(c, R(d_1 \otimes_\Dc d_2)\big) = \Hom_\Cc\big(c, R(G(d_1 \boxtimes d_2))\big) = \Hom_{\Dc \otimes \Dc}\big(F(L(c)), d_1 \boxtimes d_2\big).
    \]
    By hypothesis, $F(L(c))$ can be written as a colimit of terms $L(c_1) \boxtimes L(c_2)$ for $c_1, c_2 \in \Cc^{\omega}$.
    Thus it suffices to show that $\Hom_{\Dc \otimes \Dc}\big(L(c_1) \boxtimes L(c_2), d_1 \boxtimes d_2\big) = 0$ for all $c_1, c_2 \in \Cc^{\omega}$.
    But this is just a direct computation using the K\"unneth formula:
    \begin{align*}
        \Hom_{\Dc \otimes \Dc}\big(L(c_1) \boxtimes L(c_2), d_1 \boxtimes d_2\big) & = \Hom_\Dc\big(L(c_1), d_1\big) \otimes \Hom_\Dc\big(L(c_2), d_2\big) \\
        &= \Hom_\Cc\big(c_1, R(d_1)\big) \otimes \Hom_\Cc\big(L(c_2), d_2\big) \\
        &= 0 \otimes \Hom_\Cc\big(L(c_2), d_2\big) = 0. \qedhere
    \end{align*}
\end{proof}

\printbibliography

\end{document}